\documentclass{article} % [12pt]
\usepackage{latexsym,amsmath,amscd,amssymb,graphics}
\usepackage[backend=bibtex,sorting=none]{biblatex}
\usepackage{enumerate}
\usepackage[]{mcode}
\usepackage[toc,page]{appendix}
\usepackage[a4paper,margin=.75in]{geometry}
\usepackage{array}
\usepackage[section]{placeins}
\usepackage{afterpage}
\usepackage{verbatim}

\usepackage{graphicx}
\usepackage{subfig}
\usepackage{framed}
\usepackage[colorlinks]{hyperref}
\usepackage{url}
\usepackage{mathrsfs}
\usepackage{epstopdf}

\usepackage[all]{xy}
\usepackage{indentfirst}
\usepackage[mathscr]{euscript}
\usepackage{algorithm}
\usepackage{algpseudocode}
\usepackage{authblk}

\makeatletter
\renewcommand*{\@fnsymbol}[1]{\ensuremath{\ifcase#1\or \dagger\or *\or \ddagger\or
		\mathsection\or \mathparagraph\or \|\or **\or \dagger\dagger
		\or \ddagger\ddagger \else\@ctrerr\fi}}
\makeatother

\makeatletter

\@addtoreset{figure}{section}
\def\thefigure{\thesection.\@arabic\c@figure}
\def\fps@figure{h, t}
\@addtoreset{table}{bsection}
\def\thetable{\thesection.\@arabic\c@table}
\def\fps@table{h, t}
\@addtoreset{equation}{section}

\makeatother


\newcommand{\unaryminus}{\scalebox{0.75}[1.0]{\( - \)}}

\newtheorem{theorem}{Theorem}

\newtheorem{remark}[theorem]{Remark}

\numberwithin{theorem}{subsection}

\def\bea{\begin{eqnarray}}
\def\eea{\end{eqnarray}}
\def\ba{\begin{array}}
\def\ea{\end{array}}

\def\bomega{\boldsymbol{\omega}}
\def\bOm{\boldsymbol{\Omega}}

\def\brho{\boldsymbol{\rho}}

\def\diag{\mathrm{ \textbf{diag}}}

\def\bx{{\boldsymbol {x} }}

\let\<\langle
\let \>\rangle
\newcommand{\rem}[1]{}

\newcommand{\de}{\delta}

\newcommand{\bu}{\boldsymbol{u}}
\newcommand{\bPsi}{\boldsymbol{\Psi}}

\newcommand{\bv}{\boldsymbol{v}}
\newcommand{\bz}{\boldsymbol{z}}

\newcommand{\bGam}{\boldsymbol{\Gamma}}
\newcommand{\bGamma}{\boldsymbol{\Gamma}}
\newcommand{\bzeta}{\boldsymbol{\zeta}}
\newcommand{\bom}{\boldsymbol{\omega}}

\newcommand{\bsigma}{\boldsymbol{\sigma}}
\newcommand{\bSigma}{\boldsymbol{\Sigma}}

\newcommand{\bxi}{\boldsymbol{\xi}}

\newcommand{\bY}{\mathbf{Y}}

\newcommand{\bkappa}{\boldsymbol{\kappa}}

\newcommand{\bchi}{\boldsymbol{\chi}} 
\newcommand{\btheta}{\boldsymbol{\theta}}
\newcommand{\inertia}{\mathbb{I}}

\newcommand{\pp}[2]{\frac{\partial #1}{\partial #2}}
\newcommand{\dd}[2]{\frac{\mathrm{d} #1}{\mathrm{d} #2}}
\newcommand{\dede}[2]{\frac{\delta #1}{\delta #2}}

\newcommand{\dt}{\mathrm{d}t}

\newcommand{\mso}{\mathfrak{so}}

\newcommand{\dprime}{\prime \prime}

\newcommand{\todo}[1]{\vspace{5 mm}\par \noindent
\framebox{\begin{minipage}[c]{0.95 \textwidth}
\tt #1 \end{minipage}}\vspace{5 mm}\par}
\newcommand{\revision}[2]{#2}

\usepackage[normalem]{ulem}

\graphicspath{ {./figures_gen/} {./} {./figures_disk/} {./figures_ball/} }

\title{On the Dynamics of a Rolling Ball Actuated by Internal Point Masses}
\author[1]{Vakhtang Putkaradze\thanks{Email address: \texttt{putkarad@ualberta.ca}}}
\author[2]{Stuart Rogers\thanks{Email address: \texttt{srogers@umn.edu}}}
\affil[1]{Department of Mathematical and Statistical Sciences, Faculty of Science, University of Alberta, CAB 632, Edmonton, AB T6G 2G1, Canada}
\affil[2]{Institute for Mathematics and its Applications, College of Science and Engineering, University of Minnesota, 207 Church Street SE, 306 Lind Hall, Minneapolis, MN 55455, USA}
\date{\today}

\providecommand{\keywords}[1]{\textbf{\textit{Keywords:}} #1}
\begin{document}

\maketitle

\abstract{\noindent 
The motion of a rolling ball actuated by internal point masses that move inside the ball's frame of reference is considered. The equations of motion are derived by applying  Euler-Poincar\'e's symmetry reduction method in concert  with Lagrange-d'Alembert's principle, which accounts   for the presence of the nonholonomic rolling constraint. As a particular example, we consider the case when the masses move along internal rails,  or trajectories, of arbitrary shape and fixed within the ball's frame of reference. Our system of equations can treat most possible  methods of actuating the rolling ball with internal moving masses encountered in the literature, such as circular motion of the masses mimicking   swinging pendula or straight line motion of the masses mimicking magnets  sliding inside linear tubes embedded within a solenoid. Moreover, our method can  model arbitrary  rail shapes and  an arbitrary number of rails such as several ellipses  and/or figure eights, which may be important for future designs of rolling ball robots.  For further analytical study, we also reduce the system to a single differential equation when the motion is  planar, that is, considering  the motion of the rolling disk actuated by internal  point masses, in which case we show that the results  obtained from the variational derivation coincide with those obtained from Newton's second law. Finally, the equations of motion are solved numerically, illustrating a wealth of complex behaviors exhibited by the system's dynamics. Our results are  relevant to the dynamics of nonholonomic systems containing internal degrees of freedom and to further studies of control of such systems actuated by internal masses. 
\\
\\
\keywords{symmetry reduction, nonholonomic mechanics, rolling balls}
\tableofcontents 

\section{Introduction} \label{sec_introduction}

\subsection{Motivation and Methodology}
The first six films in the famous \textit{Star Wars} space saga starred the sidekick robot R2-D2, which locomoted via a three-wheeled tripod. However, the seventh and eighth films in that saga, \textit{The Force Awakens} and \textit{The Last Jedi}, star new, next-generation, sidekick robots called BB-8 and BB-9E. BB-8, depicted in Figure~\ref{fig_bb8_1}, and BB-8's evil nemesis BB-9E each locomote via a single rolling ball. To cash in on these new \textit{Star Wars} fan favorites, the toy company Sphero sells working toy models of BB-8 and BB-9E. But rolling ball robots are not just gimmicks used by the entertainment and toy industries. The defense, security, energy, and agricultural industries are also interested in exploiting sensor-equipped rolling ball robots, such as Rosphere shown in Figure~\ref{fig_rosphere}, for such tasks as surveillance and environmental monitoring. The goal of this paper is to study some mechanisms  for actuating the motion of rolling ball robots like BB-8, BB-9E, and Rosphere. This paper deals exclusively with the derivation and analysis of the uncontrolled equations of motion. Another paper \cite{putkaradze2017optimal} by the authors  investigates the optimal control of  rolling ball robots that are able to locomote over a prescribed trajectory, avoid obstacles, and/or perform some other maneuver by minimizing a prescribed performance index (a.k.a. cost functional). The derivation of the  uncontrolled dynamics is highly nontrivial and, as far as we know, has not been done before in the generality we present here. 
\begin{figure}[h] 
	\centering
	\subfloat[Sphero's toy incarnation of BB-8, one of \textit{Star Wars'} next-generation rolling ball robots  \cite{BB8_2}.]{\includegraphics[scale=.1025]{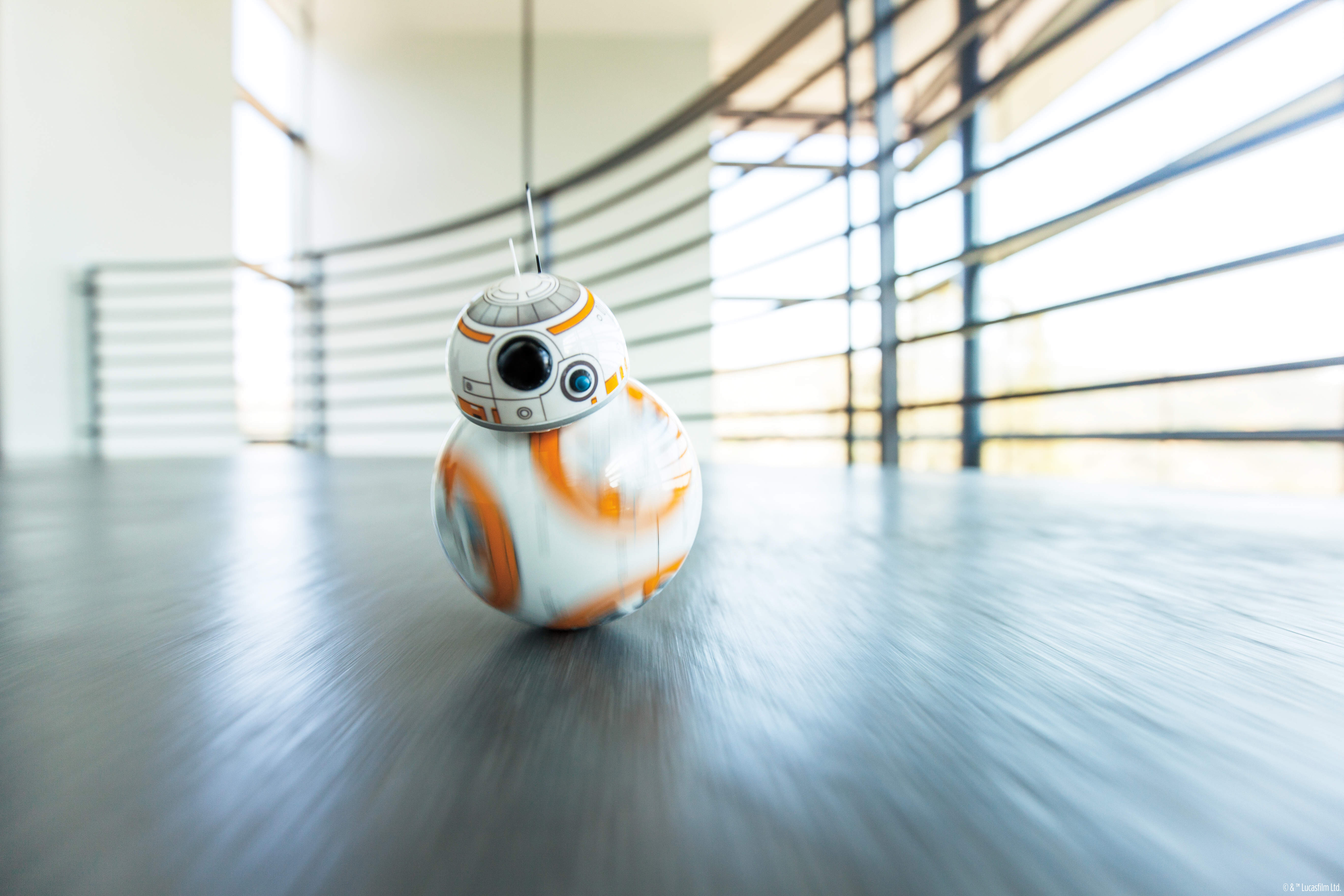}\label{fig_bb8_1}}
	\hspace{5mm}
	\subfloat[Rosphere can be used in agriculture for monitoring crops, \copyright\ 2013 Emerald \cite{hernandez2013moisture}.]{\includegraphics[scale=.151]{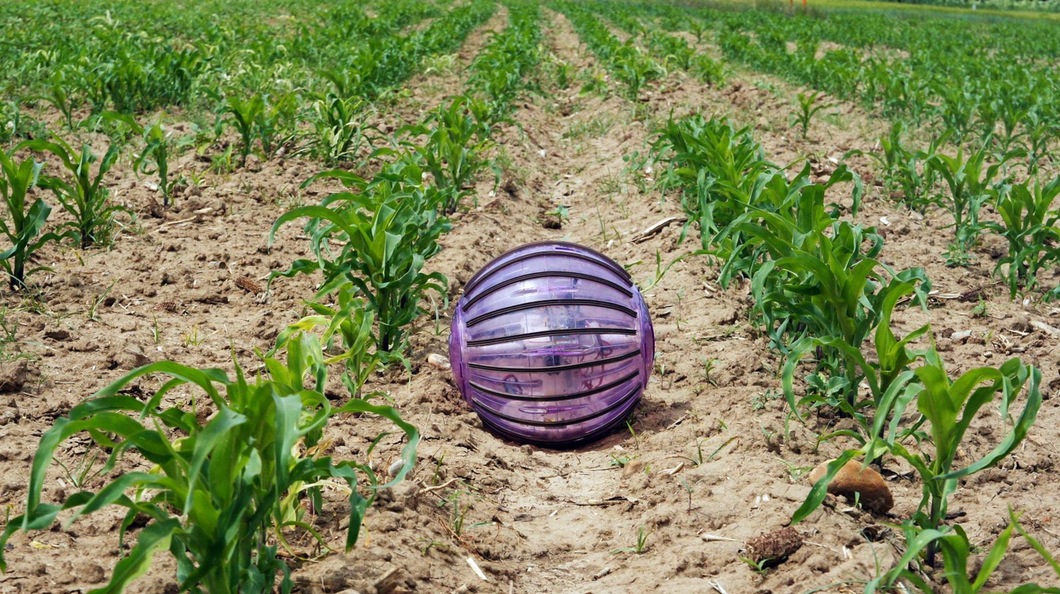}\label{fig_rosphere}}
	\caption{Examples of real rolling ball robots.}
\end{figure}

Before optimal control can be applied to the rolling ball, its ordinary differential equations of motion must be derived first; henceforth, the ordinary differential equations of motion of  the rolling ball will be referred to as the equations of motion or the uncontrolled equations of motion to distinguish them from the controlled equations of motion. To derive the uncontrolled equations of motion for the rolling ball, methods from nonholonomic mechanics must be utilized since the rolling ball is subject to a nonholonomic (as opposed to a holonomic) constraint and therefore is an example of a nonholonomic mechanical system. A constraint affine in velocity is called nonholonomic if it is ideal (i.e. virtual work on the constraint vanishes) and cannot be re-expressed as a position constraint; if the constraint can be expressed soley as a function of position, then it is said to holonomic. 
\rem{\todo{VP: Why don't you want to define the nonholonomic system above? The definition above (old one) was \sout{crossed out}. Maybe as a footnote? Technically, what we wrote there before was not quite correct: the velocity in the constraint must be involved in a non-integrable way, so the constraint cannot be written in terms of coordinates. Also, the realization of constraint must be through the LdA principle, so the constraint is not vakonomic. I changed the sentence above. \\ SMR: I made some revisions above.  
\\
VP: OK. One last thing: we have two definitions above: one precise (second) and one imprecise (first). Should we just keep the second one and drop the first one? \\ SMR: How is the above? }} The uncontrolled equations of motion governing a nonholonomic mechanical system are given by Lagrange-d'Alembert's principle, a somewhat nonintuitive method in mechanics developed by Jean d'Alembert in the 18th century. In addition, Euler-Poincar\'e's method \cite{poincare1901forme}, \revision{R2Q1}{first} published by Henri Poincar\'e in 1901 \revision{R2Q1}{and independently replicated in greater generality by Georg Hamel \cite{hamel1904,borisov2016historical} in 1904}, provides a more efficient derivation of the equations of motion of the rolling ball compared to the standard Hamilton's principle by using symmetry arguments to reduce the degrees of freedom in the dynamics. 

\subsection{Background}

Consider a ball rolling without slipping on a horizontal surface in the presence of a uniform gravitational field. Figure~\ref{fig:simple_rolling_ball} shows a ball of radius $r$ rolling without slipping on a horizontal surface in the presence of a uniform gravitational field of magnitude $g$.

\begin{figure}[h]
	\centering
	\includegraphics[width=0.5\linewidth]{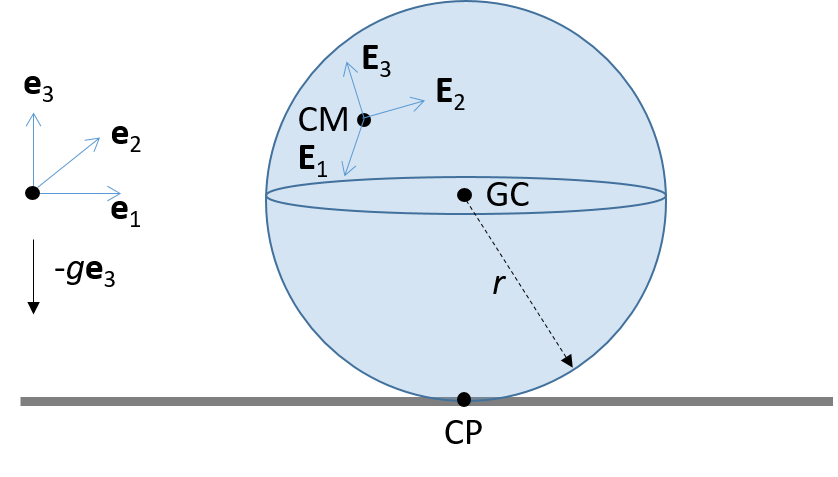}
	\caption{A ball of radius $r$ rolls without slipping on a horizontal surface in the presence of a uniform gravitational field of magnitude $g$. The ball's geometric center, center of mass, and contact point with the horizontal surface are denoted by GC, CM, and CP, respectively.}
	\label{fig:simple_rolling_ball}
\end{figure}

There are several terminologies in the literature to describe a ball rolling without slipping on a horizontal surface in the presence of a uniform gravitational field, depending on its mass distribution and the location of its center of mass. A Chaplygin sphere is a ball with an inhomogeneous mass distribution, but with its center of mass located at the ball's geometric center \cite{shen2008controllability}. A Chaplygin top is a ball with an inhomogeneous mass distribution, but with its center of mass not located at the ball's geometric center \cite{shen2008controllability}. Reference \cite{Ho2011_pII} does not distinguish between these two cases, calling a Chaplygin ball a ball with an inhomogeneous mass distribution, regardless of the location of its center of mass; as a special case of a Chaplygin ball, \cite{Ho2011_pII} calls a Chaplygin concentric sphere a ball with an inhomogeneous mass distribution with its center of mass coinciding with the ball's geometric center. Thus, the Chaplygin concentric sphere (used by \cite{Ho2011_pII}) and the Chaplygin sphere (used by \cite{shen2008controllability}) are different terms for the same mechanical system. Note that a ball with a homogeneous mass distribution (in a uniform gravitational field) necessarily has its center of mass at the ball's geometric center, and is therefore not very interesting. In this paper, these terminologies are not used, rather the mechanical system is referred to simply as a ball or a rolling ball, regardless of its mass distribution (homogeneous vs inhomogeneous) and regardless of the location of its center of mass (at the ball's geometric center vs not at the ball's geometric center). 

In this paper, the motion of the rolling ball is investigated assuming both static and dynamic internal structure. The dynamics of the rolling ball with static internal structure was first solved analytically by Chaplygin for the cylindrically symmetric rolling ball, i.e. a ball such that the line joining the ball's center of mass and geometric center forms an axis of symmetry, in 1897 \cite{chaplygin2002motion} and for the Chaplygin sphere in 1903 \cite{chaplygin2002ball}, though dynamical properties of the cylindrically symmetric rolling ball were previously investigated by Routh \cite{routh1884advanced} and Jellet \cite{jellett1872treatise}. \revision{R2Q3}{\rem{Recent studies of the dynamics of the Chaplygin top and sphere appear in \cite{borisov2016spiral,borisov2014reversal} and \cite{borisov2013problem}, respectively.} More recently, \cite{borisov2013problem} provides a detailed analysis of the trajectory of the Chaplygin sphere's contact point, and it has been shown that the dynamics of the Chaplygin top exhibit a strange attractor \cite{borisov2016spiral} and the phenomenon of reversal \cite{borisov2014reversal}. }  The dynamics of the rolling ball with dynamic internal structure is \rem{still} \revision{R2Q3}{also} an active topic in the nonholonomic mechanics literature \cite{das2001design,mojabi2002introducing,shen2008controllability,borisov2012control,bolotin2012problem,gajbhiye2016geometric,GaBa:RNC3457,kilin2015spherical,burkhardt2016reduced}. \rem{For example, the recent paper \cite{GaBa:RNC3457} explores \emph{local} controllability of a rolling ball with dynamic internal structure, in which an internal pendulum and yoke serve as the control mechanisms.}  

Many methods have been proposed (and some realized) to actuate a rolling ball, such as illustrated in Figure~\ref{fig_ball_actuation}. References \cite{borisov2012control,bolotin2012problem,gajbhiye2016geometric} actuate the rolling ball by internal rotors such as shown in Figure~\ref{fig_ball_rotors}, while \cite{burkhardt2014energy,davoodi2014moball,asama2015design,davoodi2015moball,bowkett2016combined,burkhardt2016reduced} actuate the rolling ball via $6$ internal magnets, each of which slides inside its own linear, solenoidal tube, i.e. a straight tube embedded within a solenoid that generates a magnetic field along the tube's longitudinal axis as illustrated in Figure~\ref{fig_ball_magnets}. 
References  \cite{das2001design,mojabi2002introducing} study the locomotion and trajectory-tracking of a ball with masses moving along straight rails inside the ball, as well as practical realizations of such a device. In particular, \cite{mojabi2002introducing} actuates the rolling ball by internal masses which reciprocate along spokes. Reference \cite{shen2008controllability} actuates the rolling ball by a combination of internal rotors and sliders, \cite{kilin2015spherical} actuates the rolling ball by an internal gyroscopic pendulum as shown in Figure~\ref{fig_ball_gyro_pend}, \cite{bolotin2013motion,pivovarova2014stability,ivanova2015dynamics,ivanova2018controlled} actuate the rolling ball by an internal spherical pendulum as shown in Figure~\ref{fig_ball_sph_pend}, and \cite{GaBa:RNC3457} actuates the rolling ball by an internal pendulum and yoke. This paper considers a rolling ball actuated by internal point masses that move along arbitrarily-shaped rails fixed within the ball, such as depicted in Figure~\ref{fig_intro_bsim2_control_masses_rails}. Actuating the rolling ball by moving internal point masses along general rails has not been considered yet in the literature; references \cite{burkhardt2016reduced,das2001design,mojabi2002introducing,shen2008controllability} actuate the rolling ball by moving internal masses with inertias along linear trajectories (e.g. spokes or hollow tubes) in the ball's frame. The very recent work \cite{ilin2017dynamics}  derives and simulates the dynamics of a   simplified model of a beaver ball. The results in that work (obtained independently at around the time of submission of this paper) analyze a ball  actuated  by a point mass moving  with constant angular velocity  along a  trajectory (taken to be a circle)  fixed inside the  ball. 

 This paper investigates the dynamics of the rolling ball actuated by  the general motion of internal  point masses  using the variational (Lagrange-d'Alembert's) principle of nonholonomic mechanics.  To contrast with  previous works, we have assumed maximum generality for the motion of the internal  point masses, which includes, as particular cases, all previous ways of actuating a rolling ball.  Our methods are also applicable to more complex ways of actuating the  ball,  for example  when the rails are moving relative to the ball  or when the  ball  is  driven by a double pendulum; however, we do not consider these cases here because of  their  algebraic complexity. Our paper is the first step towards the derivation of the general principle for the optimal control of such rolling balls, dealing exclusively with the dynamics.  A separate paper \cite{putkaradze2017optimal} uses the results presented here to derive optimal control techniques for such robots. 
\rem{ %%%%BEGIN REM 
In a comprehensive review of nonholonomic optimal control, \cite{Bloch2003} briefly discusses the controllability and optimal control (in the sense of Pontryagin's minimum principle) of a rolling ball, where an external control force pushes the ball's geometric center. Reference~\cite{shen2008controllability} investigates the controllability of a rolling ball actuated by internal rotors and sliders, while \cite{GaBa:RNC3457} investigates the controllability of a rolling ball actuated by an internal pendulum and yoke. Reference~\cite{borisov2012control} investigates the controllability and trajectory-tracking control of a rolling ball actuated by internal rotors, \cite{bolotin2012problem} investigates the optimal control (in the sense of Pontryagin's minimum principle) of a rolling ball actuated by internal rotors for some very specific cost functionals, and \cite{gajbhiye2016geometric} investigates the orientation and contact point trajectory-tracking control of a rolling ball actuated by internal rotors by constructing feedback control laws. Reference \cite{mojabi2002introducing} investigates the trajectory-tracking control of a rolling ball, actuated by internal masses which reciprocate along spokes, by optimizing the ball's energy expense during the motion by discretizing the trajectory and finding a global minimum of a certain function, i.e. by using the direct method. Reference~\cite{kilin2015spherical} investigates the trajectory-tracking control of a rolling ball actuated by an internal gyroscopic pendulum by piecing together special gaits, each of whose motion is determined analytically. To contrast with previous work, this paper investigates the optimal control (in the sense of Pontryagin's minimum principle) of a rolling ball actuated by internal point masses for very general cost functionals.   
} %%%END REM 
\rem{
\begin{figure}[h] 
	\centering
	\subfloat[A ball actuated by $3$ rotors, studied in \cite{borisov2012control,bolotin2012problem,gajbhiye2016geometric}, \copyright\ 2016 IFAC   \cite{gajbhiye2016geometric}.]{\includegraphics[width=0.45\textwidth]{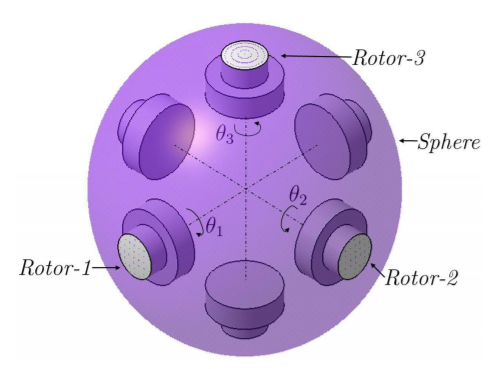}\label{fig_ball_rotors}}
	\hspace{5mm}
	\subfloat[A ball actuated by $3$ point masses, each on its own circular rail, studied in this paper.]{
		\includegraphics[width=0.45\textwidth]{bsim2_control_masses_rails_png}\label{fig_intro_bsim2_control_masses_rails}}
	\\
	\subfloat[A ball actuated by a gyroscopic pendulum, studied in \cite{kilin2015spherical}.]{
		\includegraphics[width=0.45\textwidth]{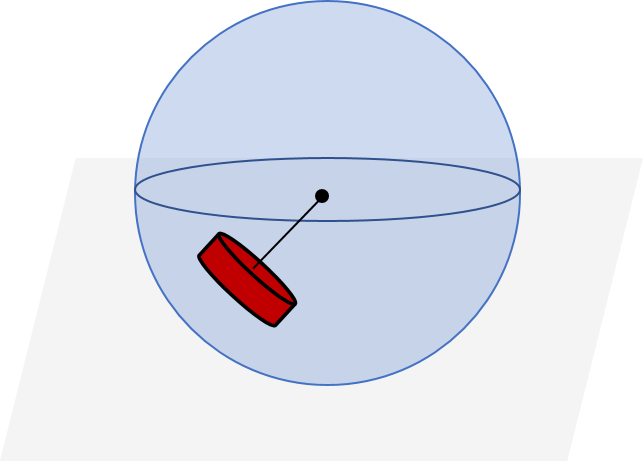}\label{fig_ball_gyro_pend}}
	\hspace{5mm}
	\subfloat[Sphero has $4$ wheels wedged inside the spherical shell, but only the lower $2$ are spun by the motor \cite{Sphero_cut}.]{
		\includegraphics[width=0.45\textwidth]{sphero_cutaway}\label{fig_ball_sphero_internal}}	
	\caption{Different methods to actuate a rolling ball.} \label{fig_ball_actuation}
\end{figure} }

\begin{figure}[h] 
	\centering
	\subfloat[A ball actuated by $3$ rotors, studied in \cite{borisov2012control,bolotin2012problem,gajbhiye2016geometric}, \copyright\ 2016 IFAC   \cite{gajbhiye2016geometric}.]{\includegraphics[width=0.3\textwidth]{ball_rotors}\label{fig_ball_rotors}}
	\hspace{5mm}
	\subfloat[A ball actuated by $6$ magnets, each in its own linear, solenoidal tube, studied in \cite{burkhardt2014energy,davoodi2014moball,asama2015design,davoodi2015moball,bowkett2016combined,burkhardt2016reduced}.]{
		\includegraphics[width=0.3\textwidth]{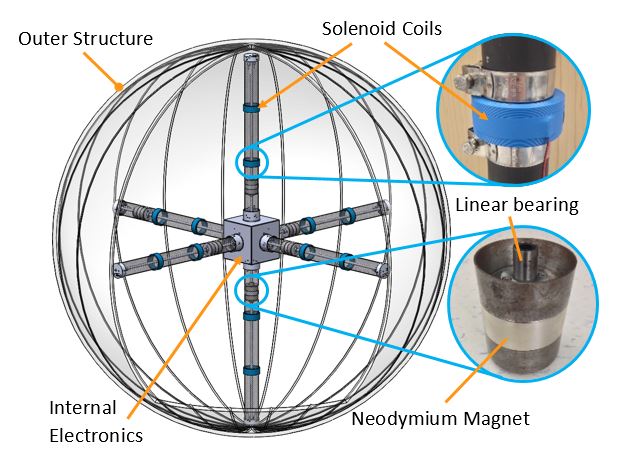}\label{fig_ball_magnets}}
	\hspace{5mm}
	\subfloat[A ball actuated by $3$ point masses, each on its own circular rail, studied in this paper.]{
		\includegraphics[width=0.3\textwidth]{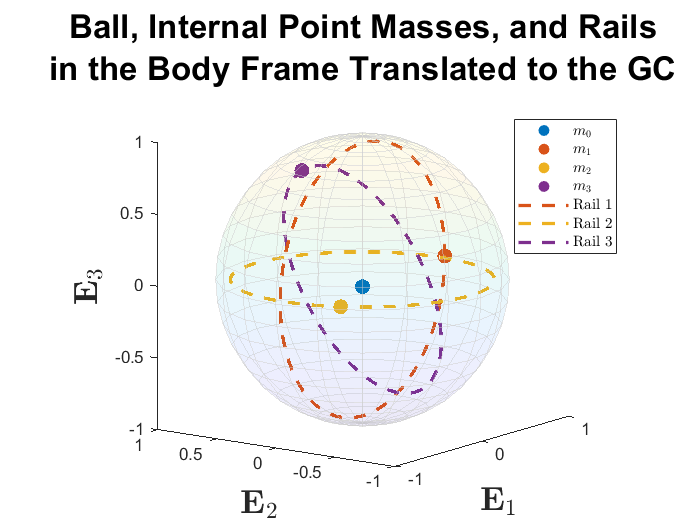}\label{fig_intro_bsim2_control_masses_rails}}
	\hspace{5mm}
	\\
	\subfloat[A ball actuated by a gyroscopic pendulum, studied in \cite{kilin2015spherical}.]{
		\includegraphics[width=0.3\textwidth]{ball_gyro_pend_custom}\label{fig_ball_gyro_pend}}
	\hspace{5mm}
		\subfloat[A ball actuated by a spherical pendulum, studied in \cite{bolotin2013motion}.]{
		\includegraphics[width=0.3\textwidth]{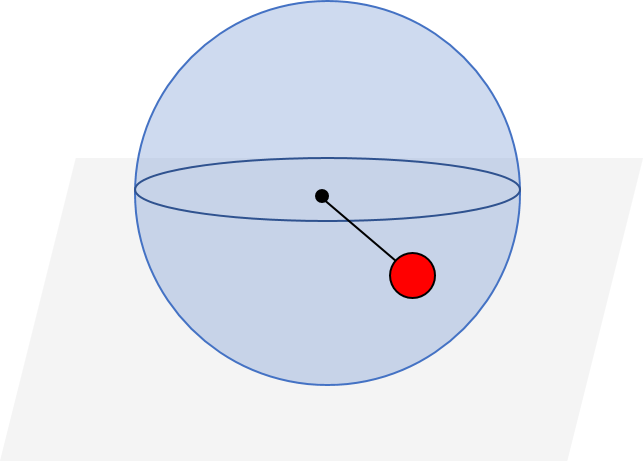}\label{fig_ball_sph_pend}}
	\hspace{5mm}
	\subfloat[Sphero has $4$ wheels wedged inside the spherical shell, but only the lower $2$ are spun by the motor \cite{Sphero_cut}.]{
		\includegraphics[width=0.3\textwidth]{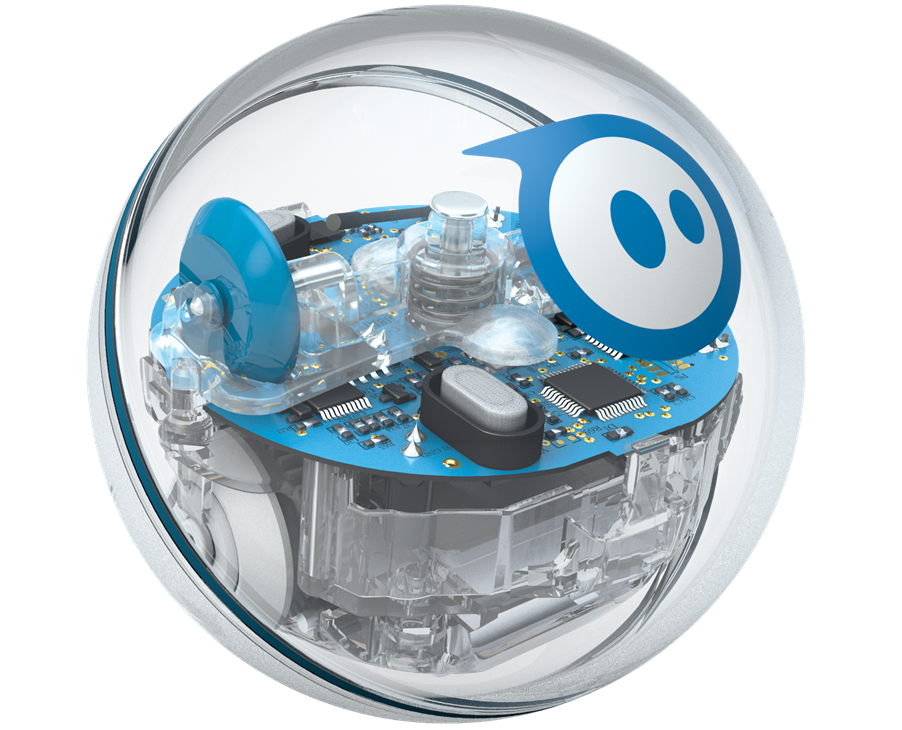}\label{fig_ball_sphero_internal}}	
	\caption{Different methods to actuate a rolling ball.} \label{fig_ball_actuation}
\end{figure}

\rem{\todo{VP: About this figure: the Sphero picture is much bigger than everything else. Can you try to match its size? Maybe using the height measurement instead of width in the pictures? \\ SMR: The first 3 figures have width .21textwidth, while the last figure now has width .16textwidth. Is the last figure good enough now? The height of the last figure better matches the other 3 figures, but the caption is not as wide as the other captions.  
\\ 
VP: It is good. The caption is a bit ugly, but OK.   } }
The paper is organized as follows. 
Section~\ref{sec_ball_mech_system} discusses the specific type of rolling ball considered, presents natural questions about this rolling ball that motivate this paper, and defines coordinates systems and notation used to describe this rolling ball. By applying Euler-Poincar\'e's method and Lagrange-d'Alembert's principle, Section~\ref{sec_ball_uncontrolled} derives the equations of motion for the rolling ball. Finally, numerical simulations of the ball's dynamics are presented in Section~\ref{sec_disk_numerical} for the case of the rolling disk (i.e. 2-d motion) and in Section~\ref{sec_ball_numerical} for general 3-d motion. In addition, Appendix~\ref{app_background} reviews Euler-Poincar\'e's method and nonholonomic mechanics since they are used to derive the equations of motion in Section~\ref{sec_ball_uncontrolled}. Also, Appendix~\ref{app_quaternions} reviews quaternions, which are utilized to formulate the equations of motion used to simulate the ball's dynamics.

\section{Mechanical System, Coordinate Systems, and Notation} \label{sec_ball_mech_system} 
Consider a rigid ball of radius $r$ containing some static internal structure as well as $n \in \mathbb{N}^0$ point masses, where $\mathbb{N}^0$ denotes the set of nonnegative integers. This ball rolls without slipping on a horizontal  surface in the presence of a uniform gravitational field. For $1 \le i \le n$, the $i^\mathrm{th}$ point mass may move within the ball along a trajectory $\bxi_i$, expressed with respect to the ball's frame of reference, as illustrated in Figure~\ref{fig:detailed_rolling_ball}. The trajectory $\bxi_i$ may be constrained in some way, such as being required to move along a 1-d rail (like a circular hoop), across a 2-d surface (like a sphere), or within a 3-d region (like a ball) fixed within the ball. The ball with its static internal structure has mass $m_0$ and the $i^\mathrm{th}$ point mass has mass $m_i$ for $1 \le i \le n$. Let $M = \sum_{i=0}^n m_i$ denote the mass of the total system. The total mechanical system consisting of the ball with its static internal structure and the $n$ point masses is referred to as the ball or the rolling ball, the ball with its static internal structure but without the $n$  point masses may also be referred to as $m_0$, and the $i^\mathrm{th}$  point mass may also be referred to as $m_i$ for $1 \le i \le n$. 

It is natural to ask the following questions for this mechanical system:
\begin{enumerate}
	\item How does the ball move if the $n$   masses are held fixed in place? 
	\item Given some prescribed motion of the $n$  masses, how does the ball move along the horizontal surface? 
	\item Suppose that it is desired to move the ball in a prescribed manner, such as moving the ball's geometric center along a prescribed trajectory parallel to the horizontal surface or performing obstacle avoidance. How might the $n$   masses be moved to realize such a motion? Figure~\ref{fig:ball_rolling_along_traj} illustrates this problem for $2$ masses.
\end{enumerate}
The remainder of this paper aims to answer questions 1 and 2. The answer to the $2^\mathrm{nd}$ question also answers the $1^\mathrm{st}$, by insisting that the prescribed motion for each  point  mass be that of holding it fixed within the ball. The $3^\mathrm{rd}$ question is the inverse of the $2^\mathrm{nd}$. Chaplygin answered the $1^\mathrm{st}$ question analytically for two special cases in his seminal 1897 and 1903 papers \cite{chaplygin2002motion,chaplygin2002ball}. In the general case, no analytical solution can be found for the $1^\mathrm{st}$ question, although the equations of motion are readily  integrated  numerically. As far as the authors know, the $2^\mathrm{nd}$ and $3^\mathrm{rd}$ questions have not been answered previously. The $1^\mathrm{st}$ and $2^\mathrm{nd}$ questions are answered in Section~\ref{sec_ball_uncontrolled}.  The answer to the $3^\mathrm{rd}$ question is highly nontrivial and is presented in a separate paper \cite{putkaradze2017optimal} by the authors. 

\begin{figure}[h]
	\centering
	\includegraphics[width=0.5\linewidth]{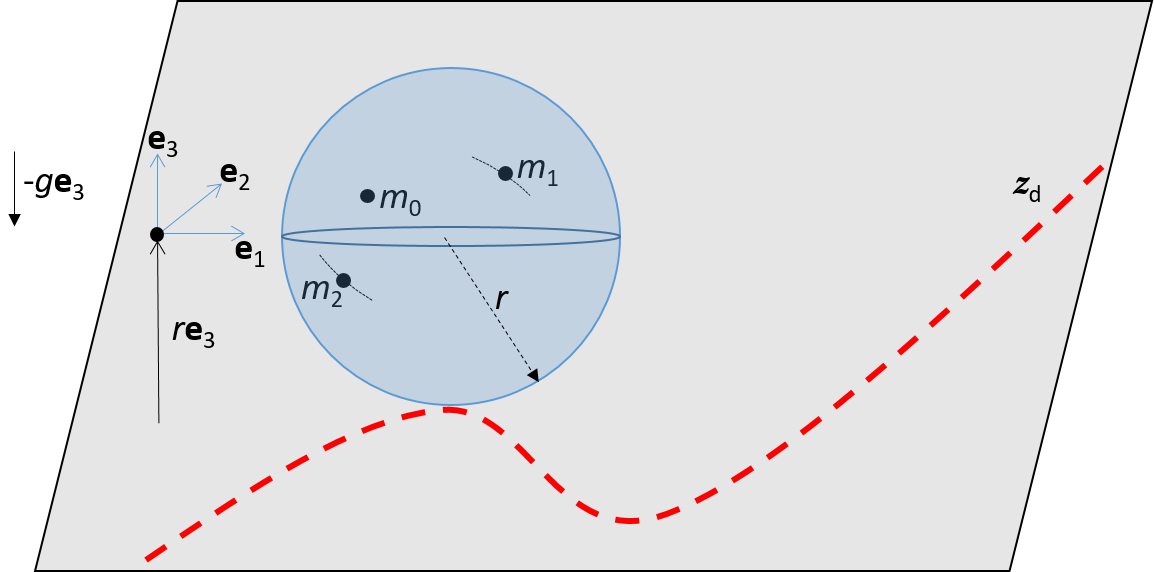}
	\caption{A ball of radius $r$ and mass $m_0$ rolls without slipping on a horizontal surface in the presence of a uniform gravitational field of magnitude $g$. The ball's center of mass is denoted by $m_0$. In addition, the ball contains $2$ internal point masses, $m_1$ and $m_2$, that may move within the ball. How must $m_1$ and $m_2$ be moved to induce the ball to follow the prescribed trajectory $\bz_\mathrm{d}$?}
	\label{fig:ball_rolling_along_traj}
\end{figure}
 
Two coordinate systems, or frames of reference, will be used to describe the motion of the rolling ball, an inertial spatial coordinate system and a body coordinate system in which each particle within the ball is always fixed. For brevity, the spatial coordinate system will be referred to as the spatial frame and the body coordinate system will be referred to as the body frame. These two frames are depicted in Figure~\ref{fig:detailed_rolling_ball}. The spatial frame has orthonormal axes $\mathbf{e}_1$, $\mathbf{e}_2$, $\mathbf{e}_3$, such that the $\mathbf{e}_1$-$\mathbf{e}_2$ plane is parallel to the horizontal surface and passes through the ball's geometric center (i.e. the $\mathbf{e}_1$-$\mathbf{e}_2$ plane is a height $r$ above the horizontal surface), such that $\mathbf{e}_3$ is vertical (i.e. $\mathbf{e}_3$ is perpendicular to the horizontal surface) and points ``upward" and away from the horizontal surface, and such that $\left(\mathbf{e}_1, \mathbf{e}_2, \mathbf{e}_3 \right)$ forms a right-handed coordinate system. For simplicity, the spatial frame axes are chosen to be
\begin{equation}
\mathbf{e}_1 = \begin{bmatrix} 1 & 0 & 0 \end{bmatrix}^\mathsf{T}, \quad \mathbf{e}_2 = \begin{bmatrix} 0 & 1 & 0 \end{bmatrix}^\mathsf{T}, \quad \mathrm{and} \quad \mathbf{e}_3 = \begin{bmatrix} 0 & 0 & 1 \end{bmatrix}^\mathsf{T}.
\end{equation}
The acceleration due to gravity in the uniform gravitational field is $\mathfrak{g} = -g \mathbf{e}_3  = \begin{bmatrix} 0 & 0 & -g  \end{bmatrix}^\mathsf{T}$ in the spatial frame.

The body frame's origin is chosen to coincide with the position of $m_0$'s center of mass. The body frame has orthonormal axes $\mathbf{E}_1$, $\mathbf{E}_2$, and $\mathbf{E}_3$, chosen to coincide with $m_0$'s principal axes, in which $m_0$'s inertia tensor $\inertia$ is diagonal, with corresponding principal moments of inertia $d_1$, $d_2$, and $d_3$. That is, in this body frame the inertia tensor is the diagonal matrix $\inertia = \diag \left( \begin{bmatrix} d_1 & d_2 & d_3 \end{bmatrix} \right)$.
\rem{\begin{equation}
\inertia_0 = \begin{bmatrix} d_1 & 0 & 0 \\ 0 & d_2 & 0 \\ 0 & 0 & d_3 \end{bmatrix}.
\end{equation}}
Moreover, $\mathbf{E}_1$, $\mathbf{E}_2$, and $\mathbf{E}_3$ are chosen so that $\left(\mathbf{E}_1, \mathbf{E}_2, \mathbf{E}_3 \right)$ forms a right-handed coordinate system. For simplicity, the body frame axes are chosen to be
\begin{equation}
\mathbf{E}_1 = \begin{bmatrix} 1 & 0 & 0 \end{bmatrix}^\mathsf{T}, \quad \mathbf{E}_2 = \begin{bmatrix} 0 & 1 & 0 \end{bmatrix}^\mathsf{T}, \quad \mathrm{and} \quad \mathbf{E}_3 = \begin{bmatrix} 0 & 0 & 1 \end{bmatrix}^\mathsf{T}.
\end{equation}

In the spatial frame, the body frame is the moving frame $\left(\Lambda \left(t\right) \mathbf{E}_1, \Lambda \left(t\right) \mathbf{E}_2, \Lambda \left(t\right) \mathbf{E}_3  \right)$, where $\Lambda \left(t\right) \in SO(3)$ defines the orientation (or attitude) of the ball at time $t$ relative to its reference configuration, for example at some initial time. For $1 \le i \le n$, it is assumed that $\bxi_i(t)$, the position of $m_i$'s center of mass, is expressed with respect to the body frame. Since $m_0$'s center of mass is always $\mathbf{0} =  \begin{bmatrix} 0 & 0 & 0  \end{bmatrix}^\mathsf{T} $ in the body frame (by choice of that frame's origin), let $\bxi_0 \equiv \mathbf{0}$; with this definition, $m_i$'s center of mass is located at $\bxi_i(t)$ for all $0 \le i \le n$.  

Let $\mathbf{z}_i(t)$ denote the position of $m_i$'s center of mass in the spatial frame. Let $\bchi_i(t)$ denote the body frame vector from the ball's geometric center to $m_i$'s center of mass. Then for $m_0$, $\bchi_0$ is the constant (time-independent) vector from the ball's geometric center to $m_0$'s center of mass. Note that the position of $m_i$'s center of mass in the body frame is $\bxi_i(t) = \bchi_i(t) -\bchi_0$ and in the spatial frame is $\mathbf{z}_i(t)=\mathbf{z}_0(t)+\Lambda(t) \bxi_i(t)=\mathbf{z}_0(t)+\Lambda(t) \left[\bchi_i(t)-\bchi_0\right]$. In general, a particle with position $\mathbf{w}(t)$ in the body frame has position $\mathbf{z}(t) = \mathbf{z}_0(t)+\Lambda(t) \mathbf{w}(t)$ in the spatial frame and has position $\mathbf{w}(t)+\bchi_0$ in the body frame translated to the ball's geometric center.

For conciseness, the ball's geometric center is often denoted GC, $m_0$'s center of mass is often denoted CM, and the ball's contact point with the surface is often denoted CP. The GC is located at $\mathbf{z}_\mathrm{GC}(t) = \mathbf{z}_0(t)-\Lambda(t) \bchi_0$ in the spatial frame, at $-\bchi_0$ in the body frame, and at $\mathbf{0}$ in the body frame translated to the GC. The CM is located at $\mathbf{z}_0(t)$ in the spatial frame, at $\mathbf{0}$ in the body frame, and at $\bchi_0$ in the body frame translated to the GC. The CP is located at $\mathbf{z}_\mathrm{CP}(t) = \mathbf{z}_0(t)-\Lambda(t) \left[r\bGam(t)+\bchi_0 \right]$ in the spatial frame, at $-\left[r\bGam(t)+\bchi_0 \right]$ in the body frame, and at $-r\bGam(t)$ in the body frame translated to the GC, where $\bGamma(t) \equiv \Lambda^{-1}(t) \mathbf{e}_3$. Since the third spatial coordinate of the ball's GC is always $0$ and of the ball's CP is always $-r$, only the first two spatial coordinates of the ball's GC and CP, denoted by $\bz(t)$, are needed to determine the spatial location of the ball's GC and CP.

For succintness, the explicit time dependence of variables is often dropped. That is, the orientation of the ball at time $t$ is denoted simply $\Lambda$ rather than $\Lambda(t)$, the position of $m_i$'s center of mass in the spatial frame at time $t$ is denoted $\mathbf{z}_i$ rather than $\mathbf{z}_i(t)$, the position of $m_i$'s center of mass in the body frame at time $t$ is denoted $\bxi_i$ rather than $\bxi_i(t)$, the position of $m_i$'s center of mass in the body frame translated to the GC at time $t$ is denoted $\bchi_i$ rather than $\bchi_i(t)$, and the spatial $\mathbf{e}_1$- and $\mathbf{e}_2$-components of the ball's GC and CP at time $t$ are denoted $\bz$ rather than $\bz(t)$.

\begin{figure}[h]
	\centering
	\includegraphics[width=0.5\linewidth]{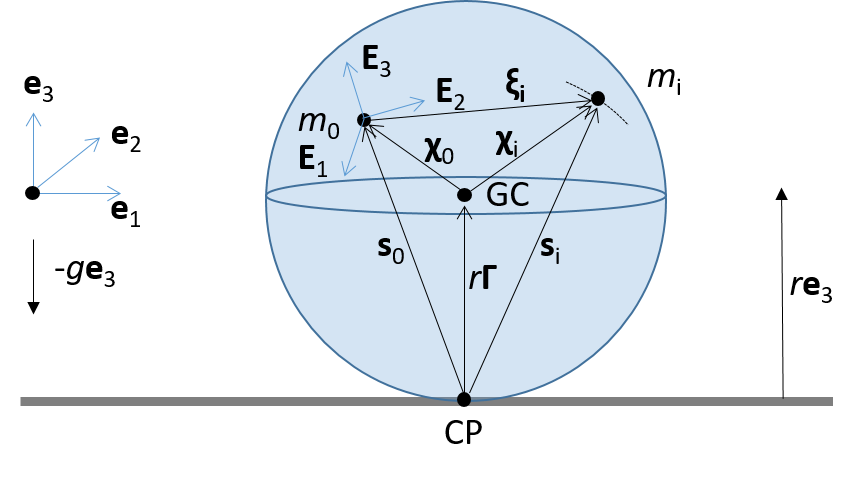}
	\caption{A ball of radius $r$ and mass $m_0$ rolls without slipping on a horizontal surface in the presence of a uniform gravitational field of magnitude $g$. The ball's geometric center, center of mass, and contact point with the horizontal surface are denoted by GC, $m_0$, and CP, respectively. The ball's motion is actuated by $n$ point masses, each of mass $m_i$, $1 \le i \le n$, that move inside the ball. The spatial frame has origin located at height $r$ above the horizontal surface and orthonormal axes $\mathbf{e}_1$, $\mathbf{e}_2$, and $\mathbf{e}_3$. The body frame has origin located at the ball's center of mass (denoted by $m_0$) and orthonormal axes $\mathbf{E}_1$, $\mathbf{E}_2$, and $\mathbf{E}_3$. All vectors inside the ball are expressed with respect to the body frame, while all vectors outside the ball are expressed with respect to the spatial frame. }
	\label{fig:detailed_rolling_ball}
\end{figure}

\section{Derivation of Lagrangian, Nonholonomic Constraint, and Variational Principle}
\label{sec_ball_uncontrolled}

This section derives the equations of motion for a rolling ball actuated by internal point masses. After developing prerequisites in Subsections~\ref{ssec_KE_PE_L} and \ref{ssec_rollconstraint_LdA}, Subsection~\ref{ssec_ball_uncontrolled} derives the equations of motion for a rolling ball actuated by internal point masses.  As special cases, Subsection~\ref{ssec_ball_uncontrolled} also derives the equations of motion for a rolling ball with static internal structure and the equations of motion for a rolling ball actuated by internal point masses that move along arbitrarily-shaped rails fixed within the ball. Finally, as an even more special case, Subsection~\ref{ssec_disk_uncontrolled} derives the  equation of motion for a rolling disk actuated by internal point masses that move along arbitrarily-shaped rails fixed within the disk.

\subsection{Kinetic Energy, Potential Energy, and Lagrangian} \label{ssec_KE_PE_L} 

As a first step to deriving the  equations of motion for the rolling ball, the ball's kinetic and potential energies must be constructed, from which the ball's Lagrangian is easily constructed.

\paragraph{Configuration Manifold and Constraints} 
Since the motion of the point masses with respect to the ball's frame are prescribed, the configuration manifold of the system consists of the group of rotations and translations, i.e. the space $SE(3)$. The ball's orientation matrix $\Lambda(t) \in SO(3)$ describes the rotation of the ball and the vector $\mathbf{z}_0(t) \in \mathbb{R}^3$ describes the translation of the ball's center of mass with respect to the fixed, spatial frame, so that $\left(\Lambda, \mathbf{z}_0\right) \in SE(3)$.  The Lagrangian depends, in general, on the variables $\Lambda$, $\dot \Lambda$, $\mathbf{z}_0$, and $\dot{\mathbf{z}}_0$. The Lagrangian  reduced with respect to the rotational symmetry, in the presence of gravity, depends on the variables $\bOm \equiv \left( \Lambda^{-1} \dot \Lambda \right)^\vee$, $\mathbf{Y}_0 \equiv \Lambda^{-1} \dot{\mathbf{z}}_0$, and $\bGamma \equiv \Lambda^{-1} \mathbf{e}_3$, where $\mathbf{e}_3$ is the unit vector along the vertical axis in the spatial frame. Here, we have used the hat map diffeomorphism ${\textvisiblespace}^\wedge$ between the vectors in $\mathbb{R}^3$ and the antisymmetric matrices in $\mso(3)$, given by $\widehat{a}_{ij}=-\epsilon_{i j k} a^k$, and ${\textvisiblespace}^\vee$ is the inverse of the hat map. For more details, we refer the reader to Appendix~\ref{app_background}, in particular, formulas \eqref{eq11}-\eqref{eq12}. For a more careful discussion of the configuration manifold and variational principles, see, for example, reference \cite{Ho2011_pII}. 
 
\paragraph{Kinetic Energy} 
For $0\le i \le n$, recall that $\mathbf{z}_i(t)$ denotes the spatial coordinates of the $i^\mathrm{th}$ mass, and  $\mathbf{Y}_i \equiv \Lambda^{-1} \dot{\mathbf{z}}_i$  is the linear  velocity of the $i^\mathrm{th}$ mass  measured in the body frame. By definition, $\bOm \equiv \left( \Lambda^{-1} \dot \Lambda \right)^\vee$ is the ball's body angular velocity. Remembering that $m_0$ and $\mathbb{I}$ denote the mass and inertia tensor, measured with respect to the center of mass, of the ball without the $n$ point masses, the  kinetic energy  of the ball without the $n$  point masses is the sum of its translational kinetic energy of and rotational kinetic energy about its center of mass:
\begin{equation}
T_0  = \frac{1}{2} m_0 \left| \dot{\mathbf{z}}_0 \right|^2+ \frac{1}{2} \left< \bOm , \inertia  \bOm \right> = \frac{1}{2} m_0 \left| \mathbf{Y}_0  \right|^2+ \frac{1}{2} \left< \bOm , \inertia  \bOm \right>.
\end{equation}
For $1\le i \le n$, since $m_i$ is a point mass, its kinetic energy is just its translational kinetic energy.  Therefore, the kinetic energy of the $i^\mathrm{th}$ point mass is given by
\begin{equation}
T_i  = \frac{1}{2} m_i \left| \dot{\mathbf{z}}_i \right|^2 =\frac{1}{2} m_i \left| \mathbf{Y}_i  \right|^2. 
\end{equation}
\rem{ %%%BEGIN REM But to be consistent with $m_0$'s kinetic energy formula, for $1\le i \le n$, $m_i$'s inertia tensor will be defined to be zero so that $\inertia_i \equiv \mathbf{0}$ and so that $m_i$'s kinetic energy is
\begin{equation}
T_i  = \frac{1}{2} m_i \left| \dot{\mathbf{z}}_i \right|^2+ \frac{1}{2} \left< \bOm , \inertia_i  \bOm \right>.
\end{equation}
Define $\mathbf{Y}_i \equiv \Lambda^{-1} \dot{\mathbf{z}}_i.$
Since $\left| \dot{\mathbf{z}}_i \right|^2 = \left| \Lambda^{-1} \dot{\mathbf{z}}_i \right|^2 = \left| \mathbf{Y}_i \right|^2$, $m_i$'s kinetic energy becomes
\begin{equation}
T_i =\frac{1}{2} m_i \left| \mathbf{Y}_i  \right|^2 + \frac{1}{2} \left< \bOm , \inertia_i  \bOm \right>.
\end{equation}
} %%%END REM 
Thus, the ball's total kinetic energy is
\begin{equation}
\begin{split}
T &=  \sum_{i=0}^n T_i   
%= \frac{1}{2} \sum_{i=0}^n m_i \left| \mathbf{Y}_i  \right|^2 + \frac{1}{2} \sum_{i=0}^n \left< \bOm , \inertia_i  \bOm \right> 
%&= \frac{1}{2} \sum_{i=0}^n m_i \left| \mathbf{Y}_i  \right|^2 + \frac{1}{2} \left< \bOm , \sum_{i=0}^n \inertia_i  \bOm \right> 
= \frac{1}{2} \sum_{i=0}^n m_i \left| \mathbf{Y}_i  \right|^2 + \frac{1}{2} \left< \bOm , \inertia  \bOm \right>.
\end{split}
\end{equation}
\rem{\begin{equation}
\inertia =  \inertia_0 = \begin{bmatrix} d_1 & 0 & 0 \\ 0 & d_2 & 0 \\ 0 & 0 & d_3 \end{bmatrix}.
\end{equation}}

\paragraph{Potential Energy}  
The potential energy due to mass $m_i$  is $V_i = m_i g \left<\bchi_i, \bGamma \right>$, where $\bGamma \equiv \Lambda^{-1} \mathbf{e}_3$. Thus, the ball's potential energy is 
\begin{equation}
\begin{split}
V = \sum_{i=0}^n V_i 
= \sum_{i=0}^n m_i g \left<\bchi_i, \bGamma \right> 
= g \left< \sum_{i=0}^n m_i \bchi_i, \bGamma \right>.
\end{split}
\end{equation}

\paragraph{Lagrangian} 
Since the spatial position of $m_i$'s center of mass is $\mathbf{z}_i=\mathbf{z}_0+\Lambda \left[\bchi_i -\bchi_0\right]$, the spatial velocity of $m_i$'s center of mass is $\dot{\mathbf{z}}_i=\dot{\mathbf{z}}_0+\dot{\Lambda} \left[\bchi_i -\bchi_0\right]+\Lambda \dot{\bchi}_i$. Hence,
\begin{equation}
\begin{split}
\mathbf{Y}_i &\equiv \Lambda^{-1} \dot{\mathbf{z}}_i 
= \Lambda^{-1} \left[ \dot{\mathbf{z}}_0+\dot{\Lambda} \left[\bchi_i-\bchi_0\right]+\Lambda \dot{\bchi}_i \right] 
%= \Lambda^{-1} \dot{\mathbf{z}}_0+\Lambda^{-1} \dot{\Lambda} \left[\bchi_i -\bchi_0\right]+\dot{\bchi}_i \\
%&= \mathbf{Y}_0 + \widehat \bOm \left[\bchi_i-\bchi_0\right]+\dot{\bchi}_i 
= \mathbf{Y}_0 + \bOm \times \left[\bchi_i-\bchi_0\right]+\dot{\bchi}_i. \\
\end{split}
\end{equation}
The ball's Lagrangian is the difference between its kinetic and potential energies:
\begin{equation}
\begin{split}
l \equiv T -  V  = \frac{1}{2} \sum_{i=0}^n m_i \left| \mathbf{Y}_i  \right|^2 + \frac{1}{2} \left< \bOm, \inertia  \bOm \right>  - g \left< \sum_{i=0}^n m_i \bchi_i, \bGamma\right>.
\end{split}
\end{equation}
Since $\mathbf{Y}_i$ can be expressed as a function of $\mathbf{Y}_0$ and $\bOm$ for $1\le i \le n$, note that the ball's Lagrangian should be expressed as $ l\left(\bOm,\mathbf{Y}_0 ,\bGamma\right)$, but this functional dependence is suppressed for concision.

\subsection{Rolling Constraint and Lagrange-d'Alembert's Principle} \label{ssec_rollconstraint_LdA}

Having constructed the rolling ball's Lagrangian, the variation of the action integral is now computed, taking into consideration the rolling constraint and Lagrange-d'Alembert's principle.

\paragraph{Rolling Constraint}  
Recall that it is assumed that the ball rolls along the horizontal surface without slipping. The vector pointing from the contact point (i.e. the point on the horizontal surface touching the bottom of the ball) to $m_0$'s center of mass (located at $\mathbf{z}_0$ in the spatial frame and at $\bxi_0 \equiv \mathbf{0}$ in the body frame) is
\begin{equation}
\bsigma_0 \equiv r \mathbf{e}_3 + \Lambda \bchi_0
\end{equation}
in the spatial frame and is 
\begin{equation}  \label{eq:s0}
\mathbf{s}_0 \equiv \Lambda^{-1}  \bsigma_0 = r \Lambda^{-1} \mathbf{e}_3 + \bchi_0 = r \bGamma + \bchi_0
\end{equation}
in the body frame. Differentiating \eqref{eq:s0} with respect to time, using the identity $\dot \bGamma=\bGamma \times \bOm$, and using the identity $-r \bGamma = \bchi_0-s_0$, which follows trivially from \eqref{eq:s0}, yields the following useful result:
\begin{equation} \label{eq:dot_s0}
{\dot {\mathbf{s}}}_0 = r \dot \bGamma = r \bGamma \times \bOm = \bOm \times \left(-r \bGamma \right) = \bOm \times \left( \bchi_0-\mathbf{s}_0 \right).
\end{equation}
Another useful result that follows trivially from \eqref{eq:s0} is 
\begin{equation} \label{eq:s0_times_bGamma}
\mathbf{s}_0 \times \bGamma = \left(r \bGamma + \bchi_0 \right) \times \bGamma
 %=r \bGamma \times \bGamma + \bchi_0 \times \bGamma
 = \bchi_0 \times \bGamma.
\end{equation}
The rolling constraint is imposed by stipulating that the contact point of the ball with the surface is at rest:
\begin{equation} \label{eq:roll_c0}
\dot {\mathbf{z}}_0 =\dot{\Lambda} \mathbf{s}_0 = \dot{\Lambda} \Lambda^{-1} \bsigma_0 = \widehat \bomega  \bsigma_0 = \bomega \times  \bsigma_0,
\end{equation}
where $\widehat {\bomega} \equiv \dot{\Lambda} \Lambda^{-1} = \Lambda \bOm \in \mso(3)$, or equivalently, by stipulating
\begin{equation} \label{eq:roll_c}
\mathbf{Y}_0  \equiv \Lambda^{-1} \dot{\mathbf{z}}_0 = \Lambda^{-1} \dot{\Lambda} \mathbf{s}_0   = \widehat \bOm \mathbf{s}_0 = \bOm \times \mathbf{s}_0.
\end{equation}
As a consequence of the rolling constraint \eqref{eq:roll_c},
\begin{equation}
\begin{split} \label{eq:roll_c_con}
\mathbf{Y}_i = \mathbf{Y}_0 + \bOm \times \left[\bchi_i-\bchi_0\right]+\dot{\bchi}_i &= \bOm \times \mathbf{s}_0 + \bOm \times \left[\bchi_i-\bchi_0\right]+\dot{\bchi}_i 
= \bOm \times \left[ \mathbf{s}_0 + \bchi_i-\bchi_0\right]+\dot{\bchi}_i.
\end{split}
\end{equation}

\paragraph{Lagrange-d'Alembert's Principle} 
\revision{R1Q3}{Letting $\delta {\mathbf{z}}_0$ denote the variational derivative of the spatial position of $m_0$'s center of mass and defining} $\bPsi \equiv \Lambda^{-1} \delta {\mathbf{z}}_0 $,
\begin{equation}
\begin{split}
\dot {\bPsi} &= \left[ \Lambda^{-1} \delta {\mathbf{z}}_0 \right]^\cdot 
= \left[ \Lambda^{-1} \right]^\cdot  \delta {\mathbf{z}}_0 + \Lambda^{-1} \left[ \delta {\mathbf{z}}_0 \right]^\cdot 
= -\Lambda^{-1} \dot{\Lambda} \Lambda^{-1}  \delta {\mathbf{z}}_0 + \Lambda^{-1} \delta {\dot {\mathbf{z}}_0} \\
&= - \widehat {\bOm} \bPsi + \Lambda^{-1} \delta {\dot {\mathbf{z}}_0} = - \bOm \times \bPsi + \Lambda^{-1} \delta {\dot {\mathbf{z}}_0}.
\end{split}
\end{equation}
Hence $\Lambda^{-1} \delta {\dot {\mathbf{z}}_0} = \dot {\bPsi} + \bOm \times \bPsi $.
Since $\mathbf{Y}_0 \equiv \Lambda^{-1} \dot{\mathbf{z}}_0 $,
\begin{equation}
\begin{split} \label{eq:de_Y0}
\delta \mathbf{Y}_0 &= \delta \left[ \Lambda^{-1} \dot{\mathbf{z}}_0  \right] 
= \left[ \delta \left( \Lambda^{-1} \right) \right] \dot{\mathbf{z}}_0 + \Lambda^{-1} \delta \dot{\mathbf{z}}_0 
= - \Lambda^{-1} \delta \Lambda \Lambda^{-1} \dot{\mathbf{z}}_0 + \Lambda^{-1} \delta \dot{\mathbf{z}}_0 \\
&= - \widehat {\bSigma} \mathbf{Y}_0 + \Lambda^{-1} \delta \dot{\mathbf{z}}_0 
= - \bSigma \times \mathbf{Y}_0 + \Lambda^{-1} \delta \dot{\mathbf{z}}_0 
= \dot {\bPsi} + \bOm \times \bPsi - \bSigma \times \mathbf{Y}_0,
\end{split}
\end{equation}
where $\widehat {\bSigma} \equiv \Lambda^{-1} \delta \Lambda \in \mso(3)$.

Since $\mathbf{Y}_i = \mathbf{Y}_0+ \bOm \times \left[\bchi_i -\bchi_0\right]+\dot{\bchi}_i$ and since the point masses move along prescribed trajectories $\left\{ \bchi_i \right\}_{i=0}^n$, so that the variation of $\mathbf{Y}_i$ is computed with respect to $\mathbf{Y}_0$ and $\bOm$, but not with respect to $\left\{ \bchi_i \right\}_{i=0}^n$: 
\begin{equation}
\begin{split} \label{eq:de_Yi}
\delta \mathbf{Y}_i &= \delta \mathbf{Y}_0 + \delta \bOm \times \left[\bchi_i-\bchi_0\right]. \\
\end{split}
\end{equation}
The variation $\de \bOm$ is still given by 
$\delta \bOm = \dot { \bSigma} + \bOm \times \bSigma$, which was derived in \eqref{sigma_dot_eq} to obtain the free rigid body equations of motion \eqref{EP_fin_rigid_body}.  
\rem{%%%BEGIN REM 
Since $\widehat {\bSigma} = \Lambda^{-1} \delta{\Lambda} $,
\begin{equation}
\begin{split}
\dot { \widehat {\bSigma}} = \left[ \Lambda^{-1} \delta \Lambda \right]^\cdot 
= \left[ \Lambda^{-1} \right]^\cdot \delta \Lambda + \Lambda^{-1} \left[ \delta{\Lambda} \right]^\cdot 
= - \Lambda^{-1} \dot \Lambda \Lambda^{-1}  \delta \Lambda + \Lambda^{-1} \delta{\dot \Lambda} 
= - \widehat {\bOm} \widehat {\bSigma} + \Lambda^{-1} \delta{\dot \Lambda}.
\end{split}
\end{equation}
Hence $\Lambda^{-1} \delta{\dot \Lambda} = \dot { \widehat \bSigma} + \widehat {\bOm} \widehat {\bSigma}$.
Since $\widehat {\bOm} = \Lambda^{-1} \dot{\Lambda} $,
\begin{equation}
\begin{split}
\delta \widehat {\bOm} &= \delta \left[ \Lambda^{-1} \dot{\Lambda} \right] 
=  \left[ \delta \left( \Lambda^{-1} \right) \right] \dot{\Lambda} + \Lambda^{-1} \delta \dot{\Lambda} 
= -\Lambda^{-1} \delta \Lambda \Lambda^{-1} \dot{\Lambda} + \Lambda^{-1} \delta \dot{\Lambda} 
= - \widehat {\bSigma} \widehat {\bOm} + \dot { \widehat \bSigma} + \widehat {\bOm} \widehat {\bSigma} \\
&= \dot { \widehat \bSigma} + \left[ \bOm \times \bSigma \right]^\wedge 
= \left[\dot \bSigma + \bOm \times \bSigma \right]^\wedge.
\end{split}
\end{equation}
Hence 
\begin{equation} \label{eq:de_Omega}
\delta \bOm = \dot { \bSigma} + \bOm \times \bSigma. 
\end{equation}
} %%%END REM 
Furthermore, since $\bGamma \equiv \Lambda^{-1} \mathbf{e}_3$, the variation $\de \bGam$ is still given by $\de \bGam=\bGamma \times \bSigma$, which was derived in \eqref{eq_de_Gamma} to obtain the heavy top equations of motion \eqref{eq_heavy_top}.
\rem{ %%%BEGIN REM 
\begin{equation}
\begin{split} \label{eq:de_Gamma}
\delta \bGamma = \left[ \delta \left( \Lambda^{-1} \right) \right] \mathbf{e}_3 
= -\Lambda^{-1} \delta \Lambda \Lambda^{-1} \mathbf{e}_3 
= -\widehat {\bSigma} \bGamma
= - \bSigma \times \bGamma
= \bGamma \times \bSigma.
\end{split}
\end{equation}
} %%%END REM

We now invoke Lagrange-d'Alembert's principle from Subappendix~\ref{ssec_LdA}. Part of Lagrange-d'Alembert's principle stipulates that due to the rolling constraint \eqref{eq:roll_c0}, which says $\dot {\mathbf{z}}_0 =\dot{\Lambda} \mathbf{s}_0$, the variations of ${\mathbf{z}}_0$ must have the form $\de {\mathbf{z}}_0 =\de {\Lambda} \mathbf{s}_0$. Hence, the variations $\bPsi \equiv \Lambda^{-1} \delta {\mathbf{z}}_0$ must take on the following form (as a consequence of the rolling constraint \eqref{eq:roll_c0} and Lagrange-d'Alembert's principle):
\begin{equation} \label{eq:ldA_rc}
\bPsi \equiv \Lambda^{-1} \delta {\mathbf{z}}_0 = \Lambda^{-1} \de {\Lambda} \mathbf{s}_0 =  \widehat \bSigma \mathbf{s}_0  = \bSigma \times \mathbf{s}_0. 
\end{equation}
The  equations of motion are derived here and in the next section from Lagrange-d'Alembert's principle. Recalling that the  point masses move along prescribed trajectories $\left\{ \bchi_i \right\}_{i=0}^n$, it is important to keep in mind that the variation of the action integral is computed with respect to $\left\{ \mathbf{Y}_i \right\}_{i=0}^n$, $\bOm$, and $\bGamma$, but \textbf{not} with respect to $\left\{ \bchi_i \right\}_{i=0}^n$. Once the variation of the action integral is computed, tedious calculations are performed to isolate $\bSigma$, after which the variation of the action integral is equated to zero in order to obtain the equations of motion. Key points in the calculations \revision{R1Q1}{after computing the variation of the action integral} are: \revision{R1Q1}{1) the rolling constraint is enforced by invoking \eqref{eq:roll_c} and \eqref{eq:roll_c_con},} 2) the variations $\bPsi$ and $\bSigma$ must satisfy \rem{$\bPsi = \bSigma \times \mathbf{s}_0$ \eqref{eq:ldA_rc}} \revision{R1Q4}{\eqref{eq:ldA_rc}}, which enforces the constraints on the variations demanded by Lagrange-d'Alembert's principle, and 3) the variation $\bSigma$ must also satisfy $\bSigma(a)=\bSigma(b)=\mathbf{0}$, which enforces the vanishing endpoint constraints. To begin the calculations, the variation of the action integral is computed as
\begin{equation}
\begin{split} \label{eq:rball_vars_p1}
\de S = \de \int_a^b l \mathrm{d} t = \int_a^b \delta l \mathrm{d} t = \int_a^b \left[ \sum_{i=0}^n m_i \left< \mathbf{Y}_i  , \delta \mathbf{Y}_i \right>
+ \left< \inertia  \bOm, \delta \bOm \right> - g \left< \sum_{i=0}^n m_i \bchi_i , \delta \bGamma \right> \right] \mathrm{d} t.
\end{split}
\end{equation}
\rem{Using the identities $\delta \mathbf{Y}_i = \delta \mathbf{Y}_0 + \delta \bOm \times \left[\bchi_i-\bchi_0\right]$ \eqref{eq:de_Yi}, $\de \mathbf{Y}_0= \dot {\bPsi} + \bOm \times \bPsi - \bSigma \times \mathbf{Y}_0$ \eqref{eq:de_Y0}, $\delta \bOm = \dot { \bSigma} + \bOm \times \bSigma$ \eqref{sigma_dot_eq}, and $\de \bGamma= \bGamma \times \bSigma$ \eqref{eq_de_Gamma}} \revision{R1Q4}{Using the identities \eqref{eq:de_Yi}, \eqref{eq:de_Y0}, \eqref{sigma_dot_eq}, and \eqref{eq_de_Gamma}} and integrating by parts, the variation of the action integral obtained in \eqref{eq:rball_vars_p1} becomes
\begin{equation}
\begin{split} \label{eq:rball_vars_p4} \hspace{-1mm}
\de S
&= \int_a^b \Bigg[ -\sum_{i=0}^n m_i \left<\left(\dd{}{t}+\bOm \times \right) \mathbf{Y}_i , \bPsi \right> \\ & \hphantom{=\int_a^b \Bigg[} +\left<-\left(\dd{}{t}+\bOm \times \right) \left[ \inertia  \bOm+\sum_{i=0}^n m_i \left[\bchi_i-\bchi_0\right] \times \mathbf{Y}_i \right] +\sum_{i=0}^n m_i \left( \mathbf{Y}_i \times \mathbf{Y}_0+ g \bGamma \times \bchi_i  \right) , \bSigma \right> \Bigg] \mathrm{d} t \\
& \hphantom{=} +\left. \sum_{i=0}^n m_i \left< \mathbf{Y}_i  , \bPsi \right>\right|_a^b + \left. \left< \inertia  \bOm+\sum_{i=0}^n m_i \left[\bchi_i-\bchi_0\right] \times \mathbf{Y}_i, \bSigma \right> \right|_a^b.
\end{split}
\end{equation}
\rem{ %%%BEGIN REM 
Using the identity $\delta \mathbf{Y}_i = \delta \mathbf{Y}_0 + \delta \bOm \times \left[\bchi_i-\bchi_0\right]$ \eqref{eq:de_Yi}, the variation of the action integral obtained in \eqref{eq:rball_vars_p1} becomes
\begin{equation}
\begin{split} \label{eq:rball_vars_p2}
\de S &= \int_a^b \left[ \sum_{i=0}^n m_i \left< \mathbf{Y}_i  , \delta \mathbf{Y}_0 + \delta \bOm \times \left[\bchi_i-\bchi_0\right] \right>
+ \left< \inertia  \bOm, \delta \bOm \right> - g \left< \sum_{i=0}^n m_i \bchi_i , \delta \bGamma \right> \right] \mathrm{d} t \\
&= \int_a^b \left[ \sum_{i=0}^n m_i \left< \mathbf{Y}_i  , \delta \mathbf{Y}_0 \right>
+ \left< \inertia  \bOm+\sum_{i=0}^n m_i \left[\bchi_i-\bchi_0\right] \times \mathbf{Y}_i, \delta \bOm \right> - g \left< \sum_{i=0}^n m_i \bchi_i , \delta \bGamma \right> \right] \mathrm{d} t.
\end{split}
\end{equation}
Using the identities $\de \mathbf{Y}_0= \dot {\bPsi} + \bOm \times \bPsi - \bSigma \times \mathbf{Y}_0$ \eqref{eq:de_Y0}, $\delta \bOm = \dot { \bSigma} + \bOm \times \bSigma$ \eqref{eq:de_Omega}, and $\de \bGamma= \bGamma \times \bSigma$ \eqref{eq:de_Gamma}, the variation of the action integral obtained in \eqref{eq:rball_vars_p2} becomes
\begin{equation}
\begin{split} \label{eq:rball_vars_p3}
\de S &= \int_a^b \Bigg[ \sum_{i=0}^n m_i \left< \mathbf{Y}_i  , \dot {\bPsi} + \bOm \times \bPsi - \bSigma \times \mathbf{Y}_0 \right>
+ \left< \inertia  \bOm+\sum_{i=0}^n m_i \left[\bchi_i-\bchi_0\right] \times \mathbf{Y}_i, \dot { \bSigma} + \bOm \times \bSigma \right> \\
& \hphantom{= \int_a^b \Bigg[}  + g \left< \sum_{i=0}^n m_i \bchi_i , \bSigma \times \bGamma \right> \Bigg] \mathrm{d} t.
\end{split}
\end{equation}
Integrating by parts, the variation of the action integral obtained in \eqref{eq:rball_vars_p3} becomes 
\begin{equation}
\begin{split} \label{eq:rball_vars_p4} \hspace{-1mm}
\de S &= \int_a^b \Bigg[ \sum_{i=0}^n m_i \left[ -\left<\left(\dd{}{t}+\bOm \times \right) \mathbf{Y}_i , \bPsi \right>+\left<\mathbf{Y}_i \times \mathbf{Y}_0 , \bSigma \right> \right] \\
& \hphantom{=\int_a^b \Bigg[} - \left<\left(\dd{}{t}+\bOm \times \right) \left[ \inertia  \bOm+\sum_{i=0}^n m_i \left[\bchi_i-\bchi_0\right] \times \mathbf{Y}_i \right], \bSigma \right> 
+ g \left< \bGamma \times \sum_{i=0}^n m_i \bchi_i , \bSigma \right> \Bigg] \mathrm{d} t \\
& \hphantom{=} +\left. \sum_{i=0}^n m_i \left< \mathbf{Y}_i  , \bPsi \right>\right|_a^b + \left. \left< \inertia  \bOm+\sum_{i=0}^n m_i \left[\bchi_i-\bchi_0\right] \times \mathbf{Y}_i, \bSigma \right> \right|_a^b   \\
&= \int_a^b \Bigg[ -\sum_{i=0}^n m_i \left<\left(\dd{}{t}+\bOm \times \right) \mathbf{Y}_i , \bPsi \right> \\ & \hphantom{=\int_a^b \Bigg[} +\left<-\left(\dd{}{t}+\bOm \times \right) \left[ \inertia  \bOm+\sum_{i=0}^n m_i \left[\bchi_i-\bchi_0\right] \times \mathbf{Y}_i \right] +\sum_{i=0}^n m_i \left( \mathbf{Y}_i \times \mathbf{Y}_0+ g \bGamma \times \bchi_i  \right) , \bSigma \right> \Bigg] \mathrm{d} t \\
& \hphantom{=} +\left. \sum_{i=0}^n m_i \left< \mathbf{Y}_i  , \bPsi \right>\right|_a^b + \left. \left< \inertia  \bOm+\sum_{i=0}^n m_i \left[\bchi_i-\bchi_0\right] \times \mathbf{Y}_i, \bSigma \right> \right|_a^b.
\end{split}
\end{equation} 
Next, \eqref{eq:rball_vars_p4} is evaluated on the constraint distribution given by $\bPsi = \bSigma \times \mathbf{s}_0$ \eqref{eq:ldA_rc}, and the boundary terms in \eqref{eq:rball_vars_p4} are eliminated since $\bSigma$ is a variation such that $\bSigma(a)=\bSigma(b)=\mathbf{0}$. With these manipulations, the variation of the action integral obtained in \eqref{eq:rball_vars_p4} becomes
\begin{equation}
\begin{split} \label{eq:rball_vars_p5}
\de S &= \int_a^b \Bigg[ - \left<\left(\dd{}{t}+\bOm \times \right) \left[ \sum_{i=0}^n m_i \mathbf{Y}_i \right], \bSigma \times \mathbf{s}_0  \right> \\ 
& \hphantom{=\int_a^b \Bigg[} +\left<-\left(\dd{}{t}+\bOm \times \right) \left[ \inertia  \bOm+\sum_{i=0}^n m_i \left[\bchi_i-\bchi_0\right] \times \mathbf{Y}_i \right] +\sum_{i=0}^n m_i \left( \mathbf{Y}_i \times \mathbf{Y}_0+ g \bGamma \times \bchi_i  \right) , \bSigma \right> \Bigg] \mathrm{d} t \\
&= \int_a^b \Bigg< \left( \left(\dd{}{t}+\bOm \times \right) \left[ \sum_{i=0}^n m_i \mathbf{Y}_i \right] \right) \times \mathbf{s}_0 -\left(\dd{}{t}+\bOm \times \right) \left[ \inertia  \bOm +\sum_{i=0}^n m_i \left[\bchi_i-\bchi_0\right] \times \mathbf{Y}_i \right] \\ & \hphantom{=\int_a^b \Bigg<} +\sum_{i=0}^n m_i \left( \mathbf{Y}_i \times \mathbf{Y}_0+ g \bGamma \times \bchi_i  \right) , \bSigma \Bigg> \mathrm{d} t \\
&= \int_a^b \Bigg<-\left(\dd{}{t}+\bOm \times \right) \left[ \inertia  \bOm +\sum_{i=0}^n m_i \left[\mathbf{s}_0+\bchi_i-\bchi_0\right] \times \mathbf{Y}_i \right] \\ & \hphantom{=\int_a^b \Bigg<} +\sum_{i=0}^n m_i \left( \mathbf{Y}_i \times \left( \mathbf{Y}_0 -  \left(\dd{}{t}+\bOm \times \right) \mathbf{s}_0  \right)+ g \bGamma \times \bchi_i  \right) , \bSigma \Bigg> \mathrm{d} t.
\end{split}
\end{equation}
Using the identity $\mathbf{Y}_0= \bOm \times \mathbf{s}_0$ \eqref{eq:roll_c}, the variation of the action integral obtained in \eqref{eq:rball_vars_p5} becomes
\begin{equation}
\begin{split} \label{eq:rball_vars_p6}
\de S &= \int_a^b \left<-\left(\dd{}{t}+\bOm \times \right) \left[ \inertia  \bOm +\sum_{i=0}^n m_i \left[\mathbf{s}_0+\bchi_i-\bchi_0\right] \times \mathbf{Y}_i \right] +\sum_{i=0}^n m_i \left( {\dot {\mathbf{s}} }_0 \times \mathbf{Y}_i + g \bGamma \times \bchi_i  \right) , \bSigma \right> \mathrm{d} t.
\end{split}
\end{equation}
Finally, using the identities $\mathbf{Y}_i= \bOm \times \left[ \mathbf{s}_0 + \bchi_i-\bchi_0\right]+\dot{\bchi}_i$ \eqref{eq:roll_c_con} and ${\dot {\mathbf{s}}}_0 = \bOm \times \left( \bchi_0-\mathbf{s}_0 \right)$ \eqref{eq:dot_s0}, the variation of the action integral obtained in \eqref{eq:rball_vars_p6} becomes
\begin{equation}
\begin{split} \label{eq:rball_vars_p7}
\de S &= \int_a^b \Bigg<-\left(\dd{}{t}+\bOm \times \right) \left[ \inertia  \bOm +\sum_{i=0}^n m_i \left[\mathbf{s}_0+\bchi_i-\bchi_0\right] \times \left[\bOm \times \left[ \mathbf{s}_0 + \bchi_i-\bchi_0\right]+\dot{\bchi}_i \right] \right] \\ & \hphantom{=\int_a^b \Bigg<} +\sum_{i=0}^n m_i \left( \left[ \left( \mathbf{s}_0-\bchi_0\right) \times \bOm \right]  \times \left[ \bOm \times \left[ \mathbf{s}_0 + \bchi_i-\bchi_0\right]+\dot{\bchi}_i \right] + g \bGamma \times \bchi_i  \right) , \bSigma \Bigg> \mathrm{d} t \\
&= \int_a^b \Bigg<-\left(\dd{}{t}+\bOm \times \right) \left[ \inertia  \bOm +\sum_{i=0}^n m_i \left[\mathbf{s}_0+\bchi_i-\bchi_0\right] \times \left[\bOm \times \left[ \mathbf{s}_0 + \bchi_i-\bchi_0\right]+\dot{\bchi}_i \right] \right] \\ & \hphantom{=\int_a^b \Bigg<} +\sum_{i=0}^n m_i \left( \left[ \left( \mathbf{s}_0-\bchi_0\right) \times \bOm \right]  \times \left[ \bOm \times \bchi_i+\dot{\bchi}_i \right] + g \bGamma \times \bchi_i  \right) , \bSigma \Bigg> \mathrm{d} t.
\end{split}
\end{equation}}
\rem{Evaluating \eqref{eq:rball_vars_p4} on the constraint distribution given by $\bPsi = \bSigma \times \mathbf{s}_0$ \eqref{eq:ldA_rc}, eliminating the boundary terms in \eqref{eq:rball_vars_p4} since $\bSigma$ is a variation such that $\bSigma(a)=\bSigma(b)=\mathbf{0}$, and using the identities $\mathbf{Y}_0= \bOm \times \mathbf{s}_0$ \eqref{eq:roll_c}, $\mathbf{Y}_i= \bOm \times \left[ \mathbf{s}_0 + \bchi_i-\bchi_0\right]+\dot{\bchi}_i$ \eqref{eq:roll_c_con}, and ${\dot {\mathbf{s}}}_0 = \bOm \times \left( \bchi_0-\mathbf{s}_0 \right)$ \eqref{eq:dot_s0}} \revision{R1Q4}{Evaluating \eqref{eq:rball_vars_p4} on the constraint distribution given by \eqref{eq:ldA_rc}, eliminating the boundary terms in \eqref{eq:rball_vars_p4} since $\bSigma$ is a variation such that $\bSigma(a)=\bSigma(b)=\mathbf{0}$, and using the identities \eqref{eq:roll_c}, \eqref{eq:roll_c_con}, and \eqref{eq:dot_s0}}, the variation of the action integral obtained in \eqref{eq:rball_vars_p4} becomes
\begin{equation}
\begin{split} \label{eq:rball_vars_p7}
\de S 
&= \int_a^b \Bigg<-\left(\dd{}{t}+\bOm \times \right) \left[ \inertia  \bOm +\sum_{i=0}^n m_i \left[\mathbf{s}_0+\bchi_i-\bchi_0\right] \times \left[\bOm \times \left[ \mathbf{s}_0 + \bchi_i-\bchi_0\right]+\dot{\bchi}_i \right] \right] \\ & \hphantom{=\int_a^b \Bigg<} +\sum_{i=0}^n m_i \left( \left[ \left( \mathbf{s}_0-\bchi_0\right) \times \bOm \right]  \times \left[ \bOm \times \bchi_i+\dot{\bchi}_i \right] + g \bGamma \times \bchi_i  \right) , \bSigma \Bigg> \mathrm{d} t.
\end{split}
\end{equation}
\rem{ %%%% BEGIN REM
	\begin{align} \label{eq:rball_vars}
	%\begin{equation}
	%\begin{split} \label{eq:rball_vars}
	0 &= \de S = \de \int_a^b l \mathrm{d} t = \int_a^b \delta l \mathrm{d} t \nonumber \\
	&= \int_a^b \left[ \sum_{i=0}^n m_i \left< \mathbf{Y}_i  , \delta \mathbf{Y}_i \right>
	+ \left< \inertia  \bOm, \delta \bOm \right> - g \left< \sum_{i=0}^n m_i \bchi_i , \delta \bGamma \right> \right] \mathrm{d} t \nonumber \\
	&= \int_a^b \left[ \sum_{i=0}^n m_i \left< \mathbf{Y}_i  , \delta \mathbf{Y}_0 + \delta \bOm \times \left[\bchi_i-\bchi_0\right] \right>
	+ \left< \inertia  \bOm, \delta \bOm \right> - g \left< \sum_{i=0}^n m_i \bchi_i , \delta \bGamma \right> \right] \mathrm{d} t \nonumber \\
	&= \int_a^b \left[ \sum_{i=0}^n m_i \left< \mathbf{Y}_i  , \delta \mathbf{Y}_0 \right>
	+ \left< \inertia  \bOm+\sum_{i=0}^n m_i \left[\bchi_i-\bchi_0\right] \times \mathbf{Y}_i, \delta \bOm \right> - g \left< \sum_{i=0}^n m_i \bchi_i , \delta \bGamma \right> \right] \mathrm{d} t \nonumber \\
	&= \int_a^b \Bigg[ \sum_{i=0}^n m_i \left< \mathbf{Y}_i  , \dot {\bPsi} + \bOm \times \bPsi - \bSigma \times \mathbf{Y}_0 \right>
	+ \left< \inertia  \bOm+\sum_{i=0}^n m_i \left[\bchi_i-\bchi_0\right] \times \mathbf{Y}_i, \dot { \bSigma} + \bOm \times \bSigma \right>  \nonumber \\
	& \hphantom{=}  + g \left< \sum_{i=0}^n m_i \bchi_i , \bSigma \times \bGamma \right> \Bigg] \mathrm{d} t \nonumber \\
	&= \int_a^b \Bigg[ \sum_{i=0}^n m_i \left[ -\left<\left(\dd{}{t}+\bOm \times \right) \mathbf{Y}_i , \bPsi \right>+\left<\mathbf{Y}_i \times \mathbf{Y}_0 , \bSigma \right> \right] \\
	& \hphantom{=} - \left<\left(\dd{}{t}+\bOm \times \right) \left[ \inertia  \bOm+\sum_{i=0}^n m_i \left[\bchi_i-\bchi_0\right] \times \mathbf{Y}_i \right], \bSigma \right> 
	+ g \left< \bGamma \times \sum_{i=0}^n m_i \bchi_i , \bSigma \right> \Bigg] \mathrm{d} t \nonumber \\
	&= \int_a^b \Bigg[ -\sum_{i=0}^n m_i \left<\left(\dd{}{t}+\bOm \times \right) \mathbf{Y}_i , \bPsi \right> \nonumber \\ & \hphantom{=} +\left<-\left(\dd{}{t}+\bOm \times \right) \left[ \inertia  \bOm+\sum_{i=0}^n m_i \left[\bchi_i-\bchi_0\right] \times \mathbf{Y}_i \right] +\sum_{i=0}^n m_i \left( \mathbf{Y}_i \times \mathbf{Y}_0+ g \bGamma \times \bchi_i  \right) , \bSigma \right> \Bigg] \mathrm{d} t \nonumber \\
	&= \int_a^b \Bigg[ - \left<\left(\dd{}{t}+\bOm \times \right) \left[ \sum_{i=0}^n m_i \mathbf{Y}_i \right], \bSigma \times \mathbf{s}_0  \right> \nonumber \\ 
	& \hphantom{=} +\left<-\left(\dd{}{t}+\bOm \times \right) \left[ \inertia  \bOm+\sum_{i=0}^n m_i \left[\bchi_i-\bchi_0\right] \times \mathbf{Y}_i \right] +\sum_{i=0}^n m_i \left( \mathbf{Y}_i \times \mathbf{Y}_0+ g \bGamma \times \bchi_i  \right) , \bSigma \right> \Bigg] \mathrm{d} t \nonumber \\
	&= \int_a^b \Bigg< \left( \left(\dd{}{t}+\bOm \times \right) \left[ \sum_{i=0}^n m_i \mathbf{Y}_i \right] \right) \times \mathbf{s}_0 -\left(\dd{}{t}+\bOm \times \right) \left[ \inertia  \bOm +\sum_{i=0}^n m_i \left[\bchi_i-\bchi_0\right] \times \mathbf{Y}_i \right]  \nonumber \nonumber \\ & \hphantom{=} +\sum_{i=0}^n m_i \left( \mathbf{Y}_i \times \mathbf{Y}_0+ g \bGamma \times \bchi_i  \right) , \bSigma \Bigg> \mathrm{d} t \nonumber \\
	&= \int_a^b \Bigg<-\left(\dd{}{t}+\bOm \times \right) \left[ \inertia  \bOm +\sum_{i=0}^n m_i \left[\mathbf{s}_0+\bchi_i-\bchi_0\right] \times \mathbf{Y}_i \right] \nonumber \\ & \hphantom{=} +\sum_{i=0}^n m_i \left( \mathbf{Y}_i \times \left( \mathbf{Y}_0 -  \left(\dd{}{t}+\bOm \times \right) \mathbf{s}_0  \right)+ g \bGamma \times \bchi_i  \right) , \bSigma \Bigg> \mathrm{d} t \nonumber \\
	&= \int_a^b \left<-\left(\dd{}{t}+\bOm \times \right) \left[ \inertia  \bOm +\sum_{i=0}^n m_i \left[\mathbf{s}_0+\bchi_i-\bchi_0\right] \times \mathbf{Y}_i \right] +\sum_{i=0}^n m_i \left( {\dot {\mathbf{s}} }_0 \times \mathbf{Y}_i + g \bGamma \times \bchi_i  \right) , \bSigma \right> \mathrm{d} t \nonumber \\
	&= \int_a^b \Bigg<-\left(\dd{}{t}+\bOm \times \right) \left[ \inertia  \bOm +\sum_{i=0}^n m_i \left[\mathbf{s}_0+\bchi_i-\bchi_0\right] \times \left[\bOm \times \left[ \mathbf{s}_0 + \bchi_i-\bchi_0\right]+\dot{\bchi}_i \right] \right] \nonumber \\ & \hphantom{=} +\sum_{i=0}^n m_i \left( \left[ \left( \mathbf{s}_0-\bchi_0\right) \times \bOm \right]  \times \left[ \bOm \times \left[ \mathbf{s}_0 + \bchi_i-\bchi_0\right]+\dot{\bchi}_i \right] + g \bGamma \times \bchi_i  \right) , \bSigma \Bigg> \mathrm{d} t \nonumber \\
	&= \int_a^b \Bigg<-\left(\dd{}{t}+\bOm \times \right) \left[ \inertia  \bOm +\sum_{i=0}^n m_i \left[\mathbf{s}_0+\bchi_i-\bchi_0\right] \times \left[\bOm \times \left[ \mathbf{s}_0 + \bchi_i-\bchi_0\right]+\dot{\bchi}_i \right] \right] \nonumber \\ & \hphantom{=} +\sum_{i=0}^n m_i \left( \left[ \left( \mathbf{s}_0-\bchi_0\right) \times \bOm \right]  \times \left[ \bOm \times \bchi_i+\dot{\bchi}_i \right] + g \bGamma \times \bchi_i  \right) , \bSigma \Bigg> \mathrm{d} t. \nonumber
	%\end{split}
	%\end{equation}
	\end{align}
	In the fifth equality, the identity $\delta \mathbf{Y}_i = \delta \mathbf{Y}_0 + \delta \bOm \times \left[\bchi_i-\bchi_0\right]$ \eqref{eq:de_Yi} is used. In the seventh equality, the identities $\de \mathbf{Y}_0= \dot {\bPsi} + \bOm \times \bPsi - \bSigma \times \mathbf{Y}_0$ \eqref{eq:de_Y0}, $\delta \bOm = \dot { \bSigma} + \bOm \times \bSigma$ \eqref{eq:de_Omega}, and $\de \bGamma= \bGamma \times \bSigma$ \eqref{eq:de_Gamma} are used. In the eighth equality, integration by parts is performed. In the tenth equality, the identity $\bPsi = \bSigma \times \mathbf{s}_0$ \eqref{eq:ldA_rc} is used. In the third to last equality, the identity $\mathbf{Y}_0= \bOm \times \mathbf{s}_0$ \eqref{eq:roll_c} is used. In the second to last equality, the identities $\mathbf{Y}_i= \bOm \times \left[ \mathbf{s}_0 + \bchi_i-\bchi_0\right]+\dot{\bchi}_i$ \eqref{eq:roll_c_con} and ${\dot {\mathbf{s}}}_0 = \bOm \times \left( \bchi_0-\mathbf{s}_0 \right)$ \eqref{eq:dot_s0} are used.
} %%%% END REM
Note the order in which the operations were performed: first variations and simplifications were computed in \eqref{eq:rball_vars_p1}-\eqref{eq:rball_vars_p4}, followed by evaluation of the result \eqref{eq:rball_vars_p4} on the constraint distribution \eqref{eq:ldA_rc} to obtain \eqref{eq:rball_vars_p7}; preserving this order is key to the correct application of Lagrange-d'Alembert's principle.

Now suppose a time-varying external force $\mathbf{F}_\mathrm{e}$ acts at the ball's geometric center. For example, this force might be due to the wind blowing on the ball when the ball rolls around outdoors. If the ball's geometric center in the spatial frame is $\mathbf{z}_\mathrm{GC}$, then the rolling constraint says that ${\dot {\mathbf{z}}}_\mathrm{GC} = \dot{\Lambda} \Lambda^{-1} r \mathbf{e}_3$  and Lagrange-d'Alembert's principle says that $\de \mathbf{z}_\mathrm{GC} = \de \Lambda \Lambda^{-1} r \mathbf{e}_3$. Application of the external force yields a new variation of the action integral, $\de S_1 = \de S+\int_a^b \left<\mathbf{F}_\mathrm{e},\de \mathbf{z}_\mathrm{GC} \right> \mathrm{d} t$, using Lagrange-d'Alembert's principle for incorporating external forces into the variational principle. Performing calculations on the new variation of the action integral to isolate $\bSigma$ gives:
\begin{equation} \label{eq_dS1}
\begin{split} 
\de S_1 
&= \de S+\int_a^b \left<\mathbf{F}_\mathrm{e},\de \mathbf{z}_\mathrm{GC} \right> \mathrm{d} t 
= \de S+\int_a^b \left<\Lambda^{-1} \mathbf{F}_\mathrm{e},\Lambda^{-1} \de \mathbf{z}_\mathrm{GC} \right> \mathrm{d} t \\
&= \de S+\int_a^b \left<\Lambda^{-1} \mathbf{F}_\mathrm{e},\Lambda^{-1} \de \Lambda \Lambda^{-1} r \mathbf{e}_3 \right> \mathrm{d} t 
= \de S+\int_a^b \left<\tilde \bGamma,\widehat \bSigma r \bGamma \right> \mathrm{d} t \\
&= \de S+\int_a^b \left<\tilde \bGamma,\bSigma \times r \bGamma \right> \mathrm{d} t 
= \de S+\int_a^b \left<r \bGamma \times \tilde \bGamma,\bSigma \right> \mathrm{d} t \\
&= \int_a^b \Bigg<-\left(\dd{}{t}+\bOm \times \right) \left[ \inertia  \bOm +\sum_{i=0}^n m_i \left[\mathbf{s}_0+\bchi_i-\bchi_0\right] \times \left[\bOm \times \left[ \mathbf{s}_0 + \bchi_i-\bchi_0\right]+\dot{\bchi}_i \right] \right] \\ & \hphantom{=\int_a^b \Bigg<} +\sum_{i=0}^n m_i \left( \left[ \left( \mathbf{s}_0-\bchi_0\right) \times \bOm \right]  \times \left[ \bOm \times \bchi_i+\dot{\bchi}_i \right] + g \bGamma \times \bchi_i  \right) +r \bGamma \times \tilde \bGamma, \bSigma \Bigg> \mathrm{d} t.
\end{split}
\end{equation}
In the fourth equality, the definitions $\tilde \bGamma \equiv \Lambda^{-1} \mathbf{F}_\mathrm{e}$, $\widehat \bSigma \equiv \Lambda^{-1} \de \Lambda$, and $\bGamma \equiv \Lambda^{-1} \mathbf{e}_3$ are used. In the final equality, the simplification of $\de S$ calculated in \eqref{eq:rball_vars_p7} is used.

\subsection{Equations of Motion for the Rolling Ball} \label{ssec_ball_uncontrolled}
Having computed the variation of the action integral \revision{R1Q1}{and having enforced the rolling and variational constraints} according to Lagrange-d'Alembert's principle, the equations of motion for the rolling ball actuated by internal point masses are obtained now. In addition, the equations of motion for two important special cases, a ball with static internal structure and a ball with 1-d parameterized rails, are derived. 

\paragraph{Equations of Motion for the Rolling Ball Actuated by Internal Point Masses}
Insisting that the variation $\de S_1$ of the action integral in \eqref{eq_dS1} is zero for all variations $\bSigma$ (i.e. completing the application of Lagrange-d'Alembert principle's by letting $0=\de S_1$) and using the identities \rem{${\dot {\mathbf{s}}}_0 = \bOm \times \left( \bchi_0-\mathbf{s}_0 \right)$ \eqref{eq:dot_s0} and $r \bGamma = \mathbf{s}_0-\bchi_0$} \revision{R1Q4}{\eqref{eq:dot_s0} and \eqref{eq:s0}}, the following equations of motion are obtained:
\begin{equation}
\begin{split} \label{uncon_ball_eqns_start}
\mathbf{0} &= \inertia \dot  \bOm+\bOm \times \inertia \bOm-r \bGamma \times \tilde \bGamma-g \sum_{i=0}^n m_i \bGamma \times \bchi_i+\sum_{i=0}^n m_i \left[r \bGamma \times \bOm+{\dot \bchi}_i \right] \times \left[\bOm \times \left[r \bGamma+\bchi_i \right]+{\dot \bchi}_i \right] \\
& \hphantom{=} + \sum_{i=0}^n m_i \left[r \bGamma+ \bchi_i \right] \times \left\{ \dot \bOm \times \left[r \bGamma+ \bchi_i \right]+ \bOm \times \left[r \bGamma \times \bOm+{\dot \bchi}_i \right] +{\ddot \bchi}_i \right\}  \\
& \hphantom{=}+\sum_{i=0}^n m_i \bOm \times \left\{ \left[r \bGamma+ \bchi_i \right] \times \left[\bOm \times \left[r \bGamma+\bchi_i \right] +{\dot \bchi}_i \right] \right\}
-\sum_{i=0}^n m_i \left[r \bGamma \times \bOm \right] \times \left[\bOm \times \bchi_i+{\dot \bchi}_i \right].
\end{split}
\end{equation}
As shown in Appendix~\ref{app_rball}, \eqref{uncon_ball_eqns_start} simplifies considerably to
\begin{equation}
\begin{split} \label{uncon_ball_eqns_explicit}
\dot \bOm &= \left[\sum_{i=0}^n m_i \widehat{\mathbf{s}_i}^2  -\inertia \right]^{-1}  \left[\bOm \times \inertia \bOm+r \tilde \bGamma \times \bGamma+ \sum_{i=0}^n m_i \mathbf{s}_i \times  \left\{ g \bGamma+ \bOm \times \left(\bOm \times \bchi_i +2 {\dot \bchi}_i \right) +{\ddot \bchi}_i \right\}  \right], \\
\dot \bGamma &= \bGamma \times \bOm,
\end{split}
\end{equation}
subject to the definitions $\mathbf{s}_i \equiv r \bGamma +\bchi_i$ for $0\le i\le n$, $\bOm \equiv \left( \Lambda^{-1} \dot \Lambda \right)^\vee$, $\bGamma \equiv \Lambda^{-1} \mathbf{e}_3$, and $\tilde \bGamma \equiv \Lambda^{-1} \mathbf{F}_\mathrm{e}$. \revision{R2Q2}{The trajectory of the spatial $\mathbf{e}_1$- and $\mathbf{e}_2$-components $\bz$ of the ball's GC and CP is obtained by replacing the second ordinary differential equation (ODE), $\dot \bGamma = \bGamma \times \bOm$, in \eqref{uncon_ball_eqns_explicit} with a pair of ODEs giving the evolution of $\Lambda$ and $\bz$. By the rolling constraint applied to the ball's GC, ${\dot {\mathbf{z}}}_\mathrm{GC} = \dot \Lambda r\bGamma =\Lambda \widehat{\bOm} r\bGamma = \Lambda \left[ \bOm \times r \bGamma \right] = \Lambda \bOm \times r \mathbf{e}_3$, and since $\bz \equiv \left(\mathbf{z}_\mathrm{GC} \right)_{12}$, the full equations of motion for the rolling ball are
\begin{equation}
\begin{split} \label{uncon_ball_eqns_explicit_full}
\dot \bOm &= \left[\sum_{i=0}^n m_i \widehat{\mathbf{s}_i}^2  -\inertia \right]^{-1}  \left[\bOm \times \inertia \bOm+r \tilde \bGamma \times \bGamma+ \sum_{i=0}^n m_i \mathbf{s}_i \times  \left\{ g \bGamma+ \bOm \times \left(\bOm \times \bchi_i +2 {\dot \bchi}_i \right) +{\ddot \bchi}_i \right\}  \right], \\
\dot \Lambda &= \Lambda \widehat{\bOm}, \\
\dot \bz &= \left( \Lambda \bOm \times r \mathbf{e}_3 \right)_{12}.
\end{split}
\end{equation} 
For $\mathbf{v} = \begin{bmatrix} v_1 & v_2 & v_3  \end{bmatrix}^\mathsf{T} \in \mathbb{R}^3$, $\mathbf{v}_{12}$ is the projected vector consisting of the first two components of $\mathbf{v}$ so that
\begin{equation}
\mathbf{v}_{12} = \begin{bmatrix}v_1 & v_2 \end{bmatrix}^\mathsf{T} \in \mathbb{R}^2.
\end{equation} 
These equations of motion \eqref{uncon_ball_eqns_explicit_full} for the rolling ball actuated by internal point masses are new and have not appeared previously in the literature, as far as we know.  }

\begin{remark}[On the parameterization of $\bchi_i$ and the final equations of motion]
{\rm 
Note that in the derivation of \revision{R2Q2}{\eqref{uncon_ball_eqns_explicit_full}}, we have not assumed any parameterization of the  mass trajectories $\bchi_i$: these equations are valid for arbitrary trajectories characterizing the motion of the masses. In what follows, we will explicitly assume that each $\bchi_i$ can be computed from one scalar parameter $\theta_i$, which occurs when the masses are moving along fixed  1-d trajectories in the ball's frame.  For example, this case can be realized when the masses are spun by a rotor on a lever of fixed length or when the masses move along rails fixed in the ball's frame, which is the case we consider below. 

One could alternatively consider the case where each $\bchi_i$ is parameterized by a set of parameters $\theta_{i,j}$, $j =1, 2, \ldots, J_i$. This can occur, for example, if the rotor spinning the lever in the ball in the example above can itself move or if the length of the lever can change. While these examples are interesting, their engineering implementations  are not readily apparent. In addition, in our opinion, taking $\bchi_i$ dependent on  multiple  parameters  introduces additional complexity into the equations  of motion without  enhancing mathematical understanding. We shall thus focus on the case when each $\bchi_i$ can be defined uniquely  by only one scalar parameter $\theta_i$. 
}
\end{remark}

\paragraph{ Equations of Motion for the Rolling Ball with Static Internal Structure}
A special case of \revision{R2Q2}{\eqref{uncon_ball_eqns_explicit_full}} gives the equations of motion for a rolling ball with static internal structure.
By fixing all the  point masses (i.e. making $\bchi_i$ constant for all $1 \le i \le n$, so that ${\dot \bchi}_i = {\ddot \bchi}_i = \mathbf{0}$)  or equivalently by setting the number of point masses $n$ to 0, \revision{R2Q2}{\eqref{uncon_ball_eqns_explicit_full}} gives the equations of motion for a rolling ball with static internal structure:
\rem{%%%BEGIN REM 
\begin{equation}
\begin{split} \label{uncon_ball_eqns_static}
\dot \bOm &= \left[\sum_{i=0}^n m_i \widehat{\mathbf{s}_i}^2  -\inertia \right]^{-1}  \left[\bOm \times \inertia \bOm+r \tilde \bGamma \times \bGamma+ \sum_{i=0}^n m_i \mathbf{s}_i \times  \left\{ g \bGamma+ \bOm \times \left(\bOm \times \bchi_i \right) \right\}  \right], \\
\dot \bGamma &= \bGamma \times \bOm.
\end{split}
\end{equation}
Alternatively, the  equations of motion for a rolling ball with static internal structure may be obtained by setting the number of point masses $n$ to 0 in \eqref{uncon_ball_eqns_explicit}:
} %%%END REM 
\revision{R2Q2}{
\begin{equation}
\begin{split} \label{uncon_ball_eqns_static2}
\dot \bOm &= \left[ m_0 \widehat{\mathbf{s}_0}^2  -\inertia \right]^{-1} \left[\bOm \times \inertia \bOm+r \tilde \bGamma \times \bGamma+  m_0 \mathbf{s}_0 \times  \left\{ g \bGamma+ \bOm \times \left(\bOm \times \bchi_0 \right) \right\}  \right], \\
\rem{\dot \bGamma &= \bGamma \times \bOm.}
\dot \Lambda &= \Lambda \widehat{\bOm}, \\
\dot \bz &= \left( \Lambda \bOm \times r \mathbf{e}_3 \right)_{12}.
\end{split}
\end{equation} }

\paragraph{Equations of Motion for the Rolling Ball Assuming 1-d Parameterizations of the Mass Trajectories}
For $1 \le i \le n$, assume now that the trajectory  $\bchi_i$ of the $i^\mathrm{th}$  point mass is required to move along a 1-d rail, like a circular hoop.
Moreover, for $1 \le i \le n$, assume that the $i^\mathrm{th}$  rail is parameterized by a 1-d parameter $\theta_i$, so that the trajectory $\bzeta_i$ of the $i^\mathrm{th}$  rail, in the body frame translated to the ball's geometric center, as a function of $\theta_i$ is $\bzeta_i(\theta_i)$. Thus, the trajectory of the $i^\mathrm{th}$   point mass as a function of time $t$ is $\bchi_i(t) \equiv \bzeta_i(\theta_i (t))$, $1 \le i \le n$. Refer to Figure~\ref{fig:detailed_1dparam_rolling_ball} for an illustration. To make notation consistent, define $\bzeta_0(\theta_0) \equiv \bchi_0$, so that the constant (time-independent) vector  $\bchi_0 = \bchi_0(t) \equiv \bzeta_0(\theta_0 (t))$ for any scalar-valued, time-varying function $\theta_0(t)$. By the chain rule and using the notation ${\textvisiblespace}^\cdot$ to denote differentiation with respect to time $t$ and $\bzeta_i^\prime$ to denote differentiation of $\bzeta_i$ with respect to $\theta_i$, for $0 \le i \le n$,
\begin{equation} \label{eq_1d_param}
\begin{split}
\bchi_i(t) &\equiv \bzeta_i(\theta_i (t)) = \bzeta_i, \\
\dot \bchi_i(t) &= \dd{\bzeta_i }{\theta_i}(\theta_i (t)) \dot \theta_i(t) = \bzeta_i^{\prime} (\theta_i (t))  \dot \theta_i(t) =  \bzeta_i^{\prime} \dot \theta_i  =  \dot \theta_i \bzeta_i^{\prime}, \\
\ddot \bchi_i(t) &=\frac{d^2 \bzeta_i }{d \theta_i^2}(\theta_i (t)) \dot \theta_i^2(t) + \dd{\bzeta_i }{\theta_i}(\theta_i (t)) \ddot \theta_i(t) \\
& =\bzeta_i^{\dprime}(\theta_i (t)) \dot \theta_i^2(t) + \bzeta_i^{\prime} (\theta_i (t))  \ddot \theta_i(t)=\bzeta_i^{\dprime} \dot \theta_i^2 + \bzeta_i^{\prime} \ddot \theta_i = \dot \theta_i^2 \bzeta_i^{\dprime} + \ddot \theta_i \bzeta_i^{\prime}. 
\end{split}
\end{equation}
By plugging the formulas for $\bchi_i$, $\dot \bchi_i$, and $\ddot \bchi_i$ given in \eqref{eq_1d_param} into \revision{R2Q2}{\eqref{uncon_ball_eqns_explicit_full}}, the equations of motion become
\revision{R2Q2}{
\begin{equation}
\begin{split} \label{uncon_ball_eqns_explicit_1d}
\dot \bOm &= \left[\sum_{i=0}^n m_i \widehat{\mathbf{s}_i}^2  -\inertia \right]^{-1}  \Bigg[\bOm \times \inertia \bOm+r \tilde \bGamma \times \bGamma 
\\ &\hphantom{=} + \sum_{i=0}^n m_i \mathbf{s}_i \times  \left\{ g \bGamma+ \bOm \times \left(\bOm \times \bzeta_i +2 \dot \theta_i \bzeta_i^{\prime} \right) + \dot \theta_i^2 \bzeta_i^{\dprime} + \ddot \theta_i  \bzeta_i^{\prime} \right\}  \Bigg], \\
\rem{ \dot \bGamma &= \bGamma \times \bOm,}
\dot \Lambda &= \Lambda \widehat{\bOm}, \\
\dot \bz &= \left( \Lambda \bOm \times r \mathbf{e}_3 \right)_{12},
\end{split}
\end{equation} }
where with this new notation, $\mathbf{s}_i \equiv r \bGamma +\bchi_i = r \bGamma +\bzeta_i$ for $0\le i\le n$.

\begin{figure}[h]
	\centering
	\includegraphics[width=0.5\linewidth]{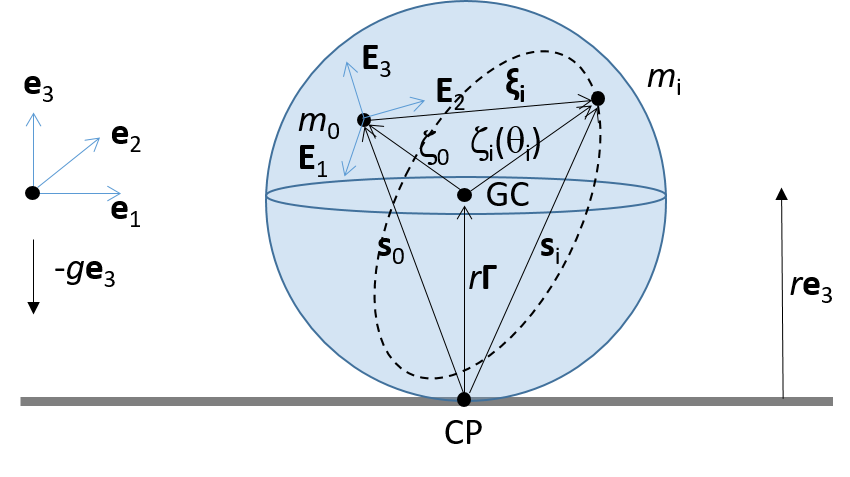}
	\caption{Each  point mass, denoted by $m_i$, $1 \le i \le n$, moves along a rail fixed inside the ball depicted here by the dashed hoop. The trajectory of the rail is denoted by $\bzeta_i$ and is parameterized by $\theta_i$. }
	\label{fig:detailed_1dparam_rolling_ball}
\end{figure}

\subsection{Equation of Motion for the Rolling Disk} \label{ssec_disk_uncontrolled}
 Let us now  demonstrate how to reduce the general equations of motion \eqref{uncon_ball_eqns_explicit_1d} for the  rolling ball  when its motion is purely planar, which  is the case of a rolling disk. \revision{R1Q2}{Unlike the rolling ball, which is a nonholonomic system, the rolling disk is a holonomic system.} While this particular 2-d \revision{R1Q2}{holonomic} case has  limited practicality, it is still useful to consider since  its equation of motion can be derived  via both variational methods and Newton's second law,   thereby providing additional validation of  \eqref{uncon_ball_eqns_explicit_1d}. In order to perform this two-dimensional reduction,  suppose that $m_0$'s inertia is such that one of $m_0$'s principal axes, say the one labeled $\mathbf{E}_2$, is orthogonal to the plane containing the GC and CM. Also assume that all the  point masses move along 1-d  rails which lie in the plane containing the GC and CM. Moreover, suppose that the ball is oriented initially so that the plane containing the GC and CM coincides with the $\mathbf{e}_1$-$\mathbf{e}_3$ plane and that the external force $\mathbf{F}_\mathrm{e}$ acts in the $\mathbf{e}_1$-$\mathbf{e}_3$ plane. Then for all time, the ball will remain oriented so that the plane containing the GC and CM coincides with the $\mathbf{e}_1$-$\mathbf{e}_3$ plane and the ball will only move in the $\mathbf{e}_1$-$\mathbf{e}_3$ plane, with the ball's rotation axis always parallel to $\mathbf{e}_2$. Note that the dynamics of this system are equivalent to that of the Chaplygin disk \cite{Ho2011_pII}, equipped with point masses, rolling in the $\mathbf{e}_1$-$\mathbf{e}_3$ plane, and where the Chaplygin disk (minus the  point masses) has polar moment of inertia $d_2$. Therefore, henceforth, this particular ball with this special inertia, orientation,  and placement of the rails and point masses, may be referred to as the disk or the rolling disk. Figure~\ref{fig:detailed_rolling_disk} depicts the rolling disk.
\begin{figure}[h]
	\centering
	\includegraphics[width=0.6\linewidth]{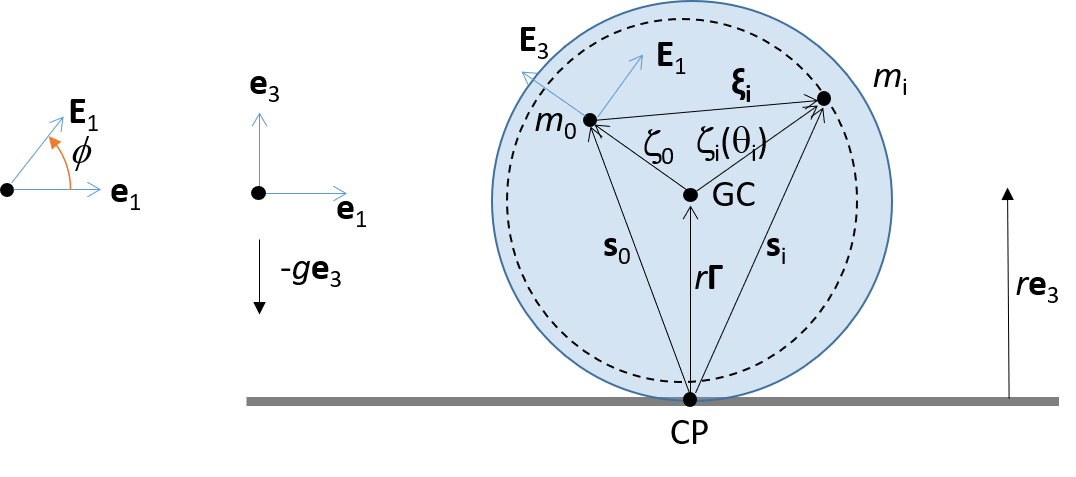}
	\caption{A disk of radius $r$ and mass $m_0$ rolls without slipping in the $\mathbf{e}_1$-$\mathbf{e}_3$ plane. $\mathbf{e}_2$ and $\mathbf{E}_2$ are directed into the page and are omitted from the figure. The disk's center of mass is denoted by $m_0$. The disk's motion is actuated by $n$ point masses, each of mass $m_i$, $1 \le i \le n$, that move along rails fixed inside the disk. The point mass depicted here by $m_i$ moves along a circular hoop in the disk that is not centered on the disk's geometric center (GC). The disk's orientation is determined by $\phi$, the angle measured counterclockwise from $\mathbf{e}_1$ to $\mathbf{E}_1$.  }
	\label{fig:detailed_rolling_disk}
\end{figure}

Let $\phi$ denote the angle between $\mathbf{e}_1$ and $\mathbf{E}_1$, measured counterclockwise from $\mathbf{e}_1$ to $\mathbf{E}_1$. Thus, if $\dot \phi > 0$, the disk rolls in the $-\mathbf{e}_1$ direction and $\bOm$ has the same direction as $-\mathbf{e}_2$, and  if $\dot \phi < 0$, the disk rolls in the $\mathbf{e}_1$ direction and $\bOm$ has the same direction as $\mathbf{e}_2$. Before constructing the equations of motion for the rolling disk using \eqref{uncon_ball_eqns_explicit_1d}, some intermediate calculations must be performed.
\rem{ 
\begin{equation}
\Lambda = \begin{bmatrix} \cos \phi & 0 & - \sin \phi \\ 0 & 1 & 0 \\ \sin \phi & 0 & \cos \phi \end{bmatrix}, \quad  \Lambda^\mathsf{T} = \begin{bmatrix} \cos \phi & 0 & \sin \phi \\ 0 & 1 & 0 \\ -\sin \phi & 0 & \cos \phi \end{bmatrix}.
\end{equation}
\begin{equation}
\hat{\bOm} = \Lambda^\mathsf{T} \dot \Lambda = \begin{bmatrix} \cos \phi & 0 & \sin \phi \\ 0 & 1 & 0 \\ -\sin \phi & 0 & \cos \phi \end{bmatrix} \begin{bmatrix} -\sin \phi & 0 & -\cos \phi \\ 0 & 0 & 0 \\ \cos \phi & 0 & -\sin \phi \end{bmatrix} \dot \phi = \begin{bmatrix} 0 & 0 & -1 \\ 0 & 0 & 0 \\ 1 & 0 & 0 \end{bmatrix} \dot \phi.
\end{equation}
\begin{equation}
\bOm = \begin{bmatrix} 0 \\ -1 \\ 0 \end{bmatrix} \dot \phi = -\dot \phi \begin{bmatrix} 0 \\ 1 \\ 0 \end{bmatrix} = - \dot \phi  \mathbf{e}_2, \quad \inertia \bOm = - d_2 \dot \phi \mathbf{e}_2, \quad \bOm \times \inertia \bOm = d_2 {\dot \phi}^2 \mathbf{e}_2 \times \mathbf{e}_2 =  \mathbf{0}.
\end{equation}
\begin{equation}
\bGam = \Lambda^\mathsf{T} \mathbf{e}_3 = \begin{bmatrix} \sin \phi \\ 0 \\ \cos \phi \end{bmatrix}, \quad \tilde {\bGam} = \Lambda^\mathsf{T} \mathbf{F}_\mathrm{e}= \begin{bmatrix} \cos \phi F_{\mathrm{e},1} + \sin \phi F_{\mathrm{e},3} \\ 0 \\ -\sin \phi F_{\mathrm{e},1} + \cos \phi F_{\mathrm{e},3} \end{bmatrix}.
\end{equation}
\begin{equation}
\begin{split}
r  \tilde \bGam \times \bGam &= r \begin{bmatrix} \cos \phi F_{\mathrm{e},1} + \sin \phi F_{\mathrm{e},3} \\ 0 \\ -\sin \phi F_{\mathrm{e},1} + \cos \phi F_{\mathrm{e},3} \end{bmatrix} \times  \begin{bmatrix} \sin \phi \\ 0 \\ \cos \phi \end{bmatrix} \\
&= r \left\{ \left( -\sin \phi F_{\mathrm{e},1} + \cos \phi F_{\mathrm{e},3} \right) \sin \phi - \left( \cos \phi F_{\mathrm{e},1} + \sin \phi F_{\mathrm{e},3}  \right) \cos \phi \right\} \mathbf{e}_2 \\
&= -r F_{\mathrm{e},1} \mathbf{e}_2.
\end{split}
\end{equation}
\begin{equation}
\dot \bOm = - \ddot \phi  \mathbf{e}_2, \quad \dot \bGam = \begin{bmatrix} \cos \phi \\ 0 \\ -\sin \phi \end{bmatrix} \dot \phi, \quad \bGam \times \bOm =  -\dot \phi \begin{bmatrix} \sin \phi \\ 0 \\ \cos \phi \end{bmatrix} \times \begin{bmatrix} 0 \\ 1 \\ 0 \end{bmatrix} = -\dot \phi \begin{bmatrix} -\cos \phi \\ 0 \\ \sin \phi \end{bmatrix} = \dot \bGam.
\end{equation}
Thus, the second equation, $\dot \bGamma = \bGamma \times \bOm$, in \eqref{uncon_ball_eqns_explicit_1d} gives no information about the dynamics and may be ignored for the disk.
\begin{equation}
\bzeta_i = \begin{bmatrix} \zeta_{i,1} \\ 0 \\ \zeta_{i,3} \end{bmatrix}, \quad \bzeta_i^{\prime} = \begin{bmatrix} \zeta_{i,1}^{\prime} \\ 0 \\ \zeta_{i,3}^{\prime} \end{bmatrix}, \quad \bzeta_i^{\dprime} = \begin{bmatrix} \zeta_{i,1}^{\dprime} \\ 0 \\ \zeta_{i,3}^{\dprime} \end{bmatrix}.
\end{equation}
\begin{equation}
\mathbf{s}_i = r \bGam + \bzeta_i = \begin{bmatrix} r \sin \phi + \zeta_{i,1} \\ 0 \\ r \cos \phi+ \zeta_{i,3} \end{bmatrix}.
\end{equation}
%\begin{multline} \label{eq_s_i_squared}
\begin{equation} \label{eq_s_i_squared}
\begin{split} 
\hat{\mathbf{s}}_i^2 = \mathbf{s}_i \mathbf{s}_i  = \hphantom{blah blah blah blah blah blah blah blah blah blah blah blah blah blah blah blah blah blah blah blah bl} 
\\ \hspace{-3mm} \begin{bmatrix}
-\left( r \cos \phi+ \zeta_{i,3} \right)^2 & 0 & \left( r \sin \phi + \zeta_{i,1} \right) \left( r \cos \phi+ \zeta_{i,3} \right) \\
0 & -\left( r \sin \phi + \zeta_{i,1} \right)^2-\left( r \cos \phi+ \zeta_{i,3} \right)^2  & 0 \\
\left( r \sin \phi + \zeta_{i,1} \right) \left( r \cos \phi+ \zeta_{i,3} \right) & 0 & -\left( r \sin \phi + \zeta_{i,1} \right)^2 
\end{bmatrix}.
\end{split}
\end{equation}
%\end{multline}
\begin{equation}
\bOm \times \bzeta_i = -\dot \phi \begin{bmatrix} 0 \\ 1 \\ 0 \end{bmatrix} \times \begin{bmatrix} \zeta_{i,1} \\ 0 \\ \zeta_{i,3} \end{bmatrix} = - \dot \phi \begin{bmatrix} \zeta_{i,3} \\ 0 \\ -\zeta_{i,1} \end{bmatrix}, \quad \bOm \times \bzeta_i + 2 \dot \theta_i \bzeta_i^{\prime}  = \begin{bmatrix} -\dot \phi\zeta_{i,3}+2 \dot \theta_i \zeta_{i,1}^{\prime} \\ 0 \\ \dot \phi \zeta_{i,1}+2 \dot \theta_i \zeta_{i,3}^{\prime} \end{bmatrix}.
\end{equation}
\begin{equation}
\bOm \times \left( \bOm \times \bzeta_i + 2 \dot \theta_i \bzeta_i^{\prime} \right)  = -\dot \phi \begin{bmatrix} 0 \\ 1 \\ 0 \end{bmatrix}  \times \begin{bmatrix} -\dot \phi\zeta_{i,3}+2 \dot \theta_i \zeta_{i,1}^{\prime} \\ 0 \\ \dot \phi \zeta_{i,1}+2 \dot \theta_i \zeta_{i,3}^{\prime} \end{bmatrix}
= -\dot \phi \begin{bmatrix} \dot \phi \zeta_{i,1}+2 \dot \theta_i \zeta_{i,3}^{\prime} \\ 0 \\ \dot \phi\zeta_{i,3}-2 \dot \theta_i \zeta_{i,1}^{\prime}  \end{bmatrix}.
\end{equation}
\begin{equation}
g \bGam + \bOm \times \left( \bOm \times \bzeta_i + 2 \dot \theta_i \bzeta_i^{\prime} \right) + \dot \theta_i^2 \bzeta_i^{\dprime} + \ddot \theta_i  \bzeta_i^{\prime} 
=  \begin{bmatrix} g \sin \phi -\dot \phi \left( \dot \phi \zeta_{i,1}+2 \dot \theta_i \zeta_{i,3}^{\prime} \right) + \dot \theta_i^2 \zeta_{i,1}^{\dprime} + \ddot \theta_i  \zeta_{i,1}^{\prime} \\ 0 \\ g \cos \phi -\dot \phi \left( \dot \phi\zeta_{i,3}-2 \dot \theta_i \zeta_{i,1}^{\prime} \right) + \dot \theta_i^2 \zeta_{i,3}^{\dprime} + \ddot \theta_i  \zeta_{i,3}^{\prime}  \end{bmatrix}.
\end{equation}
\begin{equation}
\mathbf{s}_i \times \left\{ g \bGam + \bOm \times \left( \bOm \times \bzeta_i + 2 \dot \theta_i \bzeta_i^{\prime} \right) + \dot \theta_i^2 \bzeta_i^{\dprime} + \ddot \theta_i  \bzeta_i^{\prime} \right\} = K_i \mathbf{e}_2,
\end{equation}
where
\begin{equation}
\begin{split}
K_i &=  \left(r \cos \phi + \zeta_{i,3} \right) \left(g \sin \phi - \dot \phi \left(\dot \phi \zeta_{i,1}+2 {\dot \theta}_i \zeta_{i,3}^{\prime} \right) + {\dot \theta}_i^2 \zeta_{i,1}^{\dprime} + {\ddot \theta}_i \zeta_{i,1}^{\prime} \right)\\
&\hphantom{=} - \left(r \sin \phi + \zeta_{i,1} \right) \left(g \cos \phi - \dot \phi \left(\dot \phi \zeta_{i,3}-2 {\dot \theta}_i \zeta_{i,1}^{\prime} \right) + {\dot \theta}_i^2 \zeta_{i,3}^{\dprime} + {\ddot \theta}_i \zeta_{i,3}^{\prime} \right) \\
&=  \left(g+ r {\dot \phi}^2 \right) \left(\zeta_{i,3} \sin \phi - \zeta_{i,1} \cos \phi  \right)+
\left(r \cos \phi + \zeta_{i,3} \right) \left(- 2 \dot \phi {\dot \theta}_i \zeta_{i,3}^{\prime} + {\dot \theta}_i^2 \zeta_{i,1}^{\dprime} + {\ddot \theta}_i \zeta_{i,1}^{\prime} \right)\\
&\hphantom{=}   - \left(r \sin \phi + \zeta_{i,1} \right) \left( 2 \dot \phi {\dot \theta}_i \zeta_{i,1}^{\prime}+ {\dot \theta}_i^2 \zeta_{i,3}^{\dprime} + {\ddot \theta}_i \zeta_{i,3}^{\prime} \right). 
\end{split}
\end{equation} }
For the disk, the orientation matrix $\Lambda$ is parameterized by the angle of rotation $\phi$ about the axis $\mathbf{e}_2$: 
\begin{equation}
\Lambda = \begin{bmatrix} \cos \phi & 0 & - \sin \phi \\ 0 & 1 & 0 \\ \sin \phi & 0 & \cos \phi \end{bmatrix}.
\end{equation}
Since 
\begin{equation}
\bOm \equiv \left( \Lambda^{-1} \dot \Lambda \right)^\vee = \begin{bmatrix} 0 \\ -1 \\ 0 \end{bmatrix} \dot \phi = -\dot \phi \begin{bmatrix} 0 \\ 1 \\ 0 \end{bmatrix} = - \dot \phi  \mathbf{e}_2,
\end{equation}
the cross product terms vanish identically: 
\begin{equation} \label{eq_disk_int1}
\bOm \times \inertia \bOm = d_2 {\dot \phi}^2 \mathbf{e}_2 \times \mathbf{e}_2 =  \mathbf{0}.
\end{equation}
Since
\begin{equation}
\bGam = \Lambda^\mathsf{T} \mathbf{e}_3 = \begin{bmatrix} \sin \phi \\ 0 \\ \cos \phi \end{bmatrix} \quad \mathrm{and} \quad \tilde {\bGam} = \Lambda^\mathsf{T} \mathbf{F}_\mathrm{e}= \begin{bmatrix} \cos \phi F_{\mathrm{e},1} + \sin \phi F_{\mathrm{e},3} \\ 0 \\ -\sin \phi F_{\mathrm{e},1} + \cos \phi F_{\mathrm{e},3} \end{bmatrix},
\end{equation}
there is an explicit expression for the gravity torque given by 
\begin{equation} \label{eq_disk_int2}
\begin{split}
r  \tilde \bGam \times \bGam %&= 
%r \begin{bmatrix} \cos \phi F_{\mathrm{e},1} + \sin \phi F_{\mathrm{e},3} \\ 0 \\ -\sin \phi F_{\mathrm{e},1} + \cos \phi F_{\mathrm{e},3} \end{bmatrix} \times  \begin{bmatrix} \sin \phi \\ 0 \\ \cos \phi \end{bmatrix} \\
&= r \left\{ \left( -\sin \phi F_{\mathrm{e},1} + \cos \phi F_{\mathrm{e},3} \right) \sin \phi - \left( \cos \phi F_{\mathrm{e},1} + \sin \phi F_{\mathrm{e},3}  \right) \cos \phi \right\} \mathbf{e}_2 = -r F_{\mathrm{e},1} \mathbf{e}_2.
\end{split}
\end{equation}
For the disk, note that
\begin{equation}
\bzeta_i = \begin{bmatrix} \zeta_{i,1} \\ 0 \\ \zeta_{i,3} \end{bmatrix}, \quad \bzeta_i^{\prime} = \begin{bmatrix} \zeta_{i,1}^{\prime} \\ 0 \\ \zeta_{i,3}^{\prime} \end{bmatrix}, \quad \mathrm{and} \quad \bzeta_i^{\dprime} = \begin{bmatrix} \zeta_{i,1}^{\dprime} \\ 0 \\ \zeta_{i,3}^{\dprime} \end{bmatrix}.
\end{equation}
Since
\begin{equation}
\mathbf{s}_i = r \bGam + \bzeta_i = \begin{bmatrix} r \sin \phi + \zeta_{i,1} \\ 0 \\ r \cos \phi+ \zeta_{i,3} \end{bmatrix},
\end{equation}
%\begin{multline} \label{eq_s_i_squared}
\begin{equation} \label{eq_s_i_squared}
\begin{split} 
\widehat{\mathbf{s}_i}^2 = \widehat{\mathbf{s}_i} \widehat{\mathbf{s}_i}  = \hphantom{blah blah blah blah blah blah blah blah blah blah blah blah blah blah blah blah blah blah blah blah bl} 
\\ \hspace{-3mm} \begin{bmatrix}
-\left( r \cos \phi+ \zeta_{i,3} \right)^2 & 0 & \left( r \sin \phi + \zeta_{i,1} \right) \left( r \cos \phi+ \zeta_{i,3} \right) \\
0 & -\left( r \sin \phi + \zeta_{i,1} \right)^2-\left( r \cos \phi+ \zeta_{i,3} \right)^2  & 0 \\
\left( r \sin \phi + \zeta_{i,1} \right) \left( r \cos \phi+ \zeta_{i,3} \right) & 0 & -\left( r \sin \phi + \zeta_{i,1} \right)^2 
\end{bmatrix}.
\end{split}
\end{equation}
%\end{multline}
A calculation shows that
\begin{equation} \label{eq_disk_int3}
\mathbf{s}_i \times \left\{ g \bGam + \bOm \times \left( \bOm \times \bzeta_i + 2 \dot \theta_i \bzeta_i^{\prime} \right) + \dot \theta_i^2 \bzeta_i^{\dprime} + \ddot \theta_i  \bzeta_i^{\prime} \right\} = K_i \mathbf{e}_2,
\end{equation}
where
\rem{\begin{equation} \label{eq_K_i}
\begin{split}
K_i &=  \left(r \cos \phi + \zeta_{i,3} \right) \left(g \sin \phi - \dot \phi \left(\dot \phi \zeta_{i,1}+2 {\dot \theta}_i \zeta_{i,3}^{\prime} \right) + {\dot \theta}_i^2 \zeta_{i,1}^{\dprime} + {\ddot \theta}_i \zeta_{i,1}^{\prime} \right)\\
&\hphantom{=} - \left(r \sin \phi + \zeta_{i,1} \right) \left(g \cos \phi - \dot \phi \left(\dot \phi \zeta_{i,3}-2 {\dot \theta}_i \zeta_{i,1}^{\prime} \right) + {\dot \theta}_i^2 \zeta_{i,3}^{\dprime} + {\ddot \theta}_i \zeta_{i,3}^{\prime} \right) \\
&=  \left(g+ r {\dot \phi}^2 \right) \left(\zeta_{i,3} \sin \phi - \zeta_{i,1} \cos \phi  \right)+
\left(r \cos \phi + \zeta_{i,3} \right) \left(- 2 \dot \phi {\dot \theta}_i \zeta_{i,3}^{\prime} + {\dot \theta}_i^2 \zeta_{i,1}^{\dprime} + {\ddot \theta}_i \zeta_{i,1}^{\prime} \right)\\
&\hphantom{=}   - \left(r \sin \phi + \zeta_{i,1} \right) \left( 2 \dot \phi {\dot \theta}_i \zeta_{i,1}^{\prime}+ {\dot \theta}_i^2 \zeta_{i,3}^{\dprime} + {\ddot \theta}_i \zeta_{i,3}^{\prime} \right). 
\end{split}
\end{equation}}
\begin{equation} \label{eq_K_i}
\begin{split}
K_i &\equiv \left(g+ r {\dot \phi}^2 \right) \left(\zeta_{i,3} \sin \phi - \zeta_{i,1} \cos \phi  \right)+
\left(r \cos \phi + \zeta_{i,3} \right) \left(- 2 \dot \phi {\dot \theta}_i \zeta_{i,3}^{\prime} + {\dot \theta}_i^2 \zeta_{i,1}^{\dprime} + {\ddot \theta}_i \zeta_{i,1}^{\prime} \right)\\
&\hphantom{\equiv}   - \left(r \sin \phi + \zeta_{i,1} \right) \left( 2 \dot \phi {\dot \theta}_i \zeta_{i,1}^{\prime}+ {\dot \theta}_i^2 \zeta_{i,3}^{\dprime} + {\ddot \theta}_i \zeta_{i,3}^{\prime} \right). 
\end{split}
\end{equation}
Plugging \eqref{eq_disk_int1}, \eqref{eq_disk_int2}, and \eqref{eq_disk_int3} into the first equation in \eqref{uncon_ball_eqns_explicit_1d} gives the equations of motion for the rolling disk as
\rem{\begin{equation} \label{eqmo_chap_disk_1}
\begin{split}
- \ddot \phi  \mathbf{e}_2 &=  \left[\sum_{i=0}^n m_i \widehat{\mathbf{s}_i}^2  -\inertia \right]^{-1}  \left[ -r F_{\mathrm{e},1} \mathbf{e}_2+ \sum_{i=0}^n m_i K_i \mathbf{e}_2 \right] \\
&= \left( -r F_{\mathrm{e},1}+ \sum_{i=0}^n m_i K_i \right) \left[\sum_{i=0}^n m_i \widehat{\mathbf{s}_i}^2  -\inertia \right]^{-1} \mathbf{e}_2.
\end{split}
\end{equation}}
\begin{equation} \label{eqmo_chap_disk_1}
- \ddot \phi  \mathbf{e}_2 =  \left[\sum_{i=0}^n m_i \widehat{\mathbf{s}_i}^2  -\inertia \right]^{-1}  \left[ -r F_{\mathrm{e},1} \mathbf{e}_2+ \sum_{i=0}^n m_i K_i \mathbf{e}_2 \right] 
= \left( -r F_{\mathrm{e},1}+ \sum_{i=0}^n m_i K_i \right) \left[\sum_{i=0}^n m_i \widehat{\mathbf{s}_i}^2  -\inertia \right]^{-1} \mathbf{e}_2.
\end{equation}
\rem { \begin{framed} 
Note that $\left[\sum_{i=0}^n m_i \widehat{\mathbf{s}_i}^2  -\inertia \right]^{-1} \mathbf{e}_2$ is just the middle column of the matrix inverse of $A = \sum_{i=0}^n m_i \widehat{\mathbf{s}_i}^2  -\inertia$. Denote the entries of $A$ by
\begin{equation}
A = \sum_{i=0}^n m_i \widehat{\mathbf{s}_i}^2  -\inertia = \begin{bmatrix} a_{11} & a_{12} & a_{13} \\ a_{21} & a_{22} & a_{23} \\ a_{31} & a_{32} & a_{33} \end{bmatrix}.
\end{equation}
Since $\inertia$ is diagonal and from \eqref{eq_s_i_squared}, $a_{12} = a_{21} = a_{23} = a_{32} = 0$, so that
\begin{equation}
A = \sum_{i=0}^n m_i \widehat{\mathbf{s}_i}^2  -\inertia = \begin{bmatrix} a_{11} & 0 & a_{13} \\ 0 & a_{22} & 0 \\ a_{31} & 0 & a_{33} \end{bmatrix}
\end{equation}
and the determinant of $A$ simplifies to 
\begin{equation}
\begin{split}
\det A &= a_{11} a_{22} a_{33}+a_{21} a_{32} a_{13}+a_{31} a_{12} a_{23}-a_{11} a_{32} a_{23}-a_{31} a_{22} a_{13}-a_{21} a_{12} a_{33} \\
&= a_{11} a_{22} a_{33} - a_{31} a_{22} a_{13} \\
&= a_{22} \left( a_{11} a_{33} - a_{31} a_{13} \right).
\end{split}
\end{equation}
Using the formula for the inverse of a $3 \times 3$ matrix, the middle column of the matrix inverse of $ A = \sum_{i=0}^n m_i \widehat{\mathbf{s}_i}^2  -\inertia$ is
\begin{equation} \label{eq_mid_col_inverse}
\begin{split}
\left[\sum_{i=0}^n m_i \widehat{\mathbf{s}_i}^2  -\inertia \right]^{-1} \mathbf{e}_2 &= A^{-1} \mathbf{e}_2
= \begin{bmatrix} a_{11} & 0 & a_{13} \\ 0 & a_{22} & 0 \\ a_{31} & 0 & a_{33} \end{bmatrix}^{-1} \mathbf{e}_2 
= \frac{1}{\det A} \begin{bmatrix} a_{13} a_{32} - a_{12} a_{33} \\ a_{11} a_{33} - a_{13} a_{31} \\ a_{12} a_{31} - a_{11} a_{32}  \end{bmatrix} \\
&= \frac{1}{a_{22} \left( a_{11} a_{33} - a_{31} a_{13} \right)} \begin{bmatrix} 0 \\ a_{11} a_{33} - a_{13} a_{31} \\ 0 \end{bmatrix} 
=  \frac{1}{a_{22}} \begin{bmatrix} 0 \\ 1 \\ 0 \end{bmatrix} 
= \frac{1}{a_{22}} \mathbf{e}_2.
\end{split} 
\end{equation}
Plugging \eqref{eq_mid_col_inverse} into \eqref{eqmo_chap_disk_1}, the equations of motion simplify to
\begin{equation} \label{eqmo_chap_disk_2}
- \ddot \phi  \mathbf{e}_2 = \left( -r F_{\mathrm{e},1}+ \sum_{i=0}^n m_i K_i \right) \left[\sum_{i=0}^n m_i \widehat{\mathbf{s}_i}^2  -\inertia \right]^{-1} \mathbf{e}_2 = \frac{1}{a_{22}} \left( -r F_{\mathrm{e},1}+ \sum_{i=0}^n m_i K_i \right) \mathbf{e}_2,
\end{equation}
which gives the scalar equation of motion
\begin{equation} \label{eqmo_chap_disk_3}
\ddot \phi = \frac{-1}{a_{22}} \left( -r F_{\mathrm{e},1}+ \sum_{i=0}^n m_i K_i \right).
\end{equation}
From \eqref{eq_s_i_squared}
\begin{equation} \label{eq_a_22}
a_{22} = \sum_{i=0}^n \left\{ m_i \left[-\left( r \sin \phi + \zeta_{i,1} \right)^2-\left( r \cos \phi+ \zeta_{i,3} \right)^2 \right] \right\}  -d_2.
\end{equation}
\end{framed} }
As shown in Appendix~\ref{app_rdisk}, \eqref{eqmo_chap_disk_1} simplifies to the scalar equation of motion for the rolling disk
\begin{equation} \label{eqmo_chap_disk_4}
\ddot \phi = \frac{ -r F_{\mathrm{e},1}+ \sum_{i=0}^n m_i K_i }{d_2+\sum_{i=0}^n m_i \left[\left( r \sin \phi + \zeta_{i,1} \right)^2+\left( r \cos \phi+ \zeta_{i,3} \right)^2 \right]} \equiv \kappa\left(t,\btheta,\dot \btheta,\phi,\dot \phi,\ddot \btheta\right),
\end{equation}
where $\kappa$ is a function that depends on time ($t$) through the possibly time-varying external force $F_{\mathrm{e},1}(t)$, on the  point mass parameterized positions ($\btheta$), velocities ($\dot \btheta$), and accelerations ($\ddot \btheta$), and on the disk's orientation angle ($\phi$) and its time derivative ($\dot \phi$). \revision{R2Q2}{The spatial $\mathbf{e}_1$-component $z$ of the disk's GC and CP is given by 
\begin{equation} \label{eq_disk_GC}
z = z_a-r \left(\phi-\phi_a \right),
\end{equation} 
where $z_a$ is the spatial $\mathbf{e}_1$-component of the disk's GC and CP at time $t=a$ and $\phi_a$ is the disk's angle at time $t=a$.
}

\paragraph{Verification of the Variational Equations Using Newtonian Mechanics for a Special Case of the Rolling Disk}

This paper relies on variational Lagrangian mechanics  for the derivation of the equations of motion as it is, in our opinion, much more efficient than Newtonian mechanics when applied to mechanical systems with complex internal structure. However, a special case of the rolling disk can also be analyzed using standard Newtonian mechanics, which is worthwhile to investigate in order to verify the correctness of our variational approach.  See \cite{ilin2017dynamics} for the  derivation of  the  equations  of motion in three dimensions  via Newtonian mechanics. 

\rem{ %%%BEGIN REM 
\begin{figure}[h] 
	\centering
	\subfloat[Sir Isaac Newton, 1642-1727 \cite{Newton_portrait}.]{\includegraphics[scale=.55]{GodfreyKneller-IsaacNewton-1689}\label{fig_Newton}}
	\hspace{20mm}
	\subfloat[Joseph-Louis Lagrange, 1736-1813 \cite{Lagrange_portrait}.]{
		\includegraphics[scale=.709]{Lagrange_portrait}\label{fig_Lagrange}}
	\caption{Portraits of Newton and Lagrange, progenitors of classical mechanics.} \label{fig_Newton_Lagrange}
\end{figure}
} %%%END REM 

\rem{ In order to verify the equations of motion obtained using variational principles (i.e. Lagrange-d'Alembert's principle with Euler-Poincar\'e's method),  it will be shown that they agree with the equations of motion obtained by Newton's method in one simple case where the CM and GC coincide. This calculation will also show the advantage of using variational principles as opposed to Newton's law to derive the equations of motion. }

Consider a disk of mass $m_0$ and radius $r$ whose CM and GC coincide. 
The moment of inertia of the disk computed with respect to the CM is $d_2$. The disk rolls without slipping along a horizontal surface in a uniform gravitational field of magnitude $g$. The disk is actuated by a single point mass of mass $m_1$ that moves along a circular trajectory of radius $r_1$, with $0<r_1<r$, centered on the disk's GC. The spatial $\mathbf{e}_1$-component $z$ of the disk's GC \revision{R2Q2}{and CP is given by \eqref{eq_disk_GC}.} \rem{is given by $z(t)$.} Since $z(t) = z_a-r \left(\phi(t)-\phi_a \right)$, $\dot z(t) =-r \dot \phi(t)$ and $\ddot z(t) =-r \ddot \phi(t)$. Since the CM and GC coincide, the body frame coincides with the body frame translated to the GC. The  point mass's trajectory in the body frame translated to the GC is
\begin{equation} 
\bzeta_1(t) = r_1 \begin{bmatrix} \cos{\theta_1(t)} \\ 0 \\ \sin {\theta_1(t)} \end{bmatrix}
\end{equation}
and in the spatial frame is
\begin{equation} 
\begin{split}
\mathbf{z}_1(t) =\begin{bmatrix} z(t) \\ 0 \\ 0 \end{bmatrix}+ \Lambda(t) \bzeta_1(t) &=\begin{bmatrix} z(t) \\ 0 \\ 0 \end{bmatrix}+  r_1 \begin{bmatrix} \cos \phi(t) & 0 & - \sin \phi(t) \\ 0 & 1 & 0 \\ \sin \phi(t) & 0 & \cos \phi(t) \end{bmatrix} \begin{bmatrix} \cos{\theta_1(t)} \\ 0 \\ \sin {\theta_1(t)} \end{bmatrix}\\&= \begin{bmatrix} z(t)+ r_1 \cos\left(\phi(t)+\theta_1(t)\right) \\ 0 \\ r_1 \sin \left(\phi(t)+\theta_1(t)\right) \end{bmatrix}.
\end{split}
\end{equation}
Observe that the axis of rotation passes through the CM and that the axis of rotation does not change direction. Thus, it is straightforward to determine the dynamics of this system via Newtonian mechanics. Newton's second law says that the sum of all external forces acting on the disk must equal $m_0 \ddot z \mathbf{e}_1=-m_0 r \ddot \phi \mathbf{e}_1$ and that the sum of all external torques acting on the disk about the disk's CM must equal $ -d_2 \ddot \phi \mathbf{e}_2$. The external forces acting on the disk are the force $-m_1\ddot {\mathbf{z}}_1-m_1 g \mathbf{e}_3$ exerted by the accelerating  point mass, the gravitational force $-m_0 g \mathbf{e}_3$ exerted at the CM, a horizontal static frictional force $-f_s \mathbf{e}_1$ exerted by the surface, a normal force $N \mathbf{e}_3$ exerted by the surface, and an external force $\mathbf{F}_\mathrm{e}$ exerted at the disk's GC. See Figure~\ref{fig:simple_rolling_disk} for the free body diagram depicting all the external forces acting on the disk.

\begin{figure}[h]
	\centering
	\includegraphics[width=0.5\linewidth]{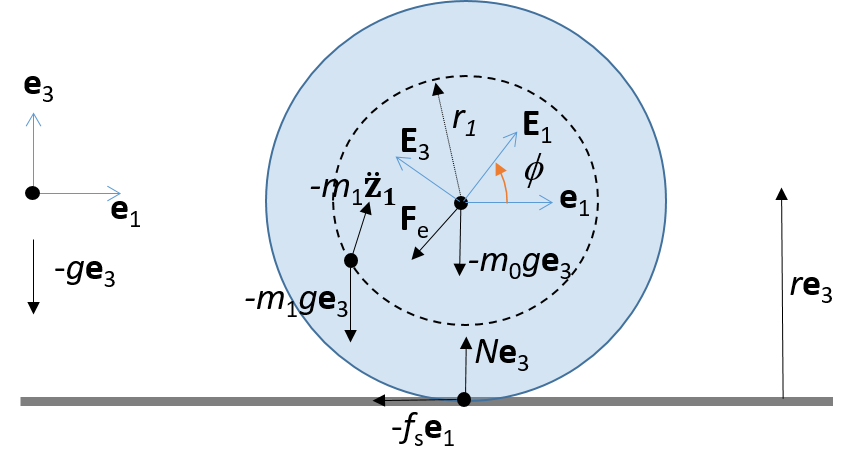}
	\caption{A disk actuated by a single point mass.}
	\label{fig:simple_rolling_disk}
\end{figure}

Application of Newton's second law to this system gives the following force and torque balance equations:
\begin{equation}
\begin{split}
\sum F_1 &=-f_\mathrm{s}+F_{\mathrm{e},1}-m_1 {\ddot {\mathrm{z}}}_{1,1} = m_0 \ddot z=-m_0 r \ddot \phi \implies f_\mathrm{s}= m_0 r \ddot \phi+F_{\mathrm{e},1}-m_1 {\ddot {\mathrm{z}}}_{1,1}\\
\sum F_3 &= N+F_{\mathrm{e},3}-m_0g-m_1g-m_1 {\ddot {\mathrm{z}}}_{1,3} = 0 \implies  N = -F_{\mathrm{e},3}+\left(m_0+m_1\right)g+m_1 {\ddot {\mathrm{z}}}_{1,3} \\
\sum \boldsymbol{\tau} &= r f_\mathrm{s} \mathbf{e}_2+m_1 g \left(\mathrm{z}_{1,1}-z \right) \mathbf{e}_2-m_1 \left(\mathbf{z}_1 - \begin{bmatrix} z \\ 0 \\ 0 \end{bmatrix}  \right) \times {\ddot {\mathbf{z}}}_1 = -d_2 \ddot \phi \mathbf{e}_2.
\end{split}
\end{equation}
Plugging the formula for the horizontal static friction force into the torque balance equation yields
\begin{equation}
r \left( m_0 r \ddot \phi+F_{\mathrm{e},1}-m_1 {\ddot {\mathrm{z}}}_{1,1} \right) \mathbf{e}_2+m_1 g \left(\mathrm{z}_{1,1}-z \right) \mathbf{e}_2-m_1 \left(\mathbf{z}_1 - \begin{bmatrix} z \\ 0 \\ 0 \end{bmatrix}  \right) \times {\ddot {\mathbf{z}}}_1 +d_2 \ddot \phi \mathbf{e}_2 = \mathbf{0},
\end{equation}
which simplifies to
\rem{ %%%BEGIN REM 
\begin{equation} \label{eq:srd_newton}
\begin{split}
\Bigg\{& d_2 \ddot \phi+ r \left( m_0 r \ddot \phi+F_{\mathrm{e},1}-m_1 \left[\ddot z- r_1 \cos\left(\phi+\theta_1\right) \left(\dot \phi+ \dot \theta_1 \right)^2-r_1 \sin \left(\phi+\theta_1\right) \left(\ddot \phi+ \ddot \theta_1 \right) \right] \right)\\
&+ m_1 g r_1 \cos\left(\phi+\theta_1\right) \Bigg\} \mathbf{e}_2\\
&-m_1 r_1 \begin{bmatrix} \cos\left(\phi+\theta_1\right) \\ 0 \\ \sin \left(\phi+\theta_1\right) \end{bmatrix} \times \begin{bmatrix} \ddot z- r_1 \cos\left(\phi+\theta_1\right) \left(\dot \phi+ \dot \theta_1 \right)^2-r_1 \sin \left(\phi+\theta_1\right) \left(\ddot \phi+ \ddot \theta_1 \right) \\ 0 \\ -r_1 \sin \left(\phi+\theta_1\right) \left(\dot \phi+ \dot \theta_1 \right)^2+r_1 \cos\left(\phi+\theta_1\right) \left(\ddot \phi+ \ddot \theta_1 \right) \end{bmatrix}  = 0.
\end{split}
\end{equation}
Equation \eqref{eq:srd_newton} is equivalent to 
\begin{equation} 
\begin{split}
&d_2 \ddot \phi+ r \left( m_0 r \ddot \phi+F_{\mathrm{e},1}-m_1 \left[-r \ddot \phi- r_1 \cos\left(\phi+\theta_1\right) \left(\dot \phi+ \dot \theta_1 \right)^2-r_1 \sin \left(\phi+\theta_1\right) \left(\ddot \phi+ \ddot \theta_1 \right) \right] \right) \\
&+ m_1 g r_1 \cos\left(\phi+\theta_1\right)\\
&-m_1 r_1 \left\{\sin \left(\phi+\theta_1\right) \left[ -r \ddot \phi- r_1 \cos\left(\phi+\theta_1\right) \left(\dot \phi+ \dot \theta_1 \right)^2-r_1 \sin \left(\phi+\theta_1\right) \left(\ddot \phi+ \ddot \theta_1 \right) \right] \right. \\
& \left. - \cos\left(\phi+\theta_1\right) \left[ -r_1 \sin \left(\phi+\theta_1\right) \left(\dot \phi+ \dot \theta_1 \right)^2+r_1 \cos\left(\phi+\theta_1\right) \left(\ddot \phi+ \ddot \theta_1 \right) \right]  \right\} = 0,
\end{split}
\end{equation}
which simplifies to
\begin{equation} 
\begin{split}
&d_2 \ddot \phi+ r \left( m_0 r \ddot \phi+F_{\mathrm{e},1}+m_1 \left[r \ddot \phi+ r_1 \cos\left(\phi+\theta_1\right) \left(\dot \phi+ \dot \theta_1 \right)^2+r_1 \sin \left(\phi+\theta_1\right) \left(\ddot \phi+ \ddot \theta_1 \right) \right] \right) \\
&+ m_1 g r_1 \cos\left(\phi+\theta_1\right)+m_1 r_1 \left\{ r \ddot \phi \sin \left(\phi+\theta_1\right)+ r_1 \left(\ddot \phi+ \ddot \theta_1 \right)  \right\} = 0,
\end{split}
\end{equation}
which further simplifies to the following equation of motion
} %%%END REM 
\begin{equation} \label{eq:srd_newton2}
\ddot \phi = -\frac{r F_{\mathrm{e},1}+m_1 r_1\left[ \cos\left(\phi+\theta_1\right)\left\{r\left(\dot \phi+\dot \theta_1\right)^2+g \right\}+\left\{r_1+r\sin\left(\phi+\theta_1\right)\right\} \ddot \theta_1\right]}{d_2+(m_0+m_1)r^2+m_1r_1 \left[r_1+2r\sin\left(\phi+\theta_1\right)\right]}.
\end{equation}
Under all these assumptions for this particular rolling disk, a calculation shows that equation \eqref{eq:srd_newton2} coincides with equation \eqref{eqmo_chap_disk_4}, which was derived earlier by variational methods (i.e. Lagrangian mechanics).

\section{Numerical Simulations of the Dynamics of the Rolling Disk}
\label{sec_disk_numerical}
	To write the equations of motion for the rolling disk in the standard ODE form, the state of the system is defined as 
\begin{equation}
\bx \equiv \begin{bmatrix} \btheta \\ \dot \btheta \\ \phi \\ \dot \phi  \end{bmatrix},
\end{equation}
where $\btheta, \dot \btheta \in \mathbb{R}^n$ and $\phi, \dot \phi  \in \mathbb{R}$.
The ODE formulation of the rolling disk's system dynamics defined for $a \le t \le b$ is
\begin{equation} \label{eq_rdisk_dynamics}
\dot {\bx} = \begin{bmatrix} \dot \btheta \\ \ddot \btheta \\ \dot \phi \\ \ddot \phi  \end{bmatrix}  = \mathbf{f}\left(t,\bx,\bu\right) \equiv \begin{bmatrix} \dot \btheta \\ \bu  \\ \dot \phi \\ \kappa\left(t,\bx,\bu\right)  \end{bmatrix},
\end{equation}
where $\bu \colon \mathbb{R} \to \mathbb{R}^n$ is a prescribed function of $t$ such that $\bu(t)=\ddot \btheta(t) \in \mathbb{R}^n$ and $\kappa\left(t,\bx,\bu \right)$ is given in \eqref{eqmo_chap_disk_4}:
\begin{equation} \label{eq_kappa}
\kappa\left(t,\bx,\bu\right) \equiv \frac{ -r F_{\mathrm{e},1}+ \sum_{i=0}^n m_i K_i }{d_2+\sum_{i=0}^n m_i \left[\left( r \sin \phi + \zeta_{i,1} \right)^2+\left( r \cos \phi+ \zeta_{i,3} \right)^2 \right]},
\end{equation}
where $K_i$ is given by \eqref{eq_K_i}. In order to simulate the rolling disk's dynamics, \eqref{eq_rdisk_dynamics} must be integrated with  prescribed initial conditions at time $t=a$:
\begin{equation} \label{eq_disk_initial_conds}
\bx\left(a\right) = \begin{bmatrix} \btheta(a) \\ \dot \btheta(a)  \\ \phi(a) \\ \dot \phi(a)  \end{bmatrix} = \begin{bmatrix} \btheta_a \\ {\dot \btheta}_a \\ \phi_a \\  \unaryminus \frac{{\dot z}_a}{r}  \end{bmatrix} \equiv \bx_a.
\end{equation}
\eqref{eq_rdisk_dynamics} and \eqref{eq_disk_initial_conds} constitute an ODE IVP.
 For   the ODE systems considered here, one can  choose $a=0$ without loss of generality;  however, we shall let $a$ be  arbitrary to keep our discussion general and consistent with the notation used in the literature on the numerical solution of boundary value problems \cite{ascher1994numerical}. \revision{R2Q2}{Given $\phi$, the spatial $\mathbf{e}_1$-component $z$ of the disk's GC and CP may be obtained from \eqref{eq_disk_GC}.} 

Consider a rolling disk of mass $m_0=1$, radius $r=1$, polar moment of inertia $d_2=1$, and with the CM coinciding with the GC (i.e. $\bzeta_0=\mathbf{0}$). The disk contains $n=4$ internal point masses, each of mass $1$ so that $m_1=m_2=m_3=m_4=1$ and each located on its own concentric circle centered on the GC of radius $r_1=.9$, $r_2=.6\overline{3}$, $r_3=.3\overline{6}$, and $r_4=.1$, respectively, as shown in Figure~\ref{fig_disk_masses_rails_bf_gc}. For $1 \le i \le n$, the position of $m_i$ in the body frame centered on the GC is:
\begin{equation}
\bzeta_i\left(\theta_i\right) = r_i \begin{bmatrix} \cos \theta_i \\ 0 \\ \sin \theta_i  \end{bmatrix}.
\end{equation}
The disk's total system mass is $M=5$, and gravity is rescaled to be $g=1$. There is no external force acting on the disk's GC so that $F_{\mathrm{e},1} = 0$ in \eqref{eq_kappa}. This disk's dynamics are simulated with initial time $a=0$ and final time $b=20$, so that the simulation time interval is $\left[0,20\right]$. The parameterized acceleration of each internal point mass is a continuous approximation of a short duration unit amplitude step function:
\rem{\todo{VP: Are you trying to say that the acceleration is a linear function on an interval and zero otherwise? \\ SMR: The acceleration is constant at 1 for a short time, then decreases linearly from 1 to 0 for a short time, and then stays constant at 0 for the rest of time. $u_i = \ddot \theta_i$ is constructed to be continuous (instead of discontinuous) so that $\dot \btheta$ and $\btheta$ are differentiable, by the Fundamental Theorem of Calculus. The derivation of the equations of motion assumed that $\btheta$ and $\dot  \btheta$ are differentiable. Should I add in a plot of $u_i$? I didn't make the plot in the interest of keeping the paper short. The coefficient $\left(\unaryminus 1\right)^i$ in $u_i$ makes the simulation a little more interesting so that odd masses rotate clockwise and even masses rotate counterclockwise in the disk's body frame. \\ 
VP: Yes, please put a plot of $u_i(t)$ to explain. \\ SMR: I added in a plot of $|u_i(t)|$.  }} 
\begin{equation} \label{eq_prescribed_u_disk}
u_i(t) = \ddot \theta_i(t) = \left(\unaryminus 1\right)^i \left\{ \begin{array}{ll} 1, & 0 \le t \le .1, \\ \unaryminus 10t+2, & .1 \le t \le .2, \\ 0, & .2 \le t \le 20, \end{array} \right. \quad \mbox{for} \quad 1 \le i \le n.
\end{equation}
\begin{figure}[h] 
	\centering
	\includegraphics[scale=.5]{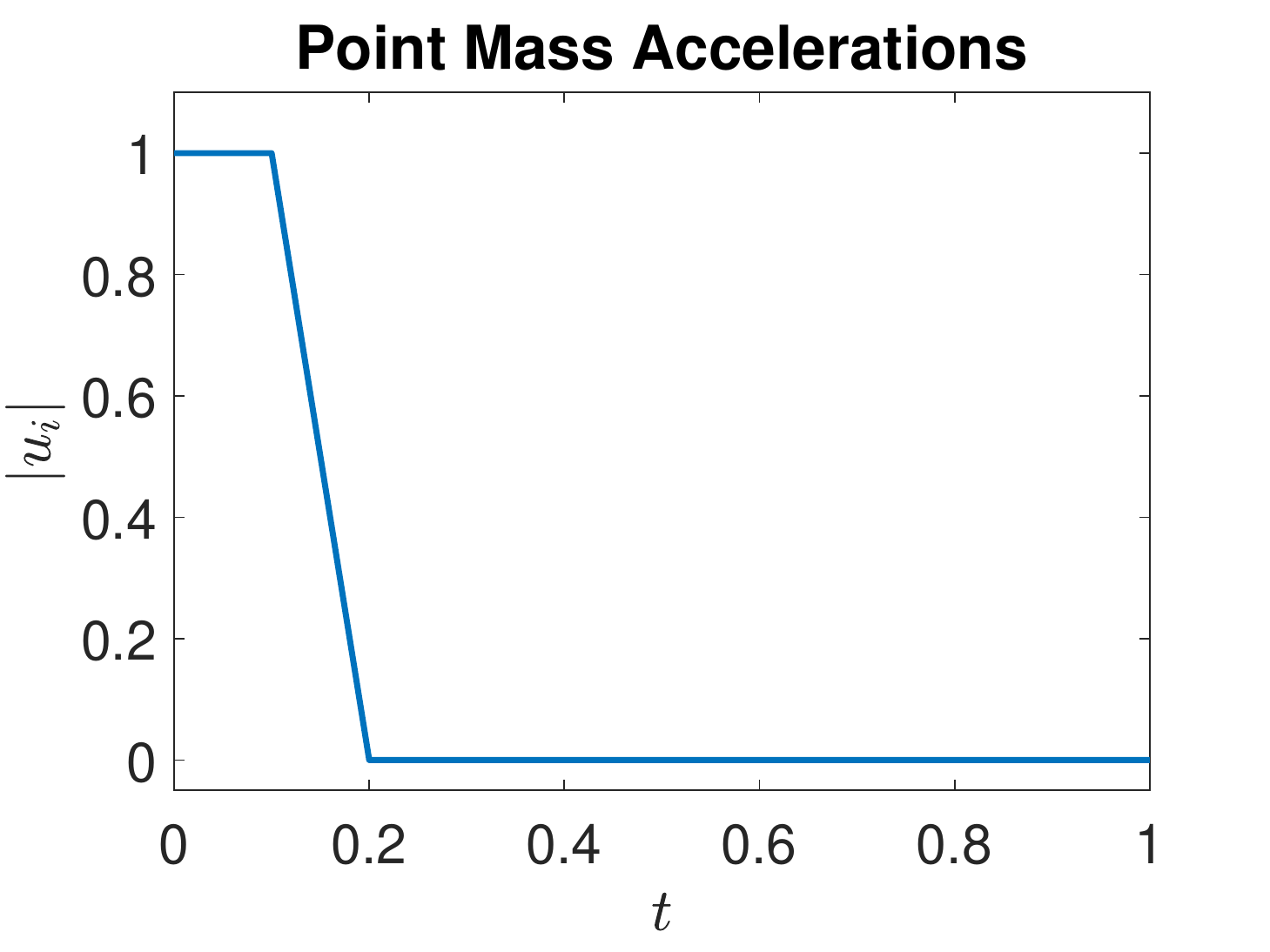}
	\caption{The magnitude of the parameterized acceleration, $u_i(t)=\ddot \theta_i(t)$, of each point mass, $1\le i \le n$.}
	\label{fig_abs_ui}
\end{figure}
The  magnitudes of the functions $u_i(t)$ are illustrated in Figure~\ref{fig_abs_ui}. For each $i$, the magnitude of the parameterized acceleration $u_i$ is chosen to be $1$ for the short time interval $0\le t \le 0.1$, then decreases linearly from $1$ to $0$ for the short time interval $0.1 \le t \le 0.2$, and finally stays constant at $0$ for the rest of time. The parameterized accelerations $u_i = \ddot \theta_i$ are constructed to be continuous (instead of discontinuous) so that $\dot \btheta$ and $\btheta$ are differentiable. We have used these parameterized accelerations since the derivation of the equations of motion \eqref{uncon_ball_eqns_explicit_1d} and \eqref{eqmo_chap_disk_4} assumed that $\btheta$ and $\dot  \btheta$ are differentiable.

The rolling disk's initial conditions are selected so that the disk starts at rest at the origin. Table~\ref{table_disk_ICs} shows parameter values used in the rolling disk's initial conditions \eqref{eq_disk_initial_conds}. Since the initial orientation of the disk is $\phi_a=0$ and since the initial configurations of the internal point masses are given by $\btheta_a=\begin{bmatrix} \unaryminus \frac{\pi}{2} & \unaryminus \frac{\pi}{2} & \unaryminus \frac{\pi}{2} & \unaryminus \frac{\pi}{2} \end{bmatrix}^\mathsf{T}$, all the internal point masses are initially located directly below the GC, so that the disk's total system CM is initially located below the GC. To ensure that the disk is initially at rest, ${\dot \btheta}_a = \begin{bmatrix} 0 & 0 & 0 & 0 \end{bmatrix}^\mathsf{T}$ and $\dot \phi_a = \unaryminus \frac{{\dot z}_a}{r} = 0$. To ensure that the disk's GC is initially located at the origin, $z_a = 0$. In summary, the rolling disk's initial conditions are
\begin{equation} \label{eq_disk_initial_conds_spec}
\bx_a = \begin{bmatrix} \unaryminus \frac{\pi}{2} & \unaryminus \frac{\pi}{2} & \unaryminus \frac{\pi}{2} & \unaryminus \frac{\pi}{2} & 0 & 0 & 0 & 0 & 0 & 0  \end{bmatrix}^\mathsf{T} .
\end{equation}

\begin{table}[h!]
	\centering 
	{ 
		\setlength{\extrarowheight}{1.5pt}
		\begin{tabular}{| c | c |} 
			\hline
			\textbf{Parameter} & \textbf{Value} \\ 
			\hline\hline 
			$\btheta_a$ & $\begin{bmatrix} \unaryminus \frac{\pi}{2} & \unaryminus \frac{\pi}{2} & \unaryminus \frac{\pi}{2} & \unaryminus \frac{\pi}{2} \end{bmatrix}^\mathsf{T}$  \\  
			\hline
			$\dot \btheta_a$ & $\begin{bmatrix} 0 & 0 & 0 & 0 \end{bmatrix}^\mathsf{T}$ \\ 
			\hline
			$\phi_a$ & $0$ \\
			\hline
			$\dot \phi_a$ & $0$ \\  
			\hline
			$z_a$ & $0$ \\
			\hline
			$\dot z_a$ & $0$ \\ 
			\hline
		\end{tabular} 
	}
	\caption{Initial condition parameter values for the rolling disk.}
	\label{table_disk_ICs}
\end{table}

The dynamics of this rolling disk are simulated by numerically integrating the ODE IVP \eqref{eq_rdisk_dynamics}, \eqref{eq_disk_initial_conds_spec} via \mcode{MATLAB} R2017b and Fortran  ODE-integration routines.  For ODE integrators, we have used the \mcode{MATLAB} R2017b routines \mcode{ode45}, \mcode{ode113}, \mcode{ode15s}, \mcode{ode23t}, and \mcode{ode23tb} and a \mcode{MATLAB} \mcode{MEX} wrapper of the Fortran routine radau5 \cite{hairer1996solving}, using the default input options except for the absolute and relative error tolerances and the Jacobian. The absolute and relative error tolerances supplied to the numerical integrators are both set to $1\mathrm{e}{-12}$. The Jacobian of $\mathbf{f}$ with respect to the state $\bx$, obtained via complex-step differentiation \cite{squire1998using,martins2001connection,martins2003complex}, is supplied to \mcode{ode15s}, \mcode{ode23t}, \mcode{ode23tb}, and radau5. Since excellent agreement was observed between all the numerical integrators, only the results obtained by numerically integrating the ODE IVP \eqref{eq_rdisk_dynamics}, \eqref{eq_disk_initial_conds_spec} with \mcode{ode45} are shown in Figure~\ref{fig_disk_ode_sims}. We shall also note that while all the numerical integrators yielded identical results, \mcode{ode113} completed the numerical integration in the shortest time.
\rem{\todo{VP: Do we need these details above you think? \\ SMR: I provided all the gory details to enable reproducibility of the simulation results. A reader of this version of the paper should be able to exactly reproduce our results, given all these details. I don't know how detailed we need to be for publication. \\
VP: I corrected above. What we can say, in principle, is this: we tried several routines, and they all gave the same results. Is there a time difference of one routine vs the other? Which one was the fastest? That is probably worth saying. \\ SMR: For the rolling disk simulation, the times are \mcode{ode45}: .23s, \mcode{ode113}: .06s, \mcode{ode15s}: .30s, \mcode{ode23t}: 5.01s, \mcode{ode23tb}: 7.18s, radau5: .27s. For the rolling ball simulation, the times are \mcode{ode45}: .25s, \mcode{ode113}: .10s, \mcode{ode15s}: .25s, \mcode{ode23t}: 2.19s, \mcode{ode23tb}: 3.72s, radau5: .41s. So for both simulations, \mcode{ode113} was the fastest.  \\ 
VP: OK, I added a sentence above. }}
\begin{figure}[h] 
	\centering
	\includegraphics[scale=.7]{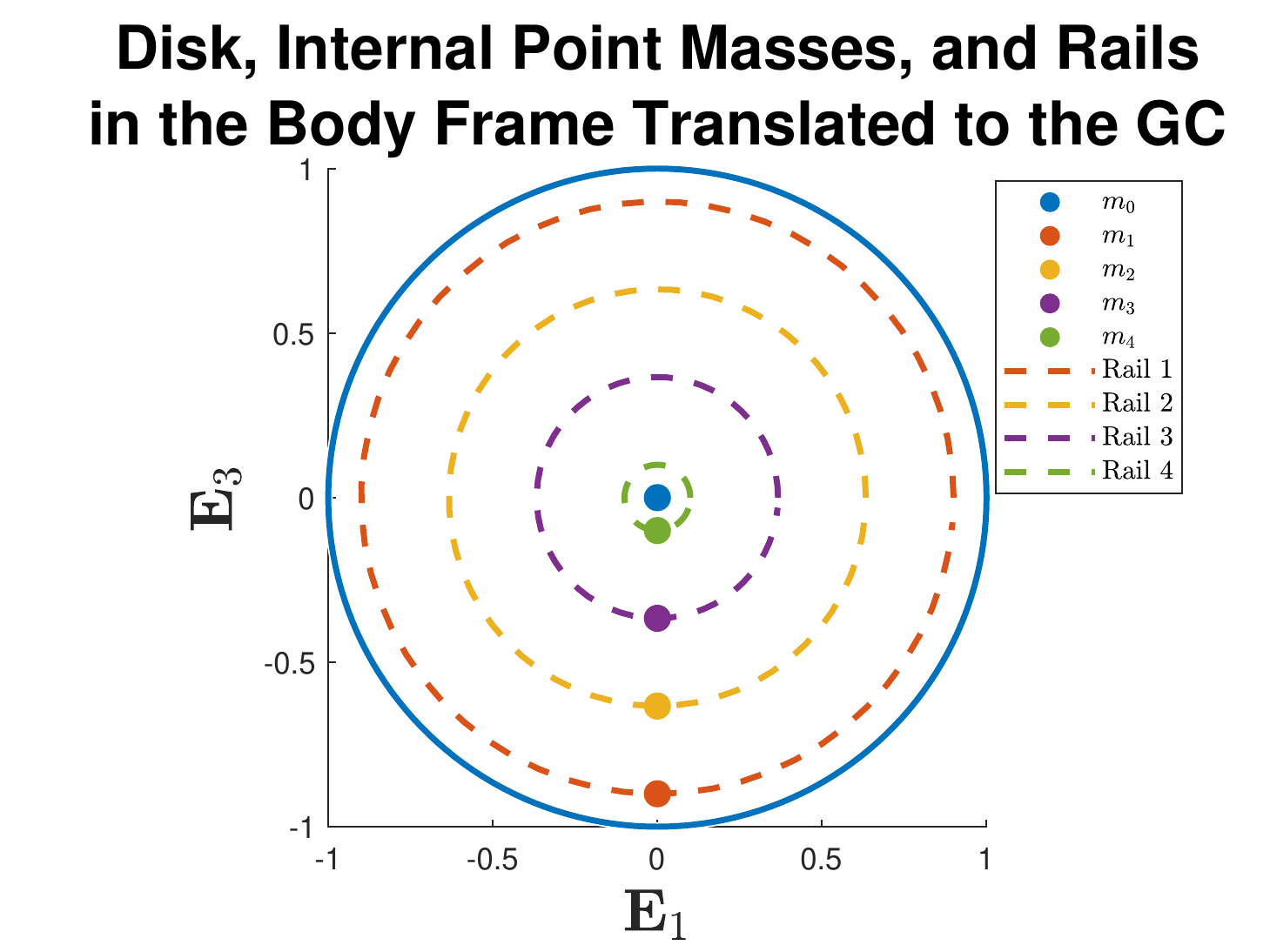}
	\caption{A disk of radius $r=1$ actuated by $4$ internal point masses, $m_1$, $m_2$, $m_3$, and $m_4$, each on its own circular rail of radius $r_1=.9$, $r_2=.6\overline{3}$, $r_3=.3\overline{6}$, and $r_4=.1$, respectively. The location of the disk's CM coincides with the GC and is denoted by $m_0$. $m_0=m_1=m_2=m_3=m_4=1$ and $g=1$. The configuration at the initial time $t=0$ is shown.}
	\label{fig_disk_masses_rails_bf_gc}
\end{figure}

\begin{figure}[!ht] 
	\centering
	\subfloat[Trajectories of the disk's internal point masses and of the total system center of mass in the body frame translated to the GC.]{\includegraphics[scale=.5]{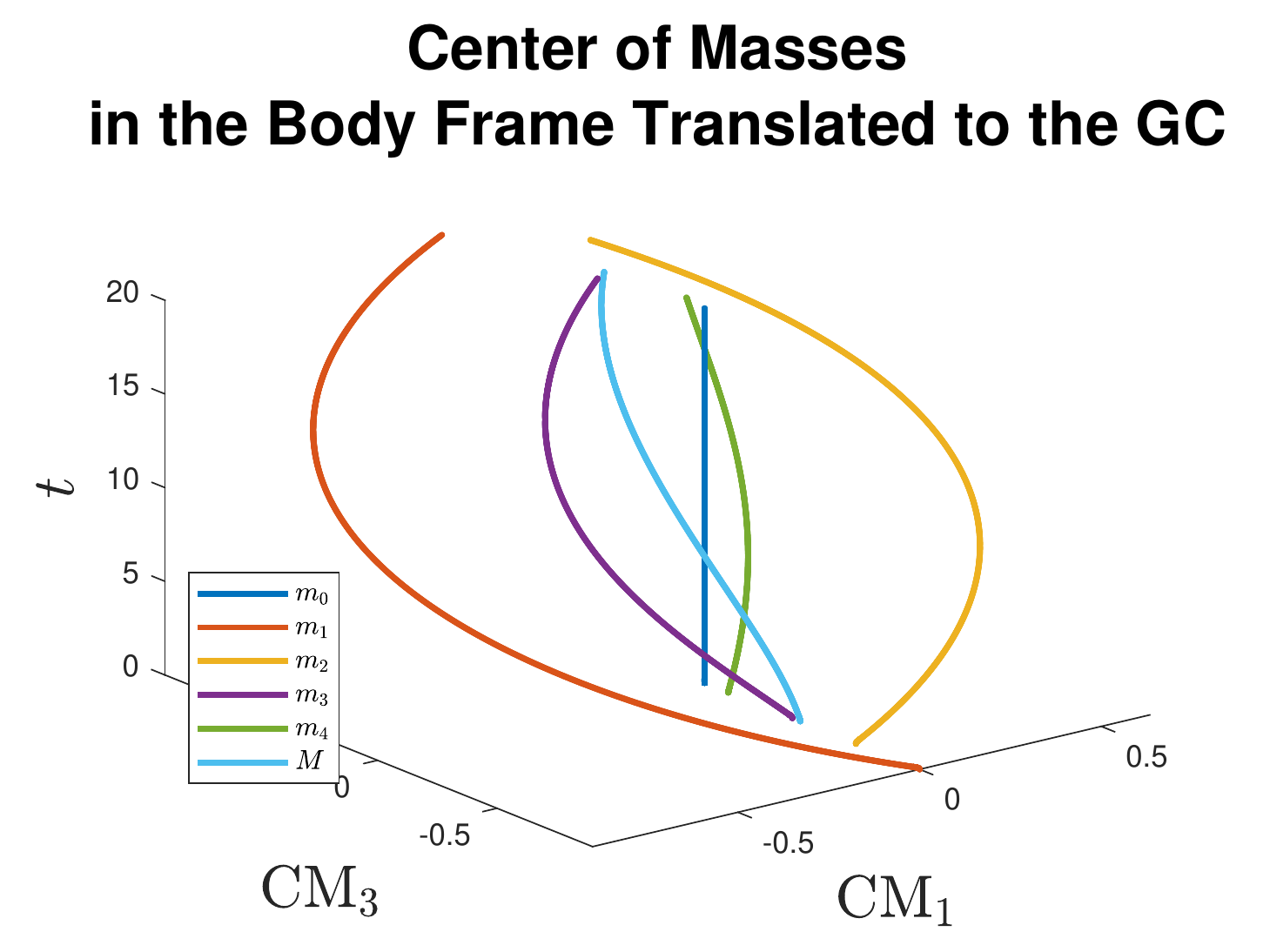}\label{fig_disk_cm_time_bf_gc}}
	\hspace{5mm}
	\subfloat[Trajectories of the disk's internal point masses and of the total system center of mass in the spatial frame translated to the GC.]{\includegraphics[scale=.5]{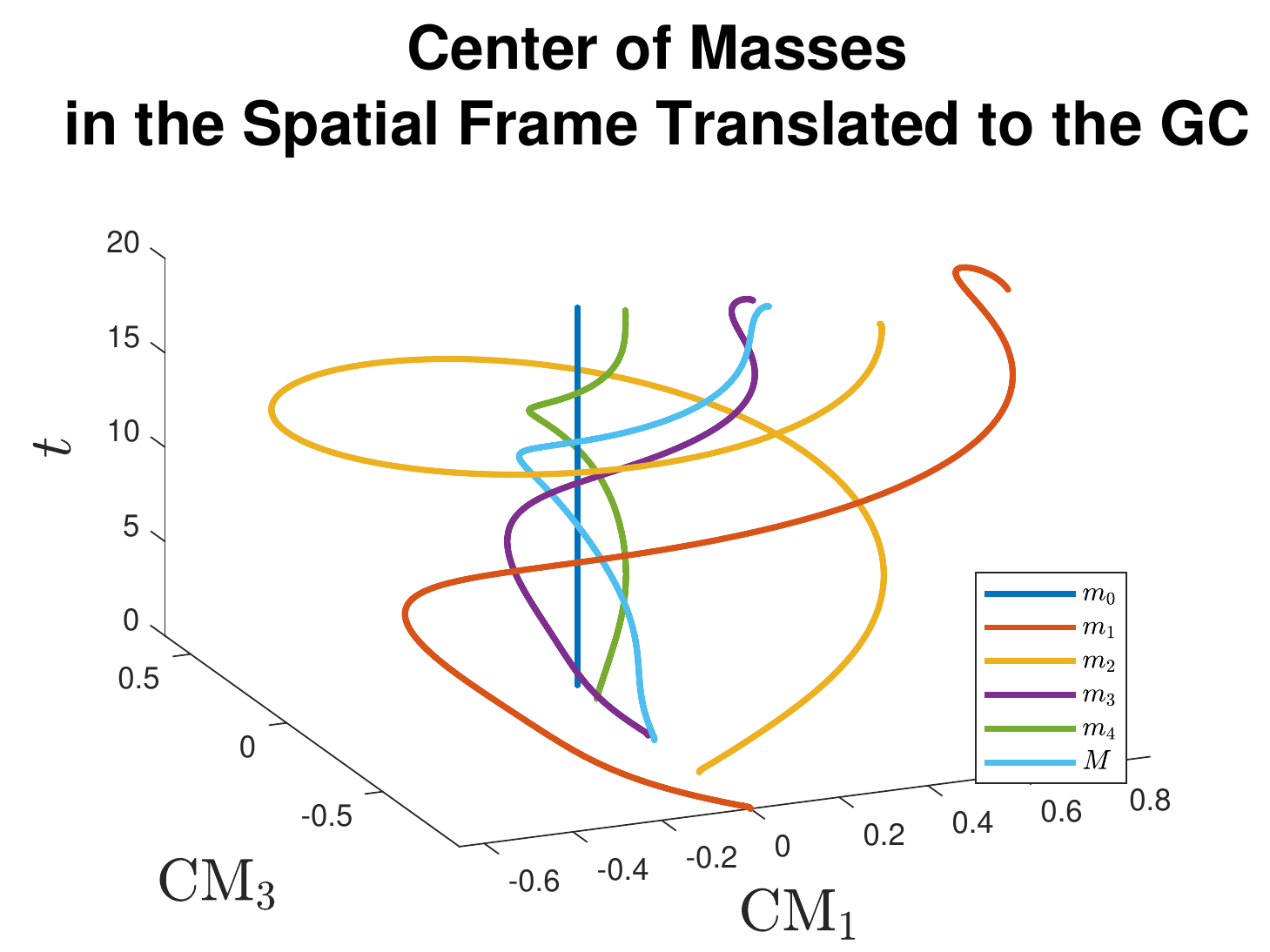}\label{fig_disk_cm_time_spatial_gc}}
	\\
	\subfloat[Evolution of the time derivative of the disk's rotation angle.]{\includegraphics[scale=.5]{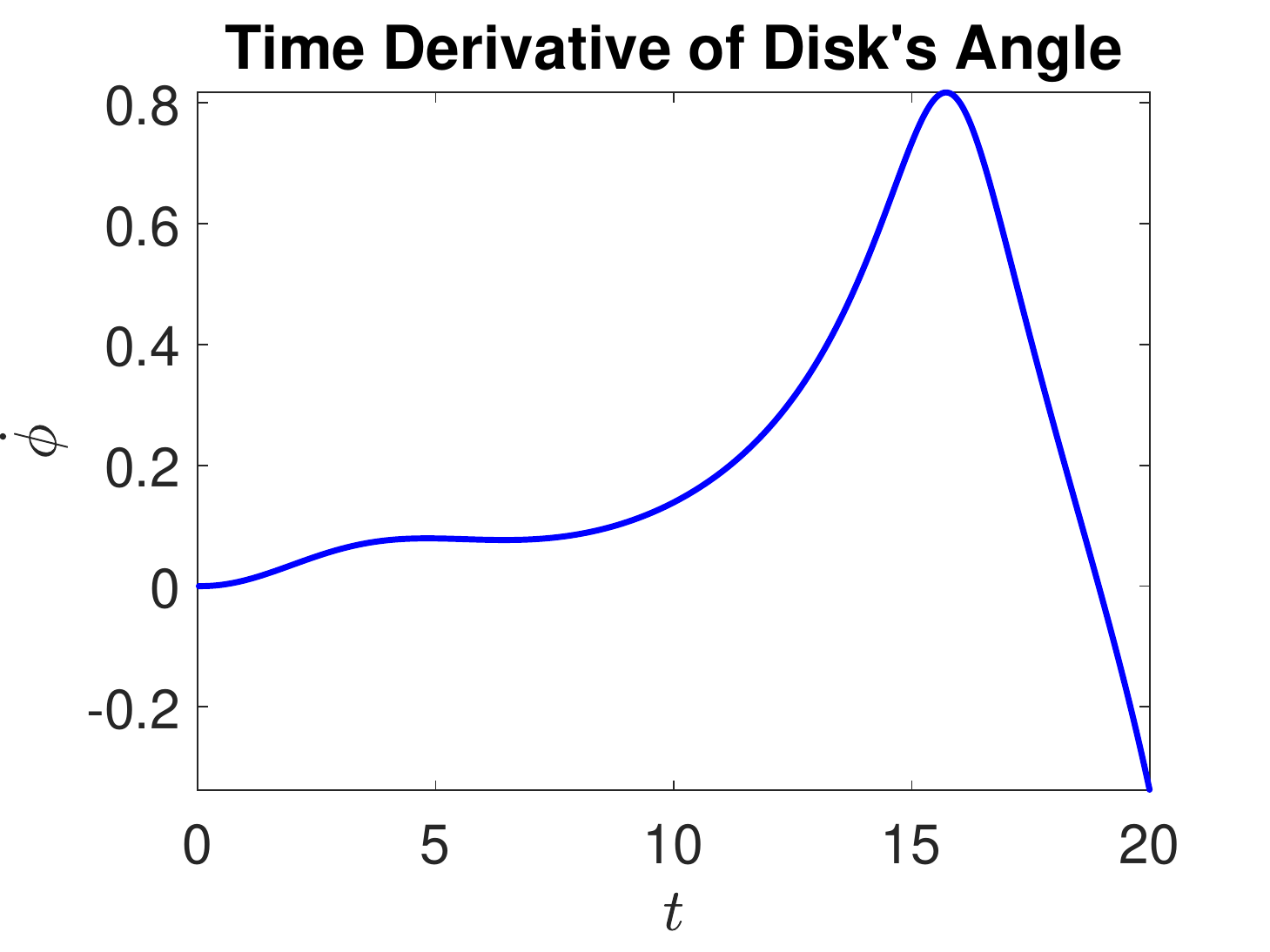}\label{fig_disk_dot_angle}}
	\hspace{5mm}
	\subfloat[Trajectory of the disk's GC and CP.]{\includegraphics[scale=.5]{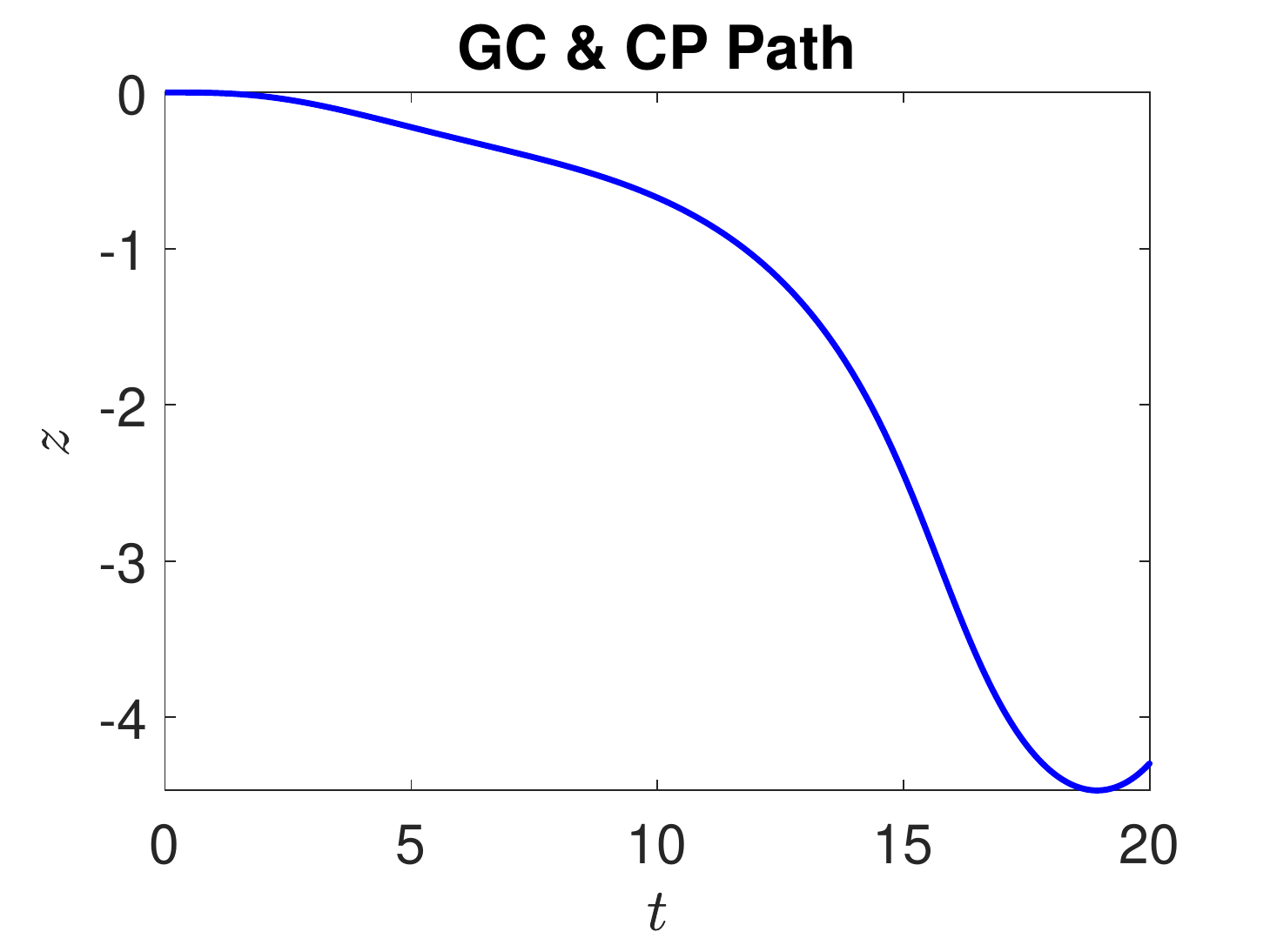}\label{fig_disk_gc}}
	\caption{Dynamics of the rolling disk shown in Figure~\ref{fig_disk_masses_rails_bf_gc} obtained by numerically integrating the ODE IVP \eqref{eq_rdisk_dynamics}, \eqref{eq_disk_initial_conds_spec} with \mcode{ode45} over the time interval $\left[0,20\right]$. The parameterized accelerations of the internal point masses are given in \eqref{eq_prescribed_u_disk}.}
	\label{fig_disk_ode_sims}
\end{figure}

\section{Numerical Simulations of the Dynamics of the Rolling Ball}
\label{sec_ball_numerical}
To write the equations of motion for the rolling ball in the standard ordinary differential/algebraic equation (ODE/DAE) form, the state of the system is defined as
\begin{equation}
\bx \equiv \begin{bmatrix} \btheta \\ \dot \btheta \\ \mathfrak{q} \\ \bOm \\ \bz  \end{bmatrix},
\end{equation}
where $\btheta, \, \dot \btheta \in \mathbb{R}^n$ encode the positions and velocities of the moving masses, the versor $\mathfrak{q} \in \mathscr{S} \cong \mathbb{S}^3 \subset \mathbb{R}^4$ encodes the orientation of the rolling ball, $\bOm \in \mathbb{R}^3$ is the body angular velocity, and $\bz \in \mathbb{R}^2$ \revision{R2Q2}{denotes the spatial $\mathbf{e}_1$- and $\mathbf{e}_2$-components of the GC and CP.} \rem{is the projection of the GC onto the horizontal rolling plane.} \revision{R1Q5}{Appendix~\ref{app_quaternions} provides a brief review of quaternions and versors. }
\color{black} 
ODE and DAE formulations of the rolling ball's system dynamics defined for $a \le t \le b$ are
\begin{equation} \label{rolling_ball_ode_f}
\dot {\bx} = \begin{bmatrix} \dot \btheta \\ \ddot \btheta \\ \dot {\mathfrak{q}} \\ \dot \bOm \\ \dot \bz  \end{bmatrix}  = \mathbf{f}\left(t,\bx,\bu\right) \equiv \begin{bmatrix} \dot \btheta \\ \bu  \\ \frac{1}{2} \mathfrak{q} \bOm^\sharp \\ \bkappa\left(t,\bx,\bu \right) \\ \left( \left[\mathfrak{q} \bOm^\sharp \mathfrak{q}^{-1} \right]^\flat \times r \mathbf{e}_3  \right)_{12}  \end{bmatrix}
\end{equation}
and
\begin{equation} \label{rolling_ball_dae_g}
\mathcal{M} \dot {\bx} = \begin{bmatrix} \dot \btheta \\ \ddot \btheta \\ 0 \\ {\dot {\mathfrak{q}}}^\flat \\ \dot \bOm \\ \dot \bz  \end{bmatrix}  = \mathbf{g}\left(t,\bx,\bu\right) \equiv \begin{bmatrix} \dot \btheta \\ \bu  \\ \left| \mathfrak{q} \right|^2-1 \\ \left[\frac{1}{2} \mathfrak{q} \bOm^\sharp\right]^\flat \\ \bkappa\left(t,\bx,\bu \right) \\ \left( \left[\mathfrak{q} \bOm^\sharp \mathfrak{q}^{-1} \right]^\flat \times r \mathbf{e}_3  \right)_{12}  \end{bmatrix},
\end{equation}
respectively, where $\bu \colon \mathbb{R} \to \mathbb{R}^n$ is a prescribed function of $t$ such that $\bu(t)=\ddot \btheta(t) \in \mathbb{R}^n$, $\bkappa\left(t,\bx,\bu \right)$ is given by the right-hand side of the formula for $\dot \bOm$ in \eqref{uncon_ball_eqns_explicit_1d}:
\begin{equation} \label{eq_ball_kappa}
\begin{split}
\bkappa\left(t,\bx,\bu \right) \equiv
\left[\sum_{i=0}^n m_i \widehat{\mathbf{s}_i}^2  -\inertia \right]^{-1}  \Bigg[&\bOm \times \inertia \bOm+r \tilde \bGamma \times \bGamma\\
&+ \sum_{i=0}^n m_i \mathbf{s}_i \times  \left\{ g \bGamma+ \bOm \times \left(\bOm \times \bzeta_i +2 \dot \theta_i \bzeta_i^{\prime} \right) + \dot \theta_i^2 \bzeta_i^{\dprime} + \ddot \theta_i  \bzeta_i^{\prime} \right\}  \Bigg],
\end{split}
\end{equation}
and
\begin{equation}
\mathcal{M} \equiv \diag \left( \begin{bmatrix} \mathbf{1}_{1 \times 2n} & 0 & \mathbf{1}_{1 \times 8} \end{bmatrix} \right)
\end{equation}
is a diagonal DAE mass matrix. Observe that \eqref{rolling_ball_dae_g} is a semi-explicit DAE of index 1.

In order to construct $\bkappa\left(t,\bx,\bu \right)$ as defined above, the variables $\bGamma \equiv \Lambda^{-1} \mathbf{e}_3$ and $\tilde \bGamma \equiv \Lambda^{-1} \mathbf{F}_\mathrm{e}$ must be  computed first. Given a versor $\mathfrak{q}$, $\bGamma$ and $\tilde \bGamma$ can be computed by first constructing $\Lambda$ from $\mathfrak{q}$ or directly from $\mathfrak{q}$ by using the Euler-Rodrigues  formulas $\bGamma \equiv \Lambda^{-1} \mathbf{e}_3=\left[\mathfrak{q}^{-1} \mathbf{e}_3^\sharp \mathfrak{q} \right]^\flat$ and $\tilde \bGamma \equiv \Lambda^{-1} \mathbf{F}_\mathrm{e} = \left[\mathfrak{q}^{-1} \mathbf{F}_\mathrm{e}^\sharp \mathfrak{q} \right]^\flat$. 

Likewise, the final formula for  computing the velocity of the GC in \eqref{rolling_ball_ode_f} and \eqref{rolling_ball_dae_g} is $\dot \bz=\left( \bom \times r \mathbf{e}_3  \right)_{12}$, where $\bom =  \Lambda \bOm=\left[\mathfrak{q} \bOm^\sharp \mathfrak{q}^{-1} \right]^\flat$  by the Euler-Rodrigues formula. Thus, given $\mathfrak{q}$, the spatial angular velocity $\bom$ can be obtained by first computing $\Lambda$ from $\mathfrak{q}$ or directly from $\mathfrak{q}$ via $\bom=\left[\mathfrak{q} \bOm^\sharp \mathfrak{q}^{-1} \right]^\flat$. The most computationally efficient method to determine the variables $\bGamma$, $\tilde \bGamma$, and $\bom$ is to use the formulas
 \begin{equation}
 \label{Variables_calc}
  \bGamma\equiv \Lambda^{-1} \mathbf{e}_3, \quad \tilde \bGamma\equiv \Lambda^{-1} \mathbf{F}_\mathrm{e}, \quad \mathrm{and} \quad \bom = \Lambda \bOm, 
  \end{equation} 
where one would  first construct $\Lambda$ from $\mathfrak{q}$, and then use this matrix and its inverse $\Lambda^{-1}=\Lambda^T$ to compute $\bGamma$, $ \tilde \bGamma$, and $\bom$ according to \eqref{Variables_calc} above. 

In order to simulate the rolling ball's dynamics, \eqref{rolling_ball_ode_f} or \eqref{rolling_ball_dae_g} must be integrated with prescribed initial conditions at time $t=a$:
\begin{equation} \label{eq_ball_initial_conds}
\bx\left(a\right) = \begin{bmatrix} \btheta(a) \\ \dot \btheta(a) \\ \mathfrak{q}(a) \\ \bOm(a) \\ \bz(a) \end{bmatrix} = \begin{bmatrix}  \btheta_a \\  {\dot \btheta}_a \\ \mathfrak{q}_a \\ \bOm_a \\ \bz_a \end{bmatrix} \equiv \bx_a.
\end{equation}
\eqref{rolling_ball_ode_f} and \eqref{eq_ball_initial_conds} constitute an ODE IVP, while \eqref{rolling_ball_dae_g} and \eqref{eq_ball_initial_conds} constitute a DAE IVP.

 In the simulations, we consider a rolling ball of mass $m_0=1$, radius $r=1$, principal moments of inertia $d_1=.9$, $d_2=1$, and $d_3=1.1$, and with the CM shifted slightly away from the GC at $\bzeta_0=\begin{bmatrix} 0 & 0 & -.05 \end{bmatrix}^\mathsf{T}$. The ball contains $n=3$ internal point masses, each of mass $1$ so that $m_1=m_2=m_3=1$ and each located on its own circular rail centered on the GC of radius $r_1=.95$, $r_2=.9$, and $r_3=.85$, respectively, oriented as shown in Figure~\ref{fig_ball_masses_rail_bf_gc}. 
The total mass of the ball's  system  is $M=4$, and  gravity is rescaled to be $g=1$.
\rem{\todo{VP: The example is quite simple, CM is in GC and inertia matrix is proportional to unity, which makes a lot of terms cancel. Of course, the total moment of inertia and CM is no longer trivial. Was there a particular reason you chose that simple case for the ball without masses? Some terms in equations will cancel, although it is still non-trivial. \\ SMR: I re-simulated the ball using a CM slightly shifted from the GC and using unequal principal moments of inertia. Are these new CM and principal moment values and the associated simulation plots ok? Do you also want me to re-simulate the disk using an off-center CM? I am a bit hesitant to re-simulate the controlled ball using these new values for the next paper on control, since the simulations will be tedious to regenerate for new ball parameters due to having to find new predictor-corrector tangent step lengths. \\ 
VP: OK, fair enough. Hope the Referees won't bug us about it.  } }
For $1 \le i \le n$, the position of $m_i$ in the body frame centered on the GC is:
\begin{equation}
\bzeta_i\left(\theta_i\right) = r_i \mathcal{B}_i \left(\boldsymbol{\varsigma} \left(\mathbf{v}_i\right) \right) \begin{bmatrix} \cos \theta_i \\ 0 \\ \sin \theta_i  \end{bmatrix},
\end{equation}
where $\mathcal{B}_i \left(\mathbf{n} \right) \in SO(3)$ is a rotation matrix whose columns are the right-handed orthonormal basis constructed from the unit vector $\mathbf{n} \in \mathbb{R}^3$ based on the algorithm given in Section 4 and Listing 2 of \cite{frisvad2012building}, $\boldsymbol{\varsigma} \colon \mathbb{R}^3 \to \mathbb{R}^3$ maps spherical coordinates to Cartesian coordinates:
\begin{equation}
\boldsymbol{\varsigma}\left(\begin{bmatrix} \phi \\ \theta \\ \rho \end{bmatrix} \right) = \begin{bmatrix} \rho \cos \theta \cos \phi  \\ \rho \cos \theta \sin \phi \\ \rho \sin \theta \end{bmatrix},
\end{equation}
and
\begin{equation}
\mathbf{v}_1 = \begin{bmatrix} 0 & 0 & 1 \end{bmatrix}^\mathsf{T}, \quad \mathbf{v}_2 = \begin{bmatrix} \frac{\pi}{2} & 0 & 1 \end{bmatrix}^\mathsf{T}, \quad \mathrm{and} \quad \mathbf{v}_3 = \begin{bmatrix} \frac{\pi}{4} & \frac{\pi}{4} & 1 \end{bmatrix}^\mathsf{T}
\end{equation}
are spherical coordinates of unit vectors in $\mathbb{R}^3$. There is no external force acting on the ball's GC so that $\mathbf{F}_\mathrm{e} = \tilde \bGamma  = \mathbf{0}$ in \eqref{eq_ball_kappa}. This ball's dynamics are simulated with initial time $a=0$ and final time $b=20$, so that the simulation time interval is $\left[0,20\right]$. The parameterized acceleration of each internal point mass is a continuous approximation of a short duration unit amplitude step function:
\begin{equation} \label{eq_prescribed_u_ball}
u_i(t) = \ddot \theta_i(t) = \left\{ \begin{array}{ll} 1, & 0 \le t \le .1, \\ \unaryminus 10t+2, & .1 \le t \le .2, \\ 0, & .2 \le t \le 20, \end{array} \right. \quad \mbox{for} \quad 1 \le i \le n.
\end{equation}
A plot of the magnitude of \eqref{eq_prescribed_u_ball} is depicted in Figure~\ref{fig_abs_ui}. The rolling ball's initial conditions are selected so that the ball starts at rest at the origin. Table~\ref{table_ball_ICs} shows parameter values used in the rolling ball's initial conditions \eqref{eq_ball_initial_conds}.   The initial orientation matrix is selected to be the identity matrix so that $\mathfrak{q}_a = \begin{bmatrix} 1 & 0 & 0 & 0 \end{bmatrix}^\mathsf{T}$ and the initial configurations of the internal point masses are given by $\btheta_a = \begin{bmatrix} 0 & 2.0369 & 0.7044 \end{bmatrix}^\mathsf{T}$, so that the ball's total system center of mass is initially located above the GC. 
 These particular initial configurations of the point masses were obtained by solving a system of algebraic equations for mass positions  based on the requirement that the ball's  total system center of mass be  directly above  or below the GC.  
To ensure that the ball is initially at rest, ${\dot \btheta}_a = \begin{bmatrix} 0 & 0 & 0 \end{bmatrix}^\mathsf{T}$ and $\bOm_a = \begin{bmatrix} 0 & 0 & 0 \end{bmatrix}^\mathsf{T}$. To ensure that the ball's GC is initially located at the origin, $\bz_a = \begin{bmatrix} 0 & 0 \end{bmatrix}^\mathsf{T}$. In summary, the rolling ball's initial conditions are
\begin{equation} \label{eq_ball_initial_conds_spec}
\bx_a = \begin{bmatrix} 0 & 2.0369 & 0.7044 & 0 & 0 & 0 & 1 & 0 & 0 & 0 & 0 & 0 & 0 & 0 & 0 \end{bmatrix}^\mathsf{T}.
\end{equation}

\begin{table}[h!]
	\centering 
	{ 
		\setlength{\extrarowheight}{1.5pt}
		\begin{tabular}{| c | c |} 
			\hline
			\textbf{Parameter} & \textbf{Value} \\ 
			\hline\hline 
			$\btheta_a$ & $\begin{bmatrix}  0 & 2.0369 & .7044 \end{bmatrix}^\mathsf{T}$  \\  
			\hline
			$\dot \btheta_a$ & $\begin{bmatrix} 0 & 0 & 0 \end{bmatrix}^\mathsf{T}$ \\ 
			\hline
			$\mathfrak{q}_a$ & $\begin{bmatrix} 1 &  0 & 0 & 0 \end{bmatrix}^\mathsf{T}$ \\
			\hline
			$\bOm_a$ & $\begin{bmatrix}  0 & 0 & 0 \end{bmatrix}^\mathsf{T}$ \\ 
			\hline
			$\bz_a$ & $\begin{bmatrix}  0 & 0 \end{bmatrix}^\mathsf{T}$ \\ 
			\hline
		\end{tabular} 
	}
	\caption{Initial condition parameter values for the rolling ball.}
	\label{table_ball_ICs}
\end{table}

The dynamics of this rolling ball are simulated by numerically integrating the ODE IVP \eqref{rolling_ball_ode_f}, \eqref{eq_ball_initial_conds_spec} or the DAE IVP \eqref{rolling_ball_dae_g}, \eqref{eq_ball_initial_conds_spec}. The ODE IVP \eqref{rolling_ball_ode_f}, \eqref{eq_ball_initial_conds_spec} is numerically integrated via the \mcode{MATLAB} R2017b routines \mcode{ode45}, \mcode{ode113}, \mcode{ode15s}, \mcode{ode23t}, and \mcode{ode23tb} and a \mcode{MATLAB} \mcode{MEX} wrapper of the Fortran routine radau5 \cite{hairer1996solving}, while the DAE IVP \eqref{rolling_ball_dae_g}, \eqref{eq_ball_initial_conds_spec} is numerically integrated via the \mcode{MATLAB} R2017b routines \mcode{ode15s} and \mcode{ode23t} and a \mcode{MATLAB} \mcode{MEX} wrapper of the Fortran routine radau5. Except for the absolute and relative error tolerances and the Jacobian, all the numerical integrators are used with the default input options. The absolute and relative error tolerances supplied to the numerical integrators are both set to $1\mathrm{e}{-10}$. Jacobions of $\mathbf{f}$ and $\mathbf{g}$ with respect to the state $\bx$, obtained via complex-step differentiation \cite{squire1998using,martins2001connection,martins2003complex}, are supplied to \mcode{ode15s}, \mcode{ode23t}, \mcode{ode23tb}, and radau5, depending on whether the ODE or DAE IVP is numerically integrated. Since excellent agreement was observed between all the numerical integrators, only the results obtained by numerically integrating the DAE IVP \eqref{rolling_ball_dae_g}, \eqref{eq_ball_initial_conds_spec} with radau5 are shown in Figure~\ref{fig_ball_dae_sims}. As was the case for the rolling disk, \mcode{ode113} completed the numerical integration of the rolling ball's equations of motion in the shortest time.

\begin{figure}[h] 
	\centering
	\includegraphics[scale=.7]{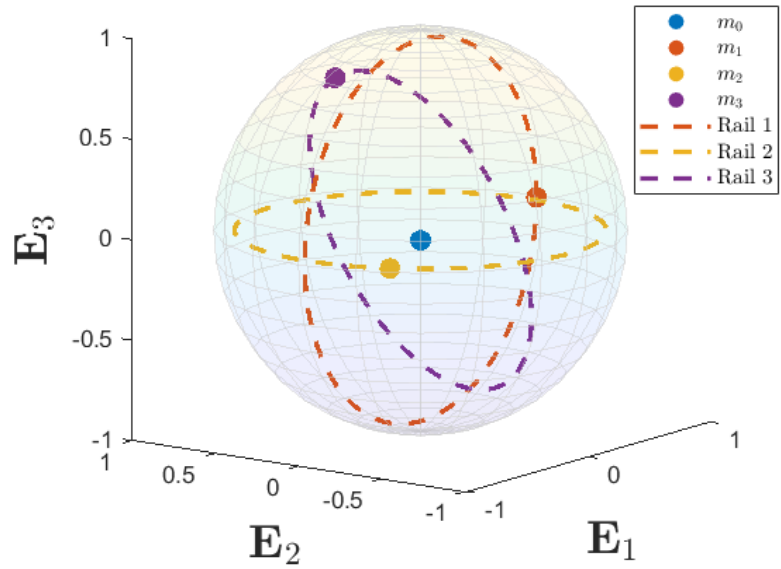}
	\caption{A ball of radius $r=1$ actuated by $3$ internal point masses, $m_1$, $m_2$, and $m_3$, each on its own circular rail of radius $r_1=.95$, $r_2=.9$, and $r_3=.85$, respectively. The location of the ball's CM is shifted slightly away from the GC and is denoted by $m_0$. $m_0=m_1=m_2=m_3=1$ and $g=1$. The configuration at the initial time $t=0$ is shown.}
	\label{fig_ball_masses_rail_bf_gc}
\end{figure}

\begin{figure}[!ht] 
	\centering
	\subfloat[Trajectories of the ball's internal point masses and of the total system center of mass in the body frame translated to the GC.]{\includegraphics[scale=.5]{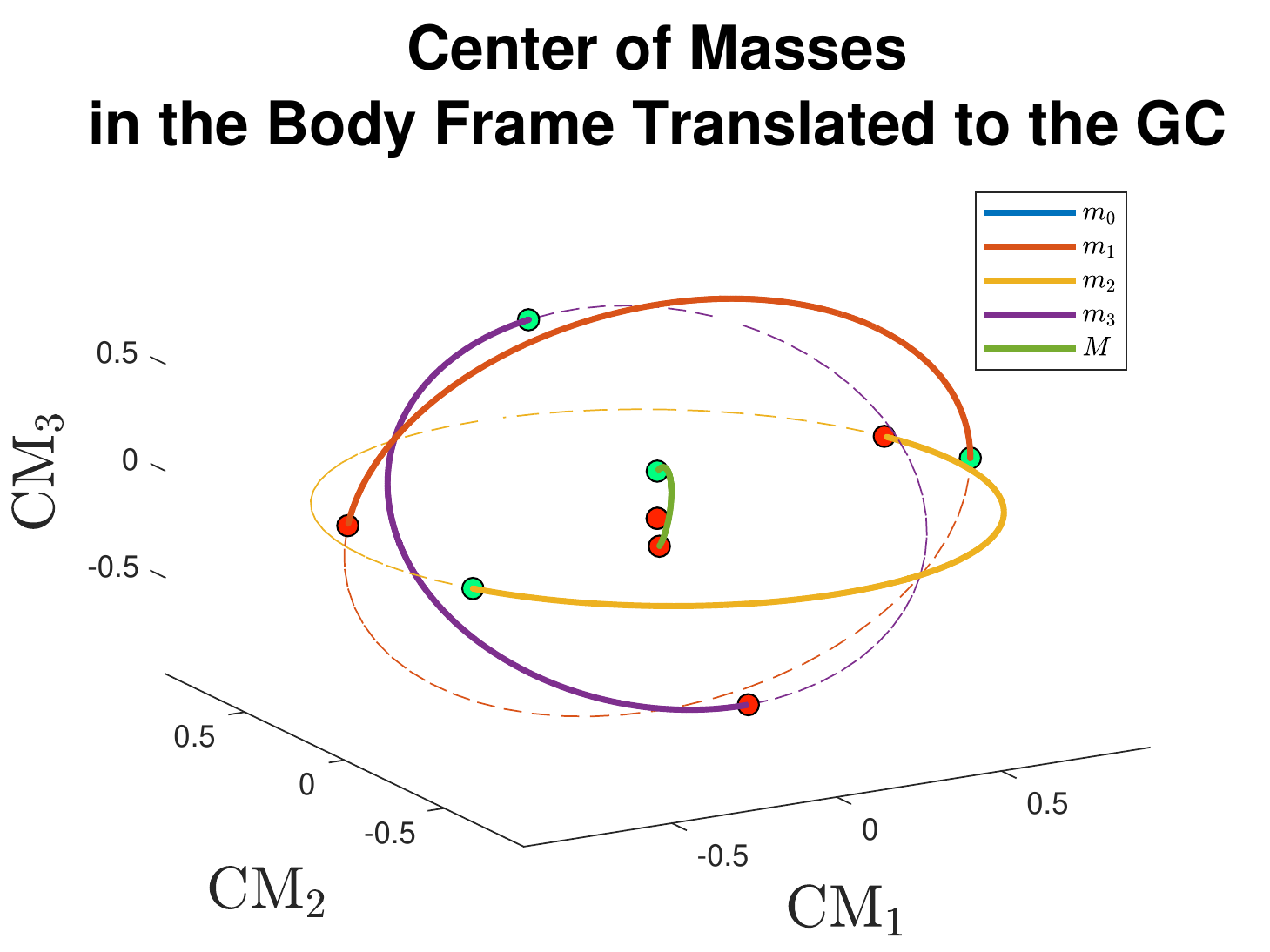}\label{fig_ball_cm_bf_gc}}
	\hspace{5mm}
	\subfloat[Trajectories of the ball's internal point masses and of the total system center of mass in the spatial frame translated to the GC.]{\includegraphics[scale=.5]{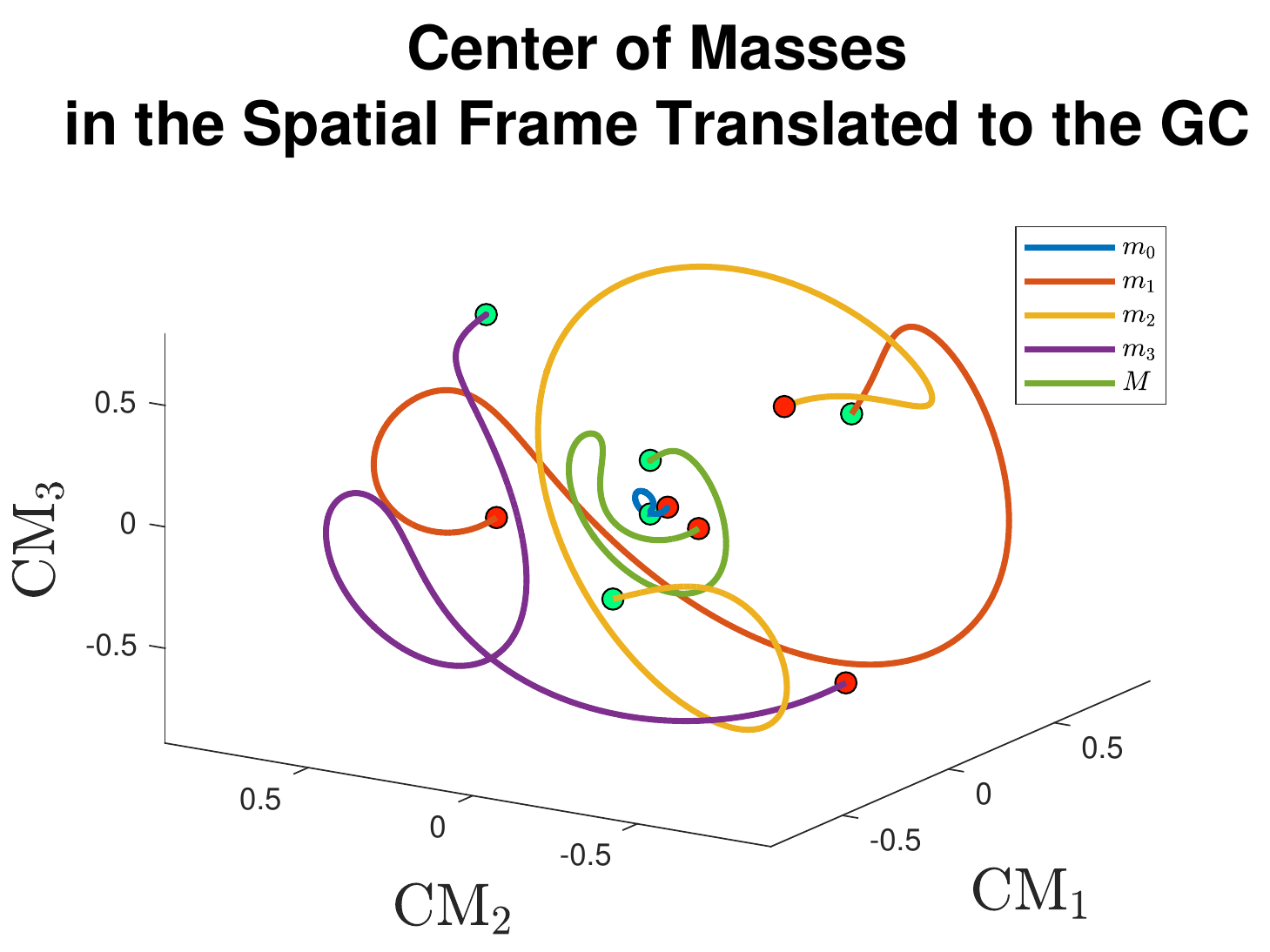}\label{fig_ball_cm_spatial_gc}}
	\\
	\subfloat[Evolution of the ball's body angular velocity.]{\includegraphics[scale=.5]{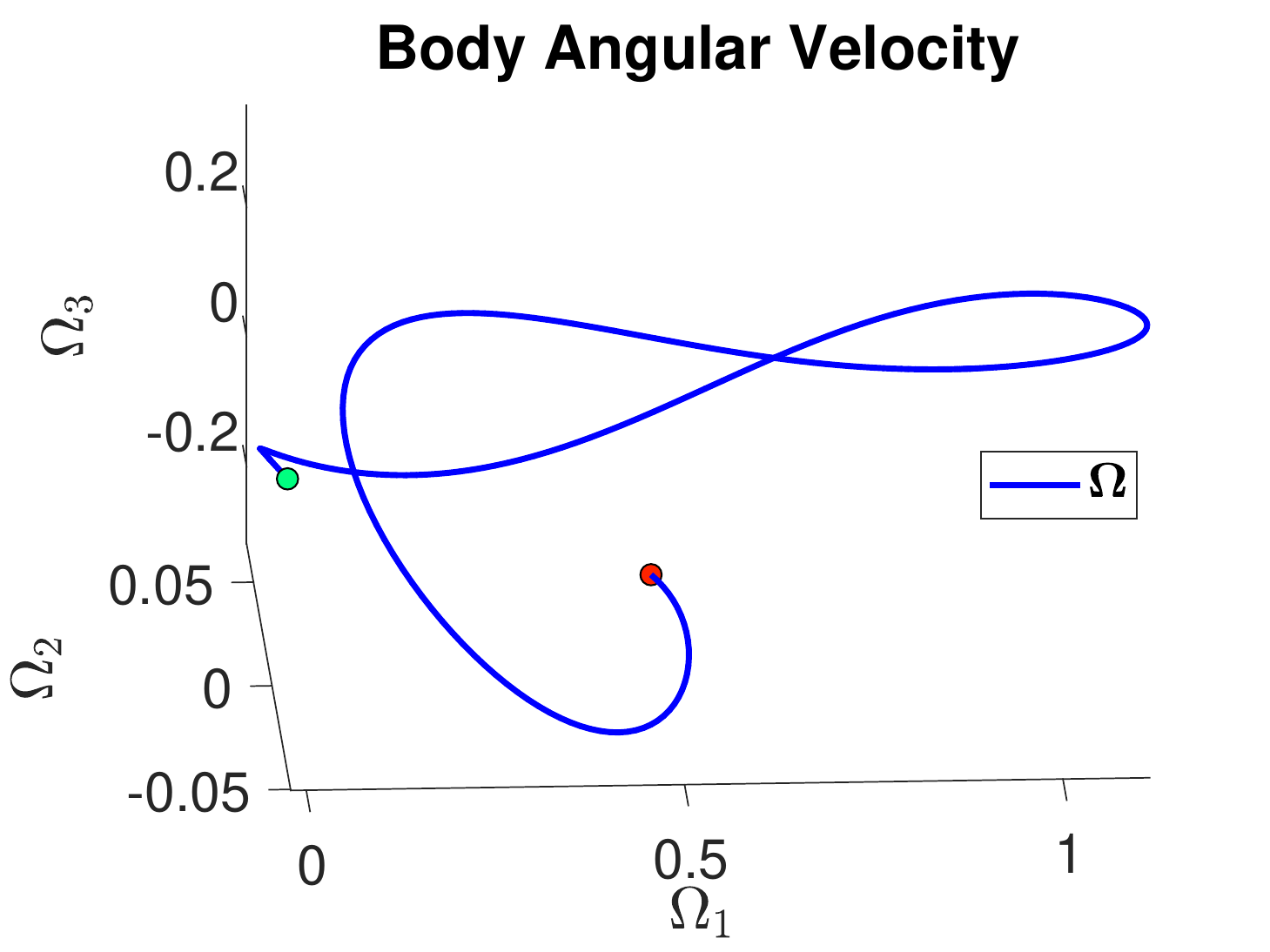}\label{fig_ball_Omega}}
	\hspace{5mm}
	\subfloat[Trajectory of the ball's GC and CP.]{\includegraphics[scale=.5]{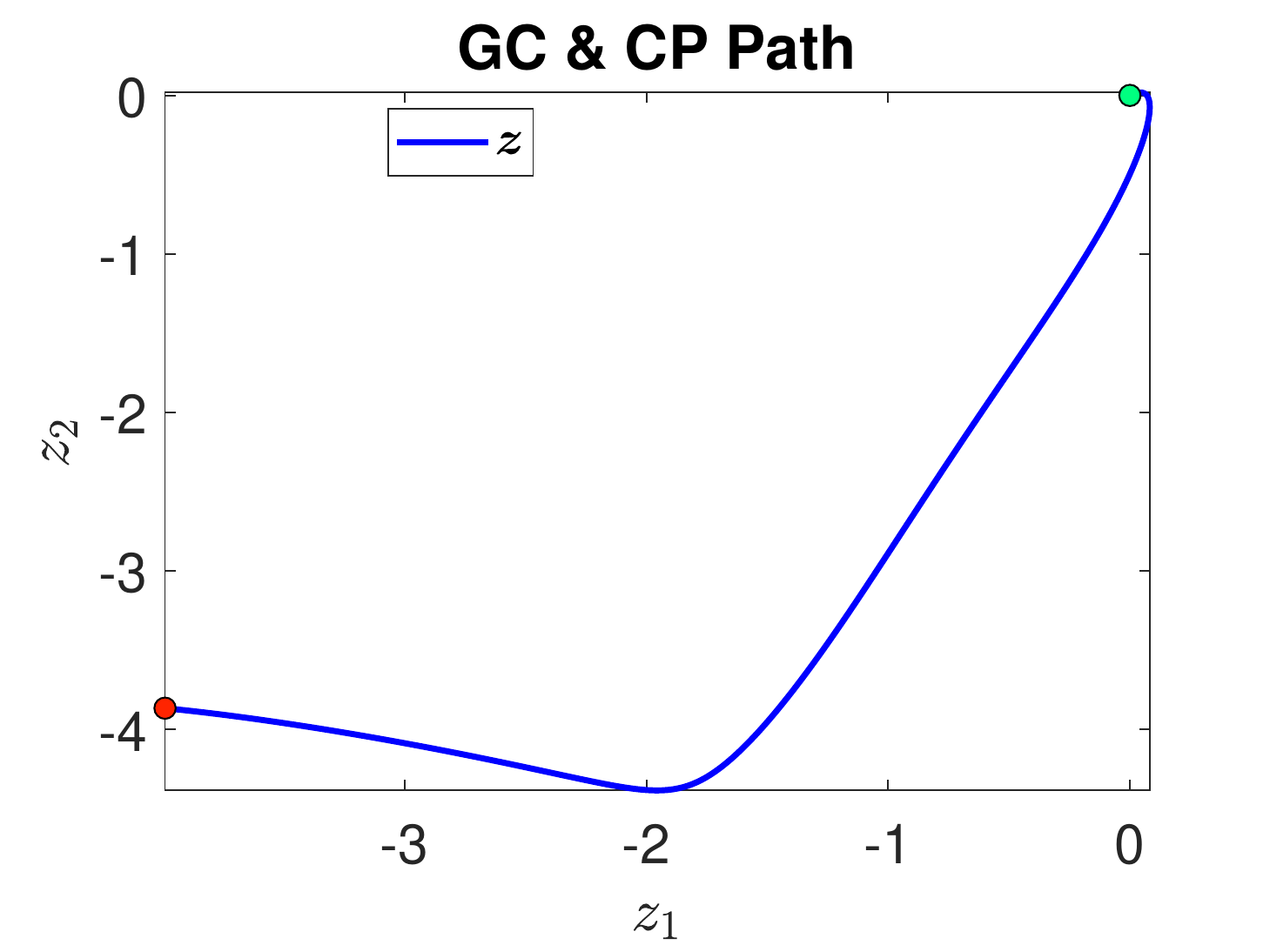}\label{fig_ball_gc_path}}
	\caption{Dynamics of the rolling ball shown in Figure~\ref{fig_ball_masses_rail_bf_gc} obtained by numerically integrating the DAE IVP \eqref{rolling_ball_dae_g}, \eqref{eq_ball_initial_conds_spec} with radau5 over the time interval $\left[0,20\right]$. The parameterized accelerations of the internal point masses are given in \eqref{eq_prescribed_u_ball}.}
	\label{fig_ball_dae_sims}
\end{figure}

\section{Conclusions and Future Work} \label{sec_conclusions}

In this paper, we have developed a consistent theory of motion for a rolling ball actuated by internal point masses moving along trajectories fixed in the ball's reference frame. 
We have described the motion of a general three-dimensional ball and the special case of a rolling disk actuated by internal point masses. For the latter case, we have shown that the  equations obtained  from the variational derivation coincide exactly with the force balance equations obtained by balancing forces in Newton's second law. For general three-dimensional motion of the ball and  its masses, the Newtonian derivation is, in our opinion, too cumbersome, and our derivation is advantageous as it leads to the derivation of the equations of motion using a direct, algorithmic approach that is difficult to reproduce by balancing forces in the non-inertial frame of the moving ball.  This is especially true for highly complex motions of  the internal  masses,  such as is the case  for a  ball  actuated by  several double  pendula. Because of the  increased complexity of   practical actuation mechanisms, we believe that  the  variational (Lagrange-d'Alembert's)  principle is  superior to the direct balance of forces in Newton's  laws, since the use of Lagrange-d'Alembert's  principle is methodical and algorithmic in nature, even for the most complex configurations. However, one should  realize  that every problem solved by  Lagrange-d'Alembert's  principle  can also be solved by  Newton's laws, provided that all the forces are  accounted for, which may be quite difficult for a complex system consisting of many interacting parts and constraints.   

A second paper \cite{putkaradze2017optimal}  on the subject  focuses on the control of the rolling ball actuated by internal point masses. Besides that, an interesting avenue would be to introduce friction  acting on the ball  caused by  friction with the substrate or dissipation  induced by the external media. The exact form  of the friction for the three-dimensional motion of the rolling ball is still rather controversial and subject to considerable discussion. We refer the reader to recent general theory derived  in \cite{kozlov2010lagrangian,kozlov2011friction,karavaev2017dynamical} which, with the right fitting of yet unknown forms for experimental  friction functions, may be used to consistently incorporate friction into our model. However, the derivation of the form of even such simple laws for nonholonomic systems from  first principles is nontrivial  and is definitely beyond the scope of this article. We shall postpone the discussion of this interesting problem for future work. 

\rem{
\todo{VP: This reference is to the Russian article, I think there should be a translation? If you find it, let us refer to the English translation. \\ SMR 1: I searched google and scholar.google.com, but I was unable to find an English translation of this Russian article. I used the bib reference for this article from scholar.google.com instead of your bib reference. Kozlov published another article on friction in English, "Friction by Painlev\'e and Lagrangian Mechanics" in 2011: \url{https://link.springer.com/article/10.1134\%2FS1028335811060115?LI=true}. Is this English article relevant? \\\ SMR 2: Above you wrote "... to introduce friction on the ball which may be caused by the friction with the substrate\sout{,} or dissipation caused by the external media." What is the difference between friction with the substrate and dissipation caused by the external media? \\ 
VP: Friction with the substrate, or so-called rolling friction, is usually referred to as interaction of the ball with the substrate. More precisely, the ball deforms the substrate which causes propagation of waves and deformations in the substrate. The friction with the media occurs because the molecules of air, or water (whatever media the ball is surrounded by), are hitting the ball's surface. So if you put the ball in vacuum, but it  rolls on the table, you will still have friction due to rolling friction, even though there is no media to interact with. 
\\
VP: I believe the paper you found  is the  translation of Kozlov's article, let us refer to that one instead. I strengthened the statement a bit. Or, we can include both the Russian and English articles. 
\\ SMR: The English article I mention above is a translation of: "Original Russian Text  V.V. Kozlov, 2011, published in Doklady Akademii Nauk, 2011, Vol. 438, No. 6, pp. 758 to 761," as stated on that webpage given by the link above. So there doesn't seem to be an English translation of the Russian article you cite in the Conclusions.
\\ 
VP: I guess what happened was he has published two very similar articles. I am OK with citing them both. As far as I could see, the material contained in the english version you found overlapped quite a bit with the Russian article. I guess Doklady is a short paper, and the other one was the longer paper, or something like that. 
  } }
\section*{Acknowledgements}

We are indebted to our colleagues  A.M. Bloch, D.M. de Diego, F. Gay-Balmaz, D.D. Holm, M. Leok, A. Lewis, T. Ohsawa, V.A. Vladimirov, and D.V. Zenkov for useful and fruitful discussions. This research was partially  supported by  the NSERC Discovery Grant, the University of Alberta Centennial Fund, and the Alberta Innovates Technology Funding (AITF) which came through the Alberta Centre for Earth Observation Sciences (CEOS). S.M. Rogers also received support from the University of Alberta Doctoral Recruitment Scholarship, the FGSR Graduate Travel Award, the IGR Travel Award, the GSA Academic Travel Award, and the AMS Fall Sectional Graduate Student Travel Grant. The authors declare that they have no conflict of interest.

\phantomsection
\addcontentsline{toc}{section}{References}
\printbibliography
\hypertarget{References}{}

\appendix
\section{Background Material: Variational Mechanics, Euler-Poincar\'e's Method, and Nonholonomic Constraints} \label{app_background} 
This appendix reviews several principles from mechanics that are useful for developing the equations of motion for the rolling ball. Hamilton's principle and Euler-Poincar\'e's method are reviewed in Subappendix~\ref{ssec_EP}, while Lagrange-d'Alembert's principle is reviewed in Subappendix~\ref{ssec_LdA}. Euler-Poincar\'e's method and Lagrange-d'Alembert's principle are later utilized to derive the equations of motion for the rolling ball in Section~\ref{sec_ball_uncontrolled}.

\subsection{Hamilton's Principle, Symmetry Reduction, and Euler-Poincar\'e's Method} \label{ssec_EP}
\paragraph{Hamilton's Principle}
A mechanical system  consists of a configuration space, which is a manifold  $M$ with tangent bundle $TM = \bigcup_{q \in M} T_q M$, and a Lagrangian $L(q,\dot q): TM \rightarrow \mathbb{R}$,  $(q,\dot q) \in TM$.  
\rem{ %%%BEGIN REM 
Figure~\ref{fig_tangent_bundles} illustrates the tangent bundles of several manifolds.
\begin{figure}[h] 
	\centering
	\subfloat[The tangent bundle of the circle \cite{Circle_TB}.]{\includegraphics[width=0.3\textwidth, height=.2\textwidth]{tangent_bundle_circle}\label{fig_tangent_bundle_circle}}
	\hspace{5mm}
	\subfloat[The tangent bundle of the sphere \cite{Sphere_TB}.]{
		\includegraphics[width=0.3\textwidth]{tangent_bundle_sphere}\label{fig_tangent_bundle_sphere}}
	\hspace{5mm}
	\subfloat[The tangent bundle of the M\"obius strip \cite{Mobius_TB}.]{
		\includegraphics[width=0.3\textwidth]{tangent_bundle_mobius}\label{fig_tangent_bundle_mobius}}	
	\caption{Illustrations of the tangent bundles of several manifolds.} \label{fig_tangent_bundles}
\end{figure}
} %%%%END REM 
The equations of motion are given by Hamilton's principle (also called the variational principle of stationary action) which states that 
\begin{equation} 
\de \int_a^b L \left(q, \dot q \right) \mathrm{d} t=0, \quad \de q(a)=\de q(b)=0,
\label{EL_var}
\end{equation} 
for all smooth functions $\de q(t)$ defined for $a\le t \le b$ and that vanish at the endpoints (i.e. $\de q(a)=\de q(b)=0$). Pushing the variational derivative inside the integral, integrating by parts, and enforcing the vanishing endpoint conditions $\de q(a)=\de q(b)=0$ yields
\begin{equation}
\begin{split}
\de \int_a^b L\left(q, \dot q\right) \mathrm{d}t=\int_a^b \de L\left(q, \dot q\right) \mathrm{d}t
&=\int_a^b \left[ \pp{L\left(q, \dot q\right)}{q}\de q + \pp{L\left(q, \dot q\right)}{\dot q}\de \dot q \right]  \mathrm{d} t \\
&=\int_a^b \left[ \pp{L\left(q, \dot q \right)}{q}\de q - \dd{}{t} \pp{L\left(q, \dot q\right)}{\dot q}\de q \right]  \mathrm{d} t
+ \left. \pp{L\left(q, \dot q \right)}{\dot q} \de q \right|_a^b \\
&=\int_a^b \left[ \pp{L\left(q, \dot q\right)}{q} - \dd{}{t} \pp{L\left(q, \dot q\right)}{\dot q} \right] \de q  \mathrm{d} t.
\end{split} 
\end{equation}
Insisting that $\de \int_a^b L\left(q, \dot q \right) \mathrm{d}t=0$ for all such smooth functions $\de q$ produces the Euler-Lagrange equations of motion: 
\begin{equation} \label{eq_euler_lagrange}
\pp{L}{q} - \dd{}{t} \pp{L}{\dot q} =0. 
\end{equation}
Recall that a Lie group is a smooth manifold which is also a group and for which the group operations of multiplication and inversion are smooth functions \cite{Ho2011_pII}. In the case when there is an intrinsic symmetry in the configuration space, in particular when $M=G$, a Lie group, and  when there is an  appropriate invariance of the Lagrangian with respect to $G$, these Euler-Lagrange equations, defined on the group's tangent bundle $TG$ (i.e. depending on both $g$ and $\dot g$), are cumbersome to use. 

\paragraph{Free Rigid Body}
For example, consider the case of a rigid body rotating about a fixed point with no external torques, so that $G=SO(3)$, $g=\Lambda \in SO(3)=G$, and the Lagrangian is $L\left(\Lambda,\dot \Lambda\right)$. This mechanical system is called a free rigid body. The Euler-Lagrange equations are 
\begin{equation} 
\pp{L}{\Lambda} - \dd{}{t} \pp{L}{\dot \Lambda} =\mathbf{0}, \quad \Lambda^\mathsf{T} \Lambda=I, 
\label{cumber}
\end{equation}
where $I  \in \mathbb{R}^{3 \times 3}$ is the $3 \times 3$ identity matrix. Equation \eqref{cumber} involves $9$ ordinary differential equations with 6  algebraic (i.e. non-differential) constraints, and \eqref{cumber} is highly counterintuitive to use. Euler devised a  description of reduced rigid body motion in terms of the body angular velocity. A more rigorous way to represent this reduction idea is through the Euler-Poincar\'e description of motion \cite{poincare1901forme}, or Euler-Poincar\'e's method. Assuming that the Lagrangian is invariant with respect to rotations on the left, which corresponds to the description of the equations of motion in the body frame, the symmetry-reduced Lagrangian should be of the form $\ell \left(\Lambda^{-1} \dot \Lambda \right)$.

Since $\Lambda \in SO(3)$, $\Lambda^{-1} \Lambda = I$ and $\Lambda^{-1} = \Lambda^\mathsf{T}$, so that
\begin{equation}
\begin{split}
{ \left(\Lambda^{-1} \Lambda \right) }^\cdot &= \Lambda^{-1} \dot \Lambda +  \left(\Lambda^{-1}\right)^\cdot \Lambda 
= \Lambda^{-1} \dot \Lambda +  \left(\Lambda^\mathsf{T} \right)^\cdot \Lambda 
= \Lambda^{-1} \dot \Lambda + {\dot \Lambda}^\mathsf{T} \Lambda 
= \Lambda^{-1} \dot \Lambda + \left(\Lambda^\mathsf{T} \dot \Lambda \right)^\mathsf{T} \\
&= \Lambda^{-1} \dot \Lambda + \left(\Lambda^{-1} \dot \Lambda \right)^\mathsf{T} = \mathbf{0}.
\end{split}
\end{equation}
Hence $\Lambda^{-1} \dot \Lambda = -\left(\Lambda^{-1} \dot \Lambda \right)^\mathsf{T} $, and so $\Lambda^{-1} \dot \Lambda \in \mso(3)$; moreover, $\Lambda^{-1} \delta \Lambda \in \mso(3)$. The isomorphic mapping from the column vectors in $\mathbb{R}^3$ to the Lie algebra $\mso(3)$, i.e. skew-symmetric matrices, is defined using the hat map ${\textvisiblespace}^\wedge : \mathbb{R}^3 \to \mso(3)$ as 
\begin{equation} \label{eq11}
\widehat{\boldsymbol\omega}  
= \begin{bmatrix}
\omega_1 \\
\omega_2 \\
\omega_3
\end{bmatrix}^\wedge
=  \begin{bmatrix}
0 & -\omega_3 & \omega_2 \\
\omega_3 & 0 & -\omega_1 \\
-\omega_2 & \omega_1 & 0
\end{bmatrix},
\end{equation} 
and the inverse mapping from $\mso(3)$ to the column vectors in $\mathbb{R}^3$  is defined using the caron map ${\textvisiblespace}^\vee : \mso(3) \to \mathbb{R}^3$ as
\begin{equation} \label{eq12}
{\begin{bmatrix}
	0 & -\omega_3 & \omega_2 \\
	\omega_3 & 0 & -\omega_1 \\
	-\omega_2 & \omega_1 & 0
	\end{bmatrix}}^\vee 
= \begin{bmatrix}
\omega_1 \\
\omega_2 \\
\omega_3
\end{bmatrix}
= \boldsymbol\omega.  
\end{equation}  
Since the hat map ${\textvisiblespace}^\wedge : \mathbb{R}^3 \to \mso(3)$ and its inverse ${\textvisiblespace}^\vee : \mso(3) \to \mathbb{R}^3$ give isomorphisms between $\mso(3)$ and $\mathbb{R}^3$ and since $\Lambda^{-1} \dot \Lambda \in \mso(3)$, the symmetry-reduced Lagrangian should also be of the form $\ell \left( \bOm \right)$, where $\bOm \equiv \left[ \Lambda^{-1} \dot \Lambda \right]^\vee \in \mathbb{R}^3$. The variation of $\bOm$ is computed as follows \cite{Ho2011_pII}: 
\begin{equation}
\label{sigma_dot_eq}
\begin{split}
\delta \bOm &
%\left( \delta \Lambda^{-1} \dot{\Lambda } + \Lambda^{-1} \delta \dot{\Lambda } \right)^\vee 
%=  \left(-\Lambda^{-1} \delta \Lambda \Lambda^{-1} \dot{\Lambda } + \left( \dot \bSigma \right)^\wedge - \left({\Lambda}^{-1}\right)^\cdot \delta \Lambda  \right)^\vee \\
%&=  \left(-\Lambda^{-1} \delta \Lambda \Lambda^{-1} \dot{\Lambda } + \left( \dot \bSigma \right)^\wedge + \Lambda^{-1} \dot \Lambda \Lambda^{-1}  \delta \Lambda  \right)^\vee 
%=  \left(-\widehat \bSigma \widehat \bOm +\left( \dot \bSigma \right)^\wedge + \widehat \bOm \widehat \bSigma  \right)^\vee \\
= \dot \bSigma + \left( \widehat \bOm \widehat \bSigma-\widehat \bSigma \widehat \bOm \right)^\vee 
= \dot \bSigma + \bOm \times \bSigma, 
\end{split}
\end{equation}
where $\bSigma \equiv \left( \Lambda^{-1} \delta \Lambda \right)^\vee \in \mathbb{R}^3$. Under the hat map isomorphism, the variations $\bSigma$ lie in the Lie algebra $\mso(3)$.
Taking the variation of the action integral, pushing the variational derivative inside the integral, integrating by parts, and enforcing the endpoint conditions $\bSigma(a) =\bSigma(b) = \mathbf{0}$ yields
\begin{equation}
\begin{split} 
\de \int_a^b \ell \left(\bOm \right) \mathrm{d} t &= \int_a^b \de \ell \left(\bOm \right) \mathrm{d} t 
= \int_a^b \left< \dede {\ell}{\bOm}, \de \bOm \right> \mathrm{d} t 
= \int_a^b \left< \dede {\ell}{\bOm},  \dot \bSigma + \bOm \times \bSigma \right> \mathrm{d} t \\
&= -\int_a^b \left< \left(\frac{\mathrm{d}}{\mathrm{d}t}+\bOm \times  \right) \dede {\ell}{\bOm},  \bSigma \right> \mathrm{d} t+ \left. \left< \dede {\ell}{\bOm}, \bSigma \right> \right|_a^b 
= -\int_a^b \left< \left(\frac{\mathrm{d}}{\mathrm{d}t}+\bOm \times  \right) \dede {\ell}{\bOm},  \bSigma \right> \mathrm{d} t.
\end{split}
\end{equation}
Insisting that $\de \int_a^b \ell \left(\bOm \right) \mathrm{d} t = 0$ for all smooth variations $\bSigma$ that vanish at the endpoints generates the well-known equations of motion for the free rigid body: 
\begin{equation} 
\frac{\mathrm{d}}{\mathrm{d} t} \dede{\ell}{\bOm} + \bOm \times  \dede{\ell}{\bOm} =\mathbf{0}. 
\label{EP_rigid_body} 
\end{equation} 
Note that in the above derivation, the functional derivative notation $\dede{\ell}{\bOm}$ is used rather than the partial derivative notation $\pp{\ell}{\bOm}$. The former is used if the Lagrangian depends functionally (e.g. involving a derivative or integral) rather than pointwise on its argument. If the Lagrangian depends only pointwise on its argument, such as is the case for the free rigid body and heavy top (to be discussed next), the two notations agree. For the free rigid body, the symmetry-reduced Lagrangian is $l\left(\bOm\right)=\frac{1}{2}\left<\mathbb{I}\bOm,\bOm\right>$, $\dede{\ell}{\bOm} = \mathbb{I}\bOm$, and the equations of motion \eqref{EP_rigid_body} become 
\begin{equation} \label{EP_fin_rigid_body} 
\dot \bOm = {\mathbb{I}}^{-1} \left[ \left( \mathbb{I}\bOm \right) \times \bOm \right].
\end{equation} 
By multiplying \eqref{EP_rigid_body} by $\Lambda$ and using the identity $\dot \Lambda = \Lambda \widehat{\bOm}$, the equations of motion for the free rigid body may be expressed in conservation law form:
\begin{equation}
\dd{}{t} \left[ \Lambda \dede{\ell}{\bOm} \right] =\mathbf{0} \Leftrightarrow \Lambda \dede{\ell}{\bOm}={\rm const}.
\end{equation}

\paragraph{Heavy Top}
As another application of Euler-Poincar\'e's method, consider the heavy top, which is a rigid body of mass $m$ rotating with a fixed point of support in a uniform gravitational field with gravitational acceleration $g$. Let $\bchi$ denote the vector in the body frame from the fixed point of support to the heavy top's center of mass. To compute the equations of motion for the heavy top, another advected variable $\bGam\equiv\Lambda^{-1} \mathbf{e}_3$ must be introduced. $\bGam$ represents the motion of the unit vector $\mathbf{e}_3$ along the spatial vertical axis, as seen from the body frame. Computing the time and variational derivatives of $\bGamma$ yields 
\begin{equation}
\dot \bGamma = \left(\Lambda^{-1} \mathbf{e}_3 \right)^\cdot=-\Lambda^{-1} \dot \Lambda \Lambda^{-1} \mathbf{e}_3=-\widehat{\bOm}\bGamma=\bGamma \times \bOm 
\end{equation}
and 
\begin{equation} \label{eq_de_Gamma}
\de \bGam= \de \left(\Lambda^{-1} \mathbf{e}_3 \right)=-\Lambda^{-1} \de \Lambda \Lambda^{-1} \mathbf{e}_3=-\widehat{\bSigma}\bGam=\bGam \times \bSigma.
\end{equation} 
The heavy top's reduced Lagrangian is $l\left(\bOm,\bGamma \right)=\frac{1}{2}\left<\mathbb{I}\bOm,\bOm\right>-\left<mg\bchi,\bGamma \right>$. Taking the variation of the action integral, pushing the variational derivative inside the integral, integrating by parts, and enforcing the endpoint conditions $\bSigma(a) = \bSigma(b) = \mathbf{0}$ yields
\begin{equation}
\label{heavy_top_derive}
\begin{split} 
\de \int_a^b l\left(\bOm,\bGamma \right) \mathrm{d} t = \int_a^b \de  l\left(\bOm,\bGamma \right) \mathrm{d} t 
&= \int_a^b \left[ \left<\mathbb{I}\bOm,\de \bOm\right>-\left<mg\bchi,\de \bGamma \right> \right] \mathrm{d} t 
\rem{ %%%BEGIN REM 
\\
&= \int_a^b \left[ \left<\mathbb{I}\bOm,\dot \bSigma + \bOm \times \bSigma \right>-\left<mg\bchi,\bGamma \times \bSigma \right> \right] \mathrm{d} t \\
&= \int_a^b \left<-\dd{}{t} \left( \mathbb{I}\bOm \right) + \left( \mathbb{I}\bOm \right) \times \bOm+mg\bGamma \times \bchi,\bSigma \right> \mathrm{d} t + \left. \left<\mathbb{I}\bOm,\bSigma \right> \right|_a^b 
} %%%END REM 
\\
&= \int_a^b \left<-\dd{}{t} \left( \mathbb{I}\bOm \right) + \left( \mathbb{I}\bOm \right) \times \bOm+mg\bGamma \times \bchi,\bSigma \right> \mathrm{d} t. 
\end{split}
\end{equation}
Insisting that $\de \int_a^b l\left(\bOm,\bGamma \right) \mathrm{d} t = 0$ for all smooth variations $\bSigma$ that vanish at the endpoints generates the equations of motion for the heavy top:
\begin{equation} \label{eq_heavy_top}
\begin{split} 
\dot \bOm &= {\mathbb{I}}^{-1} \left[ \left( \mathbb{I}\bOm \right) \times \bOm+mg\bGamma \times \bchi \right], \\
\dot \bGamma &= \bGamma \times \bOm .
\end{split}
\end{equation}

\paragraph{Adjoint and Coadjoint Operations}
In order to consider  mechanics on general groups, adjoint and coadjoint operations are defined as follows. Consider a Lie group $G$ with Lie algebra $\mathfrak{g}$, dual Lie algebra $\mathfrak{g}^*$, and a pairing $\left<\cdot,\cdot\right>: \mathfrak{g}^* \times \mathfrak{g} \to \mathbb{R} $. The ADjoint operation ${\rm AD} : G \times G \to G$ is defined by 
\begin{equation}
{\rm AD}_g h = g h g^{-1} \quad \forall g,h \in G.
\end{equation}
The Adjoint operation ${\rm Ad} : G \times \mathfrak{g} \to \mathfrak{g}$ is defined by taking a smooth curve $h(t)$ with $h(0)=e$ and $\dot h(0)=\eta \in \mathfrak{g}$ (arbitrary and fixed) and computing 
\begin{equation}
{\rm Ad}_g \eta := \left. \dd{}{t}\right|_{t=0} {\rm AD}_g h (t) = g \eta g^{-1} \quad \forall g \in G, \quad \forall \eta \in \mathfrak{g}. 
\end{equation}
The adjoint operation ${\rm ad} : \mathfrak{g} \times \mathfrak{g} \to \mathfrak{g}$ is defined by taking a smooth curve $g(t)$ with $g(0)=e$ and $\dot g(0)=\xi \in \mathfrak{g}$ (arbitrary and fixed)  and computing 
\begin{equation}
{\rm ad}_\xi \eta := \left. \dd{}{t}\right|_{t=0} {\rm Ad}_{g(t)} \eta   = \xi \eta - \eta \xi = \left[\xi,\eta \right] \quad \forall \xi,\eta \in \mathfrak{g}, 
\end{equation}
where $\left[\cdot,\cdot\right]: \mathfrak{g} \times \mathfrak{g} \to \mathbb{R} $  is the Lie bracket defined by 
\begin{equation}
\left[\xi,\eta \right] = \xi \eta - \eta \xi \quad \forall \xi,\eta \in \mathfrak{g}.
\end{equation}
The coAdjoint operation ${\rm Ad}^* : G \times \mathfrak{g}^* \to \mathfrak{g}^*$ is defined by 
\begin{equation}
\left<{\rm Ad}^*_g \mu , \eta \right> = \left< \mu , {\rm Ad}_g \eta \right> \quad \forall g \in G, \quad \forall \mu \in \mathfrak{g}^*, \quad \forall \eta \in \mathfrak{g}. 
\end{equation}
The coadjoint operation ${\rm ad}^* : \mathfrak{g} \times \mathfrak{g}^* \to \mathfrak{g}^*$ is defined by 
\begin{equation} \label{eq_coad_def}
\left<{\rm ad}^*_\xi \mu , \eta \right> = \left< \mu , {\rm ad}_\xi \eta \right> \quad \forall \xi,\eta \in \mathfrak{g}, \quad \forall \mu \in \mathfrak{g}^*.
\end{equation}

\paragraph{Euler-Poincar\'e's Method}
More generally, if the Lagrangian $L : TG \to \mathbb{R}$ is left-invariant, i.e. 
$L\left(hg,h \dot g\right)=L\left(g,\dot g\right) \quad \forall \left(g,\dot g \right) \in TG, \quad \forall h \in G$, we can define the \emph{symmetry-reduced Lagrangian} through the symmetry reduction $\ell=\ell\left(g^{-1} \dot g \right)=\ell(\xi)=L\left( e,\xi \right),$ where $\xi \equiv g^{-1} \dot g$. Then, the equations of motion \eqref{eq_euler_lagrange} are equivalent to the  Euler-Poincar\'e equations  of motion  obtained from the variational principle 
\begin{equation} 
\label{EPvar} 
\de \int_a^b \ell(\xi) \mbox{d} t =0, \quad \mbox{for  variations} \quad \de \xi= \dot \eta + {\rm ad}_\xi \eta, \quad \forall \eta(t):\, \eta(a)=\eta(b)=0. 
\end{equation} 
The variations $\eta(t)$, assumed to be sufficiently smooth, are sometimes called \emph{free} variations. Applying the variational principle \eqref{EPvar} gives
\begin{equation}  
\begin{split}
\de \int_a^b \ell(\xi) \mbox{d} t = \int_a^b \left< \dede{\ell}{\xi}, \de \xi \right> \mbox{d} t &= \int_a^b \left< \dede{\ell}{\xi}, \dot \eta + {\rm ad}_\xi \eta \right> \mbox{d} t
% \\ &= \int_a^b \left<-\dd{}{t} \dede{\ell}{\xi} + {\rm ad}^*_ \xi  \dede{\ell}{\xi},\eta \right> \mbox{d} t + \left. \left< \dede{\ell}{\xi} , \eta \right> \right|_a^b 
= \int_a^b \left<-\dd{}{t} \dede{\ell}{\xi} + {\rm ad}^*_ \xi  \dede{\ell}{\xi},\eta \right> \mbox{d} t  =   0,
\end{split}
\end{equation} 
which yields the Euler-Poincar\'e equations of motion:
\begin{equation}
\dd{}{t} \dede{\ell}{\xi} - {\rm ad}^*_ \xi  \dede{\ell}{\xi} =0. 
\label{EPeq_left}
\end{equation} 
For right-invariant Lagrangians, i.e. $L\left(gh,\dot gh\right)=L\left(g,\dot g\right) \quad \forall h \in G$, the Euler-Poincar\'e equations of motion \eqref{EPeq_left}  change by altering the sign in front of ${\rm ad}^*_\xi$ from minus to plus. For the free rigid body, $\bxi = \bOm$, $l\left(\bOm\right)=\frac{1}{2}\left<\mathbb{I}\bOm,\bOm\right>$, $\dede{\ell}{\bOm} = \mathbb{I}\bOm$, and ${\rm ad}^*_{\bOm} \dede{\ell}{\bOm} = - \bOm \times \dede{\ell}{\bOm} = \left(\mathbb{I}\bOm\right) \times \bOm$, so that the free rigid body equations of motion \eqref{EP_fin_rigid_body} derived earlier agree with the Euler-Poincar\'e equations of motion \eqref{EPeq_left}.

It is interesting that 
 \eqref{EPeq_left} implies the conservation of angular momentum. Indeed,  letting $\alpha \in \mathfrak{g}$ be arbitrary and constant in time and letting $t_0 \in \mathbb{R}$ be an arbitrary time, one can derive  that
\begin{equation} \label{eq_con_ang_mom}
\begin{split}
\left<\left.\dd{}{t} \right|_{t=t_0} {\rm Ad}^*_{g^{-1}} \dede{\ell}{\xi},  \alpha \right>&=0. 
\rem{%%%BEGIN REM 
\left.\dd{}{t}\right|_{t_0} \left<  {\rm Ad}^*_{g^{-1}} \dede{\ell}{\xi} , \alpha \right>= 
\left.\dd{}{t}\right|_{t_0} \left< \dede{\ell}{\xi} , {\rm Ad}_{g^{-1}} \alpha \right>\\&= \left< \left.\dd{}{t}\right|_{t_0} \dede{\ell}{\xi} , {\rm Ad}_{g^{-1}} \alpha \right>+\left< \dede{\ell}{\xi} , \left.\dd{}{t}\right|_{t_0} {\rm Ad}_{g^{-1}} \alpha \right>\\
&= \left< {\rm Ad}_{g^{-1}}^* \left[ \left.\dd{}{t}\right|_{t_0} \dede{\ell}{\xi} \right] , \alpha \right>+\left< \dede{\ell}{\xi} , -{\rm ad}_\xi \left[ {\rm Ad}_{g^{-1}} \alpha \right] \right>\\
&= \left< {\rm Ad}_{g^{-1}}^* \left[ \left.\dd{}{t}\right|_{t_0} \dede{\ell}{\xi} \right] , \alpha \right>-\left<{\rm ad}_\xi^* \dede{\ell}{\xi} , {\rm Ad}_{g^{-1}} \alpha \right> \\
&= \left< {\rm Ad}_{g^{-1}}^* \left[ \left.\dd{}{t}\right|_{t_0} \dede{\ell}{\xi} \right] , \alpha \right>-\left<{\rm Ad}_{g^{-1}}^* \left[ {\rm ad}_\xi^* \dede{\ell}{\xi} \right] , \alpha \right>\\
&= \left< {\rm Ad}_{g^{-1}}^* \left[ \left.\dd{}{t}\right|_{t_0} \dede{\ell}{\xi}-{\rm ad}_\xi^* \dede{\ell}{\xi} \right] , \alpha \right> = \left< {\rm Ad}_{g^{-1}}^* 0 , \alpha \right> = 0,
} %%%END REM 
\end{split}
\end{equation}
\rem{ %%%BEGIN REM 
where \eqref{EPeq_left} is used in the second to last equality. In the fourth equality of \eqref{eq_con_ang_mom}, the following result is used
\begin{equation} 
\begin{split}
\left.\dd{}{t}\right|_{t_0} {\rm Ad}_{g^{-1}} \alpha &= \left.\dd{}{t}\right|_{t=t_0} {\rm Ad}_{g(t)^{-1}} \alpha = \left.\dd{}{t}\right|_{t=t_0} {\rm Ad}_{g(t)^{-1}g(t_0)} \left[{\rm Ad}_{g(t_0)^{-1}} \alpha \right] \\
&= \left.\dd{}{t}\right|_{t=t_0} \left\{ g(t)^{-1}g(t_0) \left[{\rm Ad}_{g(t_0)^{-1}} \alpha \right] g(t_0)^{-1}g(t) \right\}   \\
&= \bigg\{  -g(t)^{-1} \dot g(t) g(t)^{-1} g(t_0) \left[{\rm Ad}_{g(t_0)^{-1}} \alpha \right] g(t_0)^{-1}g(t) \\
&\hphantom{=} + g(t)^{-1}g(t_0) \left[{\rm Ad}_{g(t_0)^{-1}} \alpha \right] g(t_0)^{-1}\dot g(t) \left. \bigg\} \right|_{t=t_0}  \\
&= -g(t_0)^{-1} \dot g(t_0) \left[{\rm Ad}_{g(t_0)^{-1}} \alpha \right]  + \left[{\rm Ad}_{g(t_0)^{-1}} \alpha \right] g(t_0)^{-1}\dot g(t_0) \\
&= -\xi(t_0) \left[{\rm Ad}_{g(t_0)^{-1}} \alpha \right]  + \left[{\rm Ad}_{g(t_0)^{-1}} \alpha \right] \xi(t_0) \\
&= -{\rm ad}_{\xi(t_0)} \left[{\rm Ad}_{g(t_0)^{-1}} \alpha \right],
\end{split}
\end{equation}
using the property ${\rm Ad}_g {\rm Ad}_h \eta = g \left(h \eta h^{-1} \right) g^{-1}= \left(gh \right) \eta \left(g h \right)^{-1}={\rm Ad}_{\left(gh\right)} \eta \quad \forall g,h \in G, \quad \forall \eta \in \mathfrak{g}$ in the second equality. Since $\alpha \in \mathfrak{g}$ is arbitrary, \eqref{eq_con_ang_mom} implies that $\left.\dd{}{t}\right|_{t_0}  {\rm Ad}^*_{g^{-1}} \dede{\ell}{\xi}=0$. Since $t_0 \in \mathbb{R}$ is an arbitrary time, $\dd{}{t} {\rm Ad}^*_{g^{-1}} \dede{\ell}{\xi}=0$, thereby proving conservation of angular momentum.  
} %%%END REM 

\rem{%%%%BEGIN REM 
\paragraph{Hamilton-Pontryagin's Principle} An alternative to Euler-Poincar\'e's method is Hamilton-Pontryagin's principle, which says that the equations of motion may be obtained from the variational principle
\begin{equation}
\de \hat{S}\left(\xi,g,\dot g \right) = \de \int_a^b \hat{\ell}\left(\xi,g,\dot g \right) \dt =0
\end{equation}
for all variations of $g$ such that $\de g(a)=\de g(b) = 0$, where $\hat{S}$ is the constrained action integral
\begin{equation}
\hat{S}\left(\xi,g,\dot g \right) = \int_a^b \hat{\ell}\left(\xi,g,\dot g \right) \dt = \int_a^b \left[ \ell(\xi)+\left<\mu,g^{-1} \dot g - \xi\right> \right] \dt
\end{equation}
and $\hat{\ell}$ is the augmented, reduced Lagrangian
\begin{equation}
\hat{\ell}\left(\xi,g,\dot g \right) = \ell(\xi)+\left<\mu,g^{-1} \dot g - \xi\right>.
\end{equation}
To see that Hamilton-Pontryagin's principle gives the same equations of motion as Euler-Poincar\'e's method, define $\eta = g^{-1} \de g$.
Since $\xi=g^{-1} \dot g$,
\begin{equation} \label{eq_de_xi}
\de \xi = \de \left(g^{-1} \dot g \right) = -g^{-1}\de g g^{-1} \dot g + g^{-1} \de \dot g = -\eta \xi + g^{-1} \de \dot g.
\end{equation}
Since $\eta = g^{-1} \de g$,
\begin{equation} \label{eq_dot_eta}
\dot \eta = \left(g^{-1} \de g \right)^\cdot = -g^{-1}\dot g g^{-1} \de g + g^{-1} \left( \de g \right)^\cdot = - \xi \eta + g^{-1} \left( \de g \right)^\cdot.
\end{equation}
Subtracting \eqref{eq_dot_eta} from \eqref{eq_de_xi} gives
\begin{equation} \label{eq_ep_hp}
\de \xi - \dot \eta = \xi \eta -\eta \xi = {\rm ad}_\xi \eta.
\end{equation}
Now compute the variation of $\hat{S}$:
\begin{equation}
\begin{split}
\de \hat{S} = \de \int_a^b \hat{\ell}\left(\xi,g,\dot g \right) \dt &= \de \int_a^b \left[ \ell(\xi)+\left<\mu,g^{-1} \dot g - \xi\right> \right] \dt \\
&= \int_a^b \left[ \left< \dede{\hat{\ell}}{\xi}-\mu,\de \xi \right> + \left<\de \mu,g^{-1} \dot g - \xi  \right>+ \left< \mu, \de \left(g^{-1} \dot g  \right) \right> \right] \dt \\
&= \int_a^b \left[ \left< \dede{\hat{\ell}}{\xi}-\mu,\de \xi \right> + \left<\de \mu,g^{-1} \dot g - \xi  \right>+ \left< \mu, \dot \eta + {\rm ad}_\xi \eta \right> \right] \dt \\
&= \int_a^b \left[ \left< \dede{\hat{\ell}}{\xi}-\mu,\de \xi \right> + \left<\de \mu,g^{-1} \dot g - \xi  \right>+ \left<-\dot \mu+{\rm ad}_\xi^* \mu, \eta  \right> \right] \dt + \left. \left<\mu ,\eta \right>\right|_a^b \\
&= \int_a^b \left[ \left< \dede{\hat{\ell}}{\xi}-\mu,\de \xi \right> + \left<\de \mu,g^{-1} \dot g - \xi  \right>+ \left<-\dot \mu+{\rm ad}_\xi^* \mu, \eta  \right> \right] \dt,
\end{split}
\end{equation}
since $\de \left(g^{-1} \dot g  \right) = \dot \eta + {\rm ad}_\xi \eta$ by \eqref{eq_ep_hp} and because $\eta(a)=\eta(b)=0$ (since $\de g(a)=\de g(b)=0$). $\de \hat{S}=0$ gives the Hamilton-Pontryagin equations of motion
\begin{equation}
\pp{\hat{\ell}}{\xi}=\mu, \quad g^{-1} \dot g - \xi, \quad \dot \mu={\rm ad}_\xi^* \mu,
\end{equation}
which agree with the Euler-Poincar\'e equations of motion \eqref{EPeq_left}. Even though Hamilton-Pontryagin's principle could be utilized, this paper instead relies on Euler-Poincar\'e's method to derive the equations of motion for the rolling ball in Section~\ref{sec_ball_uncontrolled}. }

\paragraph{Euler-Poincar\'e's Method with an Advected Parameter}
In order to further treat the effect of gravity on the heavy top and also on the rolling ball in Section~\ref{sec_ball_uncontrolled}, we let the Lagrangian depend on a parameter (gravity) which is advected with the dynamics. Formally, 
 let $V$ be a vector space. Suppose the Lagrangian $L$ depends on a parameter in the dual space $V^*$, so that the general Lagrangian has the form $L : TG \times V^* \to \mathbb{R}$. For a parameter $\alpha_0 \in V^*$, suppose that the Lagrangian $L_{\alpha_0} : TG \to \mathbb{R}$ defined by  $L_{\alpha_0}\left(g,\dot g \right) = L\left(g,\dot g, \alpha_0 \right) \quad \forall \left(g ,\dot g \right) \in TG$ is left-invariant, i.e. 
$L\left(hg,h \dot g,h \alpha_0 \right)=L \left(g,\dot g, \alpha_0 \right) \quad \forall \left(g,\dot g \right) \in TG, \quad \forall h \in G$. Then we can define the \emph{symmetry-reduced Lagrangian} through the symmetry reduction $\ell=\ell\left(g^{-1} \dot g, g^{-1} \alpha_0   \right)=\ell \left(\xi, \alpha \right)=L\left( e,\xi, \alpha \right),$ where $\xi \equiv g^{-1} \dot g$ and $\alpha \equiv g^{-1} \alpha_0$. Euler-Poincar\'e's method with an advected parameter says that the equations  of motion are  obtained from the variational principle
\begin{equation} 
\label{EPAvar} 
\de \int_a^b \ell \left(\xi,\alpha \right) \mbox{d} t =0, \quad \mbox{for  variations} \quad \de \xi= \dot \eta + {\rm ad}_\xi \eta, \, \de \alpha = -\eta \alpha, \quad \forall \eta(t):\, \eta(a)=\eta(b)=0. 
\end{equation}
Before applying this variational principle, the diamond operation  $\diamond$ is defined. $\diamond : V \times V^* \to \mathfrak{g}^*$ is defined by \begin{equation}
\left< v \diamond w , \xi \right> =  \left< w , \xi v \right> \quad \forall v \in V, \quad \forall w \in V^*, \quad \forall \xi \in \mathfrak{g}.
\end{equation} 
$\diamond : V^* \times V \to \mathfrak{g}^*$ is defined by 
\begin{equation}
\left< w \diamond v , \xi \right> = -\left< v \diamond w , \xi \right> = -\left< w , \xi v \right> \quad \forall v \in V, \quad \forall w \in V^*, \quad \forall \xi \in \mathfrak{g}.
\end{equation}
Applying the variational principle \eqref{EPAvar} gives
\begin{equation}  
\begin{split}
\de \int_a^b \ell \left(\xi,\alpha \right) \mbox{d} t &= \int_a^b \left[ \left< \dede{\ell}{\xi}, \de \xi \right>+ \left< \dede{\ell}{\alpha}, \de \alpha \right> \right] \mbox{d} t \\
&= \int_a^b \left[ \left< \dede{\ell}{\xi}, \dot \eta + {\rm ad}_\xi \eta \right>+\left<\dede{l}{\alpha},-\eta \alpha \right> \right] \mbox{d} t \\ &= \int_a^b \left[ \left<-\dd{}{t} \dede{\ell}{\xi} + {\rm ad}^*_ \xi  \dede{\ell}{\xi},\eta \right>+\left<\dede{l}{\alpha} \diamond \alpha,\eta \right> \right] \mbox{d} t + \left. \left< \dede{\ell}{\xi} , \eta \right> \right|_a^b \\ &= \int_a^b \left<-\dd{}{t} \dede{\ell}{\xi} + {\rm ad}^*_ \xi  \dede{\ell}{\xi}+\dede{l}{\alpha} \diamond \alpha,\eta \right> \mbox{d} t  =   0,
\end{split}
\end{equation} 
which yields the Euler-Poincar\'e equations of motion with an advected parameter:
\begin{equation}
\dd{}{t} \dede{\ell}{\xi} - {\rm ad}^*_ \xi  \dede{\ell}{\xi} -\dede{l}{\alpha} \diamond \alpha =0. 
\label{EPAeq_left}
\end{equation} 
The most direct application of the Euler-Poincar\'e equations of motion with an advected parameter is the heavy top, where the advected parameter is the gravity vector expressed in the heavy top's body frame. For the heavy top, $\bxi=\bOm$, $\alpha = \bGamma$, $l\left(\bOm,\bGamma \right)=\frac{1}{2}\left<\mathbb{I}\bOm,\bOm\right>-\left<mg\bchi,\bGamma \right>$, $\dede{\ell}{\bOm}=\mathbb{I}\bOm$, ${\rm ad}^*_{\bOm} \dede{\ell}{\bOm} = - \bOm \times \mathbb{I}\bOm$, $\dede{\ell}{\bGamma}=-mg\bchi$, and $\dede{l}{\bGamma} \diamond \bGamma=-mg\bchi \times \bGamma $. Plugging all these identities into \eqref{EPAeq_left} recovers the previously derived heavy top equations of motion \eqref{eq_heavy_top}.

\subsection{Nonholonomic Constraints and Lagrange-d'Alembert's Principle} \label{ssec_LdA}
Suppose a mechanical system having configuration space $M$, a manifold of dimension $n$, must satisfy $m < n$  constraints that \emph{are linear in velocity}. To express these velocity constraints formally, the notion of a distribution is needed. Given the manifold $M$, a distribution $\mathcal{D}$ on $M$  is a subset of the tangent bundle $TM = \bigcup_{q \in M} T_q M$: $\mathcal{D} = \bigcup_{q \in M} \mathcal{D}_q$, where $\mathcal{D}_q \subset T_q M$ and $m = \mathrm{dim} \, \mathcal{D}_q < \mathrm{dim} \, T_q M = n$ for each $q \in M$. A curve $q(t) \in M$ satisfies the constraints if $\dot q(t) \in \mathcal{D}_{q(t)}$. Lagrange-d'Alembert's principle states that the equations of motion are determined by 
\begin{equation} 
\de \int_a^b L(q, \dot q) \mathrm{d} t=0 \Leftrightarrow \int_a^b \left[ \dd{}{t} \pp{L}{\dot q}- \pp{L}{q}  \right] \de q \, \mbox{d} t = 0   \Leftrightarrow  \dd{}{t} \pp{L}{\dot q}- \pp{L}{q} \in \mathcal{D}_q^\circ 
\label{LdA0}
\end{equation} 
for all smooth variations $\de q(t)$ of the curve $q(t)$ such that $\de q(t) \in \mathcal{D}_{q(t)}$ for all $a\le t \le b$ and such that $\de q(a)=\de q(b)=0$, and for which $\dot q(t) \in \mathcal{D}_{q(t)}$ for all $a\le t \le b$. If one writes the nonholonomic constraint in local coordinates as $\sum_{i=1}^n A(q)^j_i \dot q^i=0$, $j=1, \ldots, m < n$, then \eqref{LdA0} is written in local coordinates as 
\begin{equation} 
\dd{}{t} \pp{L}{{\dot q}^i}- \pp{L}{q^i} = \sum_{j=1}^m \lambda_j A(q)^j_i \, , \quad i=1,\ldots,n \, , \quad \sum_{i=1}^n A(q)^j_i {\dot q}^i=0  , 
\label{LdA}
\end{equation}
where the $\lambda_j$ are Lagrange multipliers enforcing $\sum_{i=1}^n A(q)^j_i {\de q}^i=0$, $j=1, \ldots, m$. Aside from Lagrange-d'Alembert's approach, there is also an alternative \emph{vakonomic} approach to derive the equations of motion for nonholonomic mechanical systems. Simply speaking, the vakonomic approach relies on substituting the constraint into the Lagrangian before taking variations or, equivalently, enforcing the constraints using the appropriate Lagrange multiplier method \revision{R1Q5}{\cite{kozlov1982dynamics1,kozlov1982dynamics2}}. In general, it is an experimental fact that all known nonholonomic mechanical systems obey the equations of motion resulting from Lagrange-d'Alembert's principle \cite{lewis1995variational}.

\paragraph{Suslov's Problem} To illustrate Lagrange-d'Alembert's principle in conjunction with Euler-Poincar\'e's method, also known as  \emph{Euler-Poincar\'e-Suslov's method}, consider a rigid body rotating about a fixed point such that its body angular velocity $\bOm$ must be orthogonal to a prescribed body frame vector $\bxi$. Such a rigid body is called Suslov's problem in honor of the Russian mathematician who introduced and studied it in 1902 \cite{suslov1946theoretical}. Mathematically, the constraint for Suslov's problem is $\left<\bOm,\bxi\right>=0$, so that Suslov's problem is an algebraically simple example of a nonholonomic mechanical system. In Suslov's original formulation \cite{suslov1946theoretical}, $\bxi$ was assumed to be fixed in the body frame. In \cite{putkaradze2018constraint} and here, $\bxi$ is permitted to vary with time. The Lagrangian for Suslov's problem is its kinetic energy, so that the symmetry-reduced Lagrangian is $\ell(\bOm)=\frac{1}{2}\left<\inertia \bOm,\bOm\right>$ and the action integral is $S = \int_a^b \ell(\bOm) \dt = \int_a^b \frac{1}{2}\left<\inertia \bOm,\bOm\right> \dt$. Since $\bOm \equiv \left[ \Lambda^{-1} \dot \Lambda \right]^\vee$, according to \eqref{sigma_dot_eq}, $\de \bOm = \dot \bSigma + \bOm \times \bSigma$ where $\bSigma \equiv \left( \Lambda^{-1} \delta \Lambda \right)^\vee$. Part of Lagrange-d'Alembert's principle states that the nonholonomic constraint $\left<\bOm,\bxi\right>=0$ implies that the variations $\bSigma$ must satisfy $\left< \bSigma , \bxi \right>=0$ when deriving the equations of motion. Enforcing the constraint $\left< \bSigma , \bxi \right>=0$ on the variations $\bSigma$ through the time-varying Lagrange multiplier $\lambda$, Lagrange-d'Alembert's principle in conjunction with Euler-Poincar\'e's method dictate that the equations of motion for Suslov's problem are given by
\begin{equation} \label{eq_suslov_LdA}
0=\de S+\int_a^b \lambda  \left< \bSigma , \bxi \right>  \dt
\end{equation}
for variations $\de \bOm = \dot \bSigma + \bOm \times \bSigma$, for all variations $\bSigma$ such that $\bSigma(a)=\bSigma(b)=\mathbf{0}$, and such that $\left<\bOm,\bxi\right>=0$. Pushing the variational derivative inside the action integral's integration operator, using the fact that $\de \bOm = \dot \bSigma + \bOm \times \bSigma$, integrating by parts, and invoking the vanishing endpoint assumptions $\bSigma(a)=\bSigma(b)=\mathbf{0}$, \eqref{eq_suslov_LdA} simplifies to
\begin{equation} \label{eq_suslov_LdA_red}
\begin{split}
0&=\de S+\int_a^b \lambda  \left< \bSigma , \bxi \right>  \dt 
= \int_a^b \left<\inertia \bOm,\de \bOm\right> \dt+\int_a^b \lambda  \left< \bSigma , \bxi \right>  \dt \\
&= \int_a^b \left<\inertia \bOm,\dot \bSigma + \bOm \times \bSigma\right> \dt+\int_a^b \lambda  \left< \bSigma , \bxi \right>  \dt \\
&= -\int_a^b \left< \left(\frac{\mathrm{d}}{\mathrm{d}t}+\bOm \times  \right) \inertia \bOm,  \bSigma \right> \mathrm{d} t+ \left. \left< \inertia \bOm, \bSigma \right> \right|_a^b+\int_a^b \lambda  \left< \bSigma , \bxi \right> \dt \\ 
&= \int_a^b \left< -\left(\frac{\mathrm{d}}{\mathrm{d}t}+\bOm \times  \right) \inertia \bOm + \lambda \bxi,  \bSigma \right> \dt.
\end{split}
\end{equation}
Since \eqref{eq_suslov_LdA_red} must be satisfied for all variations $\bSigma$ such that $\bSigma(a)=\bSigma(b)=\mathbf{0}$, the equations of motion for Suslov's problem are given by 
\begin{equation} \label{eq_suslov_tmp}
\inertia \dot \bOm = \left( \inertia \bOm \right) \times \bOm + \lambda \bxi,
\end{equation}
where the Lagrange multiplier $\lambda$ is determined from the nonholonomic constraint $\left<\bOm,\bxi\right>=0$. Dotting both sides of \eqref{eq_suslov_tmp} by $\inertia^{-1} \bxi$, solving for $\lambda$, applying the product rule $\left<\bOm,\bxi\right>^\cdot =\left<\dot \bOm,\bxi\right>+\left<\bOm,\dot \bxi\right>$, and invoking the nonholonomic constraint $\left<\bOm,\bxi\right>=0$ yield the formula for the Lagrange multiplier $\lambda$:
\begin{equation} \label{eq_suslov_lambda}
\begin{split}
\lambda &= \frac {\left<\inertia \dot \bOm,\inertia^{-1} \bxi\right> - \left< \left( \inertia \bOm \right) \times \bOm, \inertia^{-1} \bxi \right>} { \left<\bxi, \inertia^{-1} \bxi \right> } = \frac {\left<\dot \bOm,\bxi\right> - \left< \left( \inertia \bOm \right) \times \bOm, \inertia^{-1} \bxi \right>} { \left<\bxi, \inertia^{-1} \bxi \right> } \\
&= \frac {\left<\bOm,\bxi\right>^\cdot - \left< \bOm, \dot  \bxi \right> - \left< \left( \inertia \bOm \right) \times \bOm, \inertia^{-1} \bxi \right>} { \left<\bxi, \inertia^{-1} \bxi \right> } = - \frac { \left< \bOm, \dot  \bxi \right> + \left< \left( \inertia \bOm \right) \times \bOm, \inertia^{-1} \bxi \right>} { \left<\bxi, \inertia^{-1} \bxi \right> },
\end{split}
\end{equation}
so that the equations of motion for Suslov's problem are
\begin{equation} \label{eq_suslov}
\inertia \dot \bOm = \left( \inertia \bOm \right) \times \bOm - \frac { \left< \bOm, \dot  \bxi \right> + \left< \left( \inertia \bOm \right) \times \bOm, \inertia^{-1} \bxi \right>} { \left<\bxi, \inertia^{-1} \bxi \right> } \bxi.
\end{equation}
The reader is referred to \cite{putkaradze2018constraint} for further details.

\rem{%%%BEGIN REM 
\subsection{A Simple Nonholonomic Particle} \label{ssec_bates}
Consider a free particle with unit mass moving in space subject to the nonholonomic constraint
\begin{equation} \label{eq_part_constraint}
\dot z = y \dot x.
\end{equation}
This problem was introduced and studied by Bates and \'Sniatycki in \cite{bates1993nonholonomic}. The particle's Lagrangian is 
\begin{equation} \label{eq_part_lagrangian}
l = \frac{1}{2}\left({\dot x}^2 +{\dot y}^2 + {\dot z}^2\right),
\end{equation}
the particle's action integral is
\begin{equation}
S = \int_a^b l \mathrm{d}t = \int_a^b \frac{1}{2}\left({\dot x}^2 +{\dot y}^2 + {\dot z}^2\right) \mathrm{d}t,
\end{equation}
and the variation of the particle's action integral is
\begin{equation} 
\de S = \de \int_a^b l \mathrm{d}t = \int_a^b \de l \mathrm{d}t = \int_a^b \left[ \dot x \de {\dot x} +  \dot y \de {\dot y} +  \dot z \de {\dot z}\right] \mathrm{d}t 
= -\int_a^b \left[ \ddot x \de  x +  \ddot y \de y +  \ddot z \de z\right] \mathrm{d}t.
\end{equation} 

\paragraph{Lagrange-d'Alembert's Approach}
In Lagrange-d'Alembert's approach, the constraint $\dot z = y \dot x$ implies the variational constraint $\de z = y \de x$. Lagrange-d'Alembert's principle states that the equations of motion are given by $\de S = 0$ for all variations $\de x$, $\de y$, and $\de z$ (i.e. Hamilton's principle) subject to the variational constraint $\de z = y \de x$ and the original constraint $\dot z = y \dot x$ \cite{Bloch2003}. Using the method of Lagrange multipliers to simultaneously enforce the conditions $\de S = 0$ and $\de z = y \de x$, the equations of motion must satisfy 
\begin{equation}
\begin{split}
0 &= \de S + \int_a^b \lambda \left[ \de z - y \de x \right] \mathrm{d}t \\
&=  -\int_a^b \left[ \ddot x \de  x +  \ddot y \de y +  \ddot z \de z\right] \mathrm{d}t + \int_a^b \lambda \left[ \de z - y \de x \right] \mathrm{d}t  \\
&= -\int_a^b \left[\left(\ddot x+\lambda y \right)\de x+ \ddot y \de y + \left(\ddot z - \lambda \right) \de z  \right] \mathrm{d}t, \\
\dot z &= y \dot x,
\end{split}
\end{equation}
for all variations $\de x$, $\de y$, and $\de z$ and where $\lambda$ is a Lagrange multiplier. Thus, the equations of motion are
\begin{equation}
\begin{split}
\ddot x+\lambda y &= 0, \\
\ddot y &= 0, \\
\ddot z - \lambda &= 0, \\
\dot z &= y \dot x.
\end{split}
\end{equation}
These equations may be simplified by using the equation $\ddot z - \lambda = 0$ to eliminate $\lambda$ and by then using the original constraint $\dot z = y \dot x$ to eliminate $z$; after these simplifications the equations of motion only depend on $x$ and $y$. Since $\ddot z - \lambda = 0$, $\ddot z = \lambda$ and so the equation $\ddot x+\lambda y = 0$ becomes $\ddot x+\ddot z y = 0$. Having eliminated $\lambda$, the equations of motion become
\begin{equation} \label{eq_part_LdA}
\begin{split}
\ddot x+\ddot z y &= 0, \\
\ddot y &= 0, \\
\dot z &= y \dot x.
\end{split}
\end{equation}
The original constraint $\dot z = y \dot x$ implies that $\ddot z = \dot y \dot x + y \ddot x$. Substituting $\ddot z = \dot y \dot x + y \ddot x$ into $\ddot x+\ddot z y = 0$ yields $\ddot x+\left(\dot y \dot x + y \ddot x \right) y = 0$.
The equations of motion simplify to
\begin{equation}
\begin{split}
\ddot x + \frac{y}{1+y^2}\dot x \dot y &=0, \\
\ddot y &= 0.
\end{split}
\end{equation}

For this simple problem, a more direct derivation of the above equations of motion, avoiding introduction of the Lagrange multiplier, is achieved by substituting the variational constraint $\de z = y \de x$ in for $\de z$ in the variation of the action integral. Making this substitution, the variation of the action integral becomes 
\begin{equation}
\de S = -\int_a^b \left[ \ddot x \de  x +  \ddot y \de y +  \ddot z \de z\right] \mathrm{d}t = -\int_a^b \left[ \ddot x \de  x +  \ddot y \de y +  \ddot zy \de x\right] \mathrm{d}t = -\int_a^b \left[ (\ddot x+\ddot z y) \de  x +  \ddot y \de y \right] \mathrm{d}t.
\end{equation}
Applying Hamilton's principle (i.e. demanding that $\de S=0$ for all variations $\de x$ and $\de y$) and using the original constraint $\dot z = y \dot x$ recovers the previously obtained equations of motion
\begin{equation} 
\begin{split}
\ddot x + \frac{y}{1+y^2}\dot x \dot y &=0, \\
\ddot y &= 0.
\end{split}
\end{equation}
These equations of motion may be solved analytically. Since $\ddot y = 0$, 
\begin{equation} \label{eq_part_LdA_y}
y(t) = ct+d,
\end{equation}
for integration constants $c$ and $d$. Using this result, the equation $\ddot x + \frac{y}{1+y^2}\dot x \dot y =0$ becomes 
\begin{equation}
\ddot x + \frac{ct+d}{1+\left(ct+d \right)^2}c \dot x =0.
\end{equation}
If $c = 0$, 
\begin{equation}
x(t) = \alpha t+\beta, 
\end{equation}
for integration constants $\alpha$ and $\beta$. If $c \ne 0$, 
\begin{equation} \label{eq_part_LdA_x}
x(t) = B \ln\left[ct+d + \sqrt{1+\left(ct+d \right)^2} \right]+E,
\end{equation}
for integration constants $B$ and $E$. If the particle's initial conditions $x(a)$, $\dot x(a)$, $y(a)$, and $\dot y(a)$ are given at time $t=a$, then the integration constants may be readily determined. $c=\dot y(a)$ and $d=y(a)-\dot y(a) a$. If $c = 0$, $\alpha=\dot x(a)$ and $\beta=x(a)-\dot x(a) a$. If $c \ne 0$, $B = \frac{\dot x(a) \sqrt{1+{y(a)}^2}}{\dot y(a)}$ and $E=x(a)- \frac{\dot x(a) \sqrt{1+{y(a)}^2}}{\dot y(a)} \ln \left[y(a)+\sqrt{1+{y(a)}^2} \right]$. Since $\dot y = c$ is constant, the Lagrange-d'Alembert solution conserves the $y$-component of momentum, which is contrary to a na{\"i}ve application of Noether's theorem. The Lagrangian \eqref{eq_part_lagrangian} and nonholonomic constraint \eqref{eq_part_constraint} for this particle are invariant under the translational transformations $\brho_1 \left(x,y,z \right)=\left(x+C_1,y,z \right)$ and $\brho_3 \left(x,y,z \right)=\left(x,y,z+C_3 \right)$  for constants $C_1$ and $C_3$. But the nonholonomic constraint \eqref{eq_part_constraint} is not invariant under the translational transformation $\brho_2 \left(x,y,z \right)=\left(x,y+C_2,z \right)$ for a constant $C_2$. Based on these observations, (wrong application of) Noether's theorem would predict that the particle's $x$- and $z$-momenta, $p_x$ and $p_z$, should be conserved and that the particle's $y$-momentum, $p_y$, should not be conserved, which is not in agreement with \eqref{eq_part_LdA_x}, \eqref{eq_part_LdA_y}, and \eqref{eq_part_LdA}. That is, \eqref{eq_part_LdA_x} says that $\dot x$ is nonconstant,  \eqref{eq_part_LdA_y} says that $\dot y=c$ is constant, and \eqref{eq_part_LdA} says that $\dot z$, which equals  $y \dot x$, is nonconstant.

\paragraph{Vakonomic Approach}
In the vakonomic approach, the constraint $\dot z = y \dot x$ is added to the Lagrangian via the method of Lagrange multipliers to obtain the modified Lagrangian
\begin{equation}
\tilde{l}= l+\lambda \left(\dot z - y \dot x\right)=\frac{1}{2}\left({\dot x}^2 +{\dot y}^2 + {\dot z}^2\right)+\lambda \left(\dot z - y \dot x\right),
\end{equation}
for a Lagrange multiplier $\lambda$ \cite{Bloch2003,arnold2007mathematical}. The modified action integral is 
\begin{equation}
\begin{split}
\tilde{S} = \int_a^b \tilde{l} \mathrm{d}t = \int_a^b \left[ \frac{1}{2}\left({\dot x}^2 +{\dot y}^2 + {\dot z}^2\right)+\lambda \left(\dot z - y \dot x\right) \right] \mathrm{d}t &= \int_a^b \frac{1}{2}\left({\dot x}^2 +{\dot y}^2 + {\dot z}^2\right) \mathrm{d}t +\int_a^b \lambda \left(\dot z - y \dot x\right) \mathrm{d}t \\ &= S+\int_a^b \lambda \left(\dot z - y \dot x\right) \mathrm{d}t. 
\end{split}
\end{equation}
The variation of the modified action integral is
\begin{equation}
\begin{split}
\de \tilde{S} = \de \int_a^b \tilde{l} \mathrm{d}t &= \de S+ \de \int_a^b \left[ \lambda \left(\dot z - y \dot x\right) \right] \mathrm{d}t \\
&= \de S+ \int_a^b \de \left[ \lambda \left(\dot z - y \dot x\right) \right] \mathrm{d}t \\
&= \int_a^b \left[\dot x \de \dot x +\dot y \de \dot y + \dot z \de \dot z \right] \mathrm{d}t + \int_a^b \de \left[ \lambda \left(\dot z - y \dot x\right) \right] \mathrm{d}t \\
&= \int_a^b \left[\dot x \de \dot x +\dot y \de \dot y + \dot z \de \dot z + \lambda \left(\de \dot z -y \de \dot x - \de y \dot x \right) + \de \lambda \left(\dot z - y \dot x\right)  \right] \mathrm{d}t \\
&= \int_a^b \left[-\ddot x \de x-\ddot y \de y-\ddot z \de z-\dot \lambda \de z-\lambda \dot x \de y + \left(\dot \lambda y + \lambda \dot y\right) \de x + \left(\dot z - y \dot x\right) \de \lambda \right] \mathrm{d}t \\
&= \int_a^b \left[ \left(-\ddot x+\dot \lambda y + \lambda \dot y \right) \de x - \left(\ddot y + \lambda \dot x \right) \de y - \left(\ddot z+\dot \lambda  \right) \de z + \left(\dot z - y \dot x\right) \de \lambda \right] \mathrm{d}t.
\end{split}
\end{equation}
Demanding that $\de \tilde{S}=0$ for all variations $\de x$, $\de y$, $\de z$, and $\de \lambda$ (i.e. applying Hamilton's principle) yields the equations of motion
\begin{equation}
\begin{split}
-\ddot x+\dot \lambda y + \lambda \dot y &= 0, \\
\ddot y + \lambda \dot x &= 0, \\
\ddot z+\dot \lambda &= 0, \\
\dot z - y \dot x &= 0.
\end{split}
\end{equation}
The equation $\ddot z+\dot \lambda = 0$ implies that $\lambda = -\dot z + \tilde{c}$, for an integration constant $\tilde{c}$. Taking $\tilde{c} = 0$, $\lambda = -\dot z$. Making the substitutions $\dot \lambda = -\ddot z$ and $\lambda = -\dot z$ eliminates $\lambda$ from the equations of motion, simplifying them to
\begin{equation}
\begin{split}
-\ddot x-\ddot z y - \dot z \dot y &= 0, \\
\ddot y - \dot z \dot x &= 0, \\
\dot z - y \dot x &= 0.
\end{split}
\end{equation}
The equation $\dot z - y \dot x = 0$ implies that $\ddot z = \dot y \dot x + y \ddot x$. Making the substitutions $\dot z = y \dot x$ and $\ddot z = \dot y \dot x + y \ddot x$ eliminates $z$ from the equations of motion, simplifying them to 
\begin{equation}
\begin{split}
-\ddot x-\left(\dot y \dot x + y \ddot x \right) y - y \dot x \dot y &= 0, \\
\ddot y - y \dot x \dot x &= 0, \\
\end{split}
\end{equation}
which is equivalent to
\begin{equation}
\begin{split}
\ddot x + \frac{2y}{1+y^2}\dot x \dot y &= 0, \\
\ddot y - y {\dot x}^2 &= 0. \\
\end{split}
\end{equation}

For this simple problem, a more direct derivation of the equations of motion substitutes the constraint $\dot z = y \dot x$ into the Lagrangian. With this substitution, the Lagrangian becomes
\begin{equation}
l = \frac{1}{2}\left({\dot x}^2 +{\dot y}^2 + {\dot z}^2\right) = \frac{1}{2}\left({\dot x}^2 +{\dot y}^2 + y^2 {\dot x}^2\right) 
= \frac{1}{2}\left({\left(1+y^2 \right) \dot x}^2 +{\dot y}^2 \right) .
\end{equation}
The action integral is 
\begin{equation}
S = \int_a^b l \mathrm{d}t = \int_a^b  \frac{1}{2}\left({\left(1+y^2 \right) \dot x}^2 +{\dot y}^2 \right) \mathrm{d}t.
\end{equation}
The variation of the action integral is 
\begin{equation}
\begin{split}
\de S = \de \int_a^b l \mathrm{d}t = \int_a^b \de l \mathrm{d}t &= \int_a^b \left[\left(1+y^2 \right) \dot x \de \dot x + y \de y {\dot x}^2 + \dot y \de \dot y \right] \mathrm{d}t \\
&= \int_a^b \left[- \frac{\mathrm{d}}{\mathrm{d}t} \left[ \left(1+y^2 \right) \dot x \right] \de x + y {\dot x}^2 \de y - \ddot y \de y \right] \mathrm{d}t \\
&= \int_a^b \left[ -\left(\left(1+y^2 \right)\ddot x+2y \dot y \dot x \right) \de x + \left(y {\dot x}^2-\ddot y \right) \de y \right] \mathrm{d}t.
\end{split}
\end{equation}
Demanding that $\de S = 0$ for all variations $\de x$ and $\de y$ (i.e. applying Hamilton's principle) recovers the previously obtained equations of motion
\begin{equation}
\begin{split}
\ddot x + \frac{2y}{1+y^2}\dot x \dot y &=0, \\
\ddot y - y {\dot x}^2 &= 0.
\end{split}
\end{equation}

\paragraph{Comparison}
Figure~\ref{fig_LdA_vak} compares the solutions of the Lagrange-d'Alembert and vakonomic equations of motion during the time interval $0 \le t \le 1$ assuming that the initial conditions at time $t=0$ are $x(0)=1$, $\dot x(0)=1$, $y(0)=1$, and $\dot y(0)=1$. The \mcode{MATLAB} routine \mcode{ode45} was used to numerically solve the vakonomic equations of motion using the default error tolerances. Figures~\ref{fig_x_LdA_vak} and \ref{fig_y_LdA_vak} show that the $x$- and $y$-components of the solutions of the Lagrange-d'Alembert and vakonomic equations of motion disagree. Figure~\ref{fig_cony_LdA_vak} shows that the Lagrange-d'Alembert solution conserves the $y$-component of momentum, while the vakonomic solution does not. 

\begin{figure}[h]
	\centering
	\subfloat[The $x$-component of the Lagrange-d'Alembert and vakonomic solutions disagree.]{\includegraphics[scale=.33]{x_LdA_vs_vak_eps}\label{fig_x_LdA_vak}}
	\hspace{5mm}
	\subfloat[The $y$-component of the Lagrange-d'Alembert and vakonomic solutions disagree.]{\includegraphics[scale=0.33]{y_LdA_vs_vak_eps}\label{fig_y_LdA_vak}}
	\hspace{5mm}
	\subfloat[The Lagrange-d'Alembert solution conserves the $y$-component of momentum, while the vakonomic solution does not.]{\includegraphics[scale=0.33]{y_momentum_eps}\label{fig_cony_LdA_vak}}
	\caption{Lagrange-d'Alembert vs vakonomic solutions for a simple nonholonomic particle.}
	\label{fig_LdA_vak}
\end{figure}

While the vakonomic approach does not give the correct equations of motion for a mechanical system with nonholonomic constraints, the vakonomic approach does give the correct equations of motion for a mechanical system with holonomic constraints. The vakonomic approach is also applicable in optimal control, where a performance index must be minimized subject to satisfying equations of motion, which may have been obtained through the vakonomic or Lagrange-d'Alembert approaches depending on whether the mechanical constraints are holonomic or nonholonomic. When the vakonomic approach is applied to an optimal control problem as part of Pontryagin's Minimum Principle (discussed in Section~\ref{sec_optimal_control}), the resulting equations are called the controlled equations of motion.
} %%%END REM 

\section{Details for Deriving the Equations of Motion for the Rolling Ball} \label{app_rball}
By defining $\mathbf{s}_i \equiv r \bGamma +\bchi_i$ for $0\le i\le n$ and combining the summations, the equations of motion \eqref{uncon_ball_eqns_start} for the rolling ball become
\begin{equation}
\begin{split} \label{uncon_ball_eqns1}
\mathbf{0} &= \inertia \dot  \bOm+\bOm \times \inertia \bOm-r \bGamma \times \tilde \bGamma+ \sum_{i=0}^n m_i \Big\{-g \bGamma \times \bchi_i+ \left[r \bGamma \times \bOm+{\dot \bchi}_i \right] \times \left[\bOm \times \mathbf{s}_i +{\dot \bchi}_i \right] \\
& \hphantom{=} + \mathbf{s}_i \times \left\{ \dot \bOm \times \mathbf{s}_i+ \bOm \times \left[r \bGamma \times \bOm+{\dot \bchi}_i \right] +{\ddot \bchi}_i \right\}  \\
& \hphantom{=}+ \bOm \times \left\{ \mathbf{s}_i \times \left[\bOm \times \mathbf{s}_i +{\dot \bchi}_i \right] \right\}
- \left[r \bGamma \times \bOm \right] \times \left[\bOm \times \bchi_i+{\dot \bchi}_i \right] \Big\}.
\end{split}
\end{equation}
Since $\mathbf{s}_i \equiv r \bGamma +\bchi_i$ for $0\le i\le n$,
\begin{equation} \label{eq_imr4}
-g \bGamma \times \bchi_i = \bchi_i \times \left(g \bGamma \right)  = \left(\mathbf{s}_i - r \bGamma \right) \times \left(g \bGamma \right) = \mathbf{s}_i \times \left(g \bGamma \right)
\end{equation}
and
\begin{equation} \label{eq_imr_1}
\begin{split}
\left[r \bGamma \times \bOm+{\dot \bchi}_i \right] \times \left[\bOm \times \mathbf{s}_i +{\dot \bchi}_i \right] \\
-\left[r \bGamma \times \bOm \right] \times \left[\bOm \times \bchi_i+{\dot \bchi}_i \right]
&= \left[r \bGamma \times \bOm+{\dot \bchi}_i \right] \times \left[\bOm \times \mathbf{s}_i +{\dot \bchi}_i \right] \\
&\hphantom{=} -\left[r \bGamma \times \bOm \right] \times \left[\bOm \times \left\{ \mathbf{s}_i-r \bGamma \right\}+{\dot \bchi}_i \right] \\
&= \left[r \bGamma \times \bOm+{\dot \bchi}_i \right] \times \left[\bOm \times \mathbf{s}_i +{\dot \bchi}_i \right]-
\left[r \bGamma \times \bOm \right] \times \left[\bOm \times \mathbf{s}_i+{\dot \bchi}_i \right] \\
&= {\dot \bchi}_i \times \left[\bOm \times \mathbf{s}_i \right]. \\
\end{split}
\end{equation}
Moreover, by exploiting Jacobi's identity for the sum of permuted triple cross products, we find 
\begin{equation} \label{eq_imr_2}
\begin{split}
{\dot \bchi}_i \times \left[\bOm \times \mathbf{s}_i \right]+\mathbf{s}_i \times \left[\bOm \times {\dot \bchi}_i  \right]+\bOm \times \left[\mathbf{s}_i \times {\dot \bchi}_i  \right] = 2 \mathbf{s}_i \times \left[\bOm \times {\dot \bchi}_i  \right].
\end{split}
\end{equation}
By using \eqref{eq_imr4}, \eqref{eq_imr_1}, and \eqref{eq_imr_2} in \eqref{uncon_ball_eqns1}, the equations of motion \eqref{uncon_ball_eqns1} simplify to
\begin{equation}
\begin{split} \label{uncon_ball_eqns2}
\mathbf{0} &= \inertia \dot  \bOm+\bOm \times \inertia \bOm-r \bGamma \times \tilde \bGamma+ \sum_{i=0}^n m_i \Big\{\mathbf{s}_i \times \left(g \bGamma \right)+ 2 \mathbf{s}_i \times \left[\bOm \times {\dot \bchi}_i  \right] \\
& \hphantom{=} + \mathbf{s}_i \times \left\{ \dot \bOm \times \mathbf{s}_i+ \bOm \times \left[r \bGamma \times \bOm \right] +{\ddot \bchi}_i \right\} + \bOm \times \left\{ \mathbf{s}_i \times \left[\bOm \times \mathbf{s}_i \right] \right\} \Big\}.
\end{split}
\end{equation}
\rem{Interestingly, these equations can be simplified further using expressions for quadruple cross products. 
Indeed, given arbitrary vectors $\mathbf{a},\mathbf{b} \in \mathbb{R}^3$, the following identity holds: 
$$\mathbf{a} \times \left\{ \mathbf{b} \times \left[\mathbf{a} \times \mathbf{b} \right] \right\}=\mathbf{a} \times \left\{ \left(\mathbf{b} \cdot \mathbf{b} \right)\mathbf{a} -\left(\mathbf{a} \cdot \mathbf{b} \right)\mathbf{b} \right\}=-\left(\mathbf{a} \cdot \mathbf{b} \right) \mathbf{a} \times \mathbf{b}. $$
This gives the identity
\begin{equation} \label{eq_cross_iden}
\mathbf{a} \times \left\{ \mathbf{b} \times \left[\mathbf{a} \times \mathbf{b} \right] \right\}+\mathbf{b} \times \left\{ \mathbf{a} \times \left[\mathbf{b} \times \mathbf{a} \right] \right\} =-\left(\mathbf{a} \cdot \mathbf{b} \right) \mathbf{a} \times \mathbf{b}-\left(\mathbf{a} \cdot \mathbf{b} \right) \mathbf{b} \times \mathbf{a}= \mathbf{0}.
\end{equation}}
For arbitrary vectors $\mathbf{a},\mathbf{b} \in \mathbb{R}^3$,  Jacobi's identity  yields the following identity for the  sum of quadruple cross products: 
\begin{equation} \label{eq_cross_iden}
\mathbf{a} \times \left\{ \mathbf{b} \times \left[\mathbf{a} \times \mathbf{b} \right] \right\}+\mathbf{b} \times \left\{ \mathbf{a} \times \left[\mathbf{b} \times \mathbf{a} \right] \right\} = \mathbf{0}.
\end{equation}}
Since $\mathbf{s}_i \equiv r \bGamma +\bchi_i$ for $0\le i\le n$ and using the identity \eqref{eq_cross_iden}, it follows that 
\begin{equation} \label{eq_imr_3}
\begin{split}
\mathbf{s}_i \times \left\{ \bOm \times \left[r \bGamma \times \bOm \right] \right\}
+ \bOm \times \left\{ \mathbf{s}_i \times \left[\bOm \times \mathbf{s}_i \right] \right\} &= \mathbf{s}_i \times \left\{ \bOm \times \left[ \left(\mathbf{s}_i-\bchi_i \right) \times \bOm \right] \right\} 
+ \bOm \times \left\{ \mathbf{s}_i \times \left[\bOm \times \mathbf{s}_i \right] \right\} \\
&= \mathbf{s}_i \times \left\{ \bOm \times \left[ \mathbf{s}_i \times \bOm \right] \right\} - \mathbf{s}_i \times \left\{ \bOm \times \left[ \bchi_i \times \bOm \right] \right\} \\
&\hphantom{=} + \bOm \times \left\{ \mathbf{s}_i \times \left[\bOm \times \mathbf{s}_i \right] \right\} \\
&= - \mathbf{s}_i \times \left\{ \bOm \times \left[ \bchi_i \times \bOm \right] \right\}.
\end{split}
\end{equation}
Using \eqref{eq_imr_3}, the equations of motion \eqref{uncon_ball_eqns2} simplify to
\begin{equation}
\label{uncon_ball_eqns5}
\inertia \dot  \bOm+\bOm \times \inertia \bOm+r \tilde \bGamma \times \bGamma+ \sum_{i=0}^n m_i \mathbf{s}_i \times  \left\{ g \bGamma +\dot \bOm \times \mathbf{s}_i+ \bOm \times \left(\bOm \times \bchi_i +2 {\dot \bchi}_i \right) +{\ddot \bchi}_i \right\} = \mathbf{0} . 
\end{equation} 
\rem{\begin{equation}
\begin{split} \label{uncon_ball_eqns3}
\mathbf{0} &= \inertia \dot  \bOm+\bOm \times \inertia \bOm-r \bGamma \times \tilde \bGamma+ \sum_{i=0}^n m_i \Big\{-g \bGamma \times \bchi_i+ 2 \mathbf{s}_i \times \left[\bOm \times {\dot \bchi}_i  \right] \\
& \hphantom{=} + \mathbf{s}_i \times \left\{ \dot \bOm \times \mathbf{s}_i- \bOm \times \left[ \bchi_i \times \bOm \right] +{\ddot \bchi}_i \right\} \Big\}.
\end{split}
\end{equation}
Since  $\mathbf{s}_i \equiv r \bGamma +\bchi_i$ for $0\le i\le n$, the gravity terms in \eqref{uncon_ball_eqns3} can be reduced further as follows: 
\begin{equation} \label{eq_imr4}
-g \bGamma \times \bchi_i = \bchi_i \times \left(g \bGamma \right)  = \left(\mathbf{s}_i - r \bGamma \right) \times \left(g \bGamma \right) = \mathbf{s}_i \times \left(g \bGamma \right).
\end{equation}
Using \eqref{eq_imr4} in \eqref{uncon_ball_eqns3}, the equations of motion \eqref{uncon_ball_eqns3} become 
\rem{ %%%%BEGIN REM 
	\begin{equation}
	\begin{split} \label{uncon_ball_eqns4}
	\mathbf{0} &= \inertia \dot  \bOm+\bOm \times \inertia \bOm-r \bGamma \times \tilde \bGamma+ \sum_{i=0}^n m_i \Big\{\mathbf{s}_i \times \left(g \bGamma \right)+ 2 \mathbf{s}_i \times \left[\bOm \times {\dot \bchi}_i  \right] \\
	& \hphantom{=} + \mathbf{s}_i \times \left\{ \dot \bOm \times \mathbf{s}_i- \bOm \times \left[ \bchi_i \times \bOm \right] +{\ddot \bchi}_i \right\} \Big\}.
	\end{split}
	\end{equation}
	Combining terms and eliminating minus signs by re-ordering cross product pairings in \eqref{uncon_ball_eqns4}, the  equations of motion simplify further to
} %%%END REM 
\begin{equation}
\label{uncon_ball_eqns5}
\inertia \dot  \bOm+\bOm \times \inertia \bOm+r \tilde \bGamma \times \bGamma+ \sum_{i=0}^n m_i \mathbf{s}_i \times  \left\{ g \bGamma +\dot \bOm \times \mathbf{s}_i+ \bOm \times \left(\bOm \times \bchi_i +2 {\dot \bchi}_i \right) +{\ddot \bchi}_i \right\} = \mathbf{0} . 
\end{equation} }
\rem{We can express \eqref{uncon_ball_eqns5} as an explicit system of ODEs for $\bOm$ as follows. First, note that}
Finally, since 
\begin{equation} \label{eq_imr_5}
\mathbf{s}_i \times \left\{ \dot \bOm \times \mathbf{s}_i \right\} = -\mathbf{s}_i \times \left\{\mathbf{s}_i \times \dot \bOm  \right\}= -\mathbf{s}_i \times \left\{\widehat{\mathbf{s}_i} \dot \bOm  \right\}=  -\widehat{\mathbf{s}_i} \widehat{\mathbf{s}_i} \dot \bOm = -  \widehat{\mathbf{s}_i}^2 \dot \bOm,
\end{equation}
where for $\mathbf{v} = \begin{bmatrix} v_1 & v_2 & v_3  \end{bmatrix}^\mathsf{T} \in \mathbb{R}^3$, $\widehat{\mathbf{v}}^2=\widehat{\mathbf{v}}\widehat{\mathbf{v}}$ is the symmetric matrix given by
\begin{equation}
\widehat{\mathbf{v}}^2 = \begin{bmatrix}
-(v_2^2+v_3^2) & v_1 v_2 & v_1 v_3 \\
v_1 v_2 & -(v_1^2+v_3^2)  & v_2 v_3 \\
v_1 v_3 & v_2 v_3 & -(v_1^2+v_2^2) 
\end{bmatrix},
\end{equation}
we can solve explicitly for $\dot \bOm$ in \eqref{uncon_ball_eqns5} to obtain the equations of motion \eqref{uncon_ball_eqns_explicit} for the rolling ball.

\section{Details for Deriving the Equation of Motion for the Rolling Disk} \label{app_rdisk}
\rem{Recall the equations of motion for the rolling disk given in \eqref{eqmo_chap_disk_1}:
\begin{equation} \label{eqmo_chap_disk_1_rep}
- \ddot \phi  \mathbf{e}_2 =  \left[\sum_{i=0}^n m_i \widehat{\mathbf{s}_i}^2  -\inertia \right]^{-1}  \left[ -r F_{\mathrm{e},1} \mathbf{e}_2+ \sum_{i=0}^n m_i K_i \mathbf{e}_2 \right]
= \left( -r F_{\mathrm{e},1}+ \sum_{i=0}^n m_i K_i \right) \left[\sum_{i=0}^n m_i \widehat{\mathbf{s}_i}^2  -\inertia \right]^{-1} \mathbf{e}_2,
\end{equation}
where $\widehat{\mathbf{s}_i}^2$ is given in \eqref{eq_s_i_squared}.
Note that $\left[\sum_{i=0}^n m_i \widehat{\mathbf{s}_i}^2  -\inertia \right]^{-1} \mathbf{e}_2$ is just the middle column of the matrix inverse of $A = \sum_{i=0}^n m_i \widehat{\mathbf{s}_i}^2  -\inertia$. Denote the entries of $A$ by
\begin{equation}
A = \sum_{i=0}^n m_i \widehat{\mathbf{s}_i}^2  -\inertia = \begin{bmatrix} a_{11} & a_{12} & a_{13} \\ a_{21} & a_{22} & a_{23} \\ a_{31} & a_{32} & a_{33} \end{bmatrix}.
\end{equation}
Since $\inertia$ is diagonal and from \eqref{eq_s_i_squared}, $a_{12} = a_{21} = a_{23} = a_{32} = 0$, so that
\begin{equation}
A = \sum_{i=0}^n m_i \widehat{\mathbf{s}_i}^2  -\inertia = \begin{bmatrix} a_{11} & 0 & a_{13} \\ 0 & a_{22} & 0 \\ a_{31} & 0 & a_{33} \end{bmatrix}
\end{equation}
and the determinant of $A$ simplifies to 
\begin{equation}
\begin{split}
\det A &= a_{11} a_{22} a_{33}+a_{21} a_{32} a_{13}+a_{31} a_{12} a_{23}-a_{11} a_{32} a_{23}-a_{31} a_{22} a_{13}-a_{21} a_{12} a_{33} \\
&= a_{11} a_{22} a_{33} - a_{31} a_{22} a_{13}
= a_{22} \left( a_{11} a_{33} - a_{31} a_{13} \right).
\end{split}
\end{equation}
Using the formula for the inverse of a $3 \times 3$ matrix, the middle column of the matrix inverse of $ A = \sum_{i=0}^n m_i \widehat{\mathbf{s}_i}^2  -\inertia$ is
\begin{equation} \label{eq_mid_col_inverse}
\begin{split}
\left[\sum_{i=0}^n m_i \widehat{\mathbf{s}_i}^2  -\inertia \right]^{-1} \mathbf{e}_2 &= A^{-1} \mathbf{e}_2
= \begin{bmatrix} a_{11} & 0 & a_{13} \\ 0 & a_{22} & 0 \\ a_{31} & 0 & a_{33} \end{bmatrix}^{-1} \mathbf{e}_2 
= \frac{1}{\det A} \begin{bmatrix} a_{13} a_{32} - a_{12} a_{33} \\ a_{11} a_{33} - a_{13} a_{31} \\ a_{12} a_{31} - a_{11} a_{32}  \end{bmatrix} \\
&= \frac{1}{a_{22} \left( a_{11} a_{33} - a_{31} a_{13} \right)} \begin{bmatrix} 0 \\ a_{11} a_{33} - a_{13} a_{31} \\ 0 \end{bmatrix} 
=  \frac{1}{a_{22}} \begin{bmatrix} 0 \\ 1 \\ 0 \end{bmatrix} 
= \frac{1}{a_{22}} \mathbf{e}_2.
\end{split} 
\end{equation}
Plugging \eqref{eq_mid_col_inverse} into \eqref{eqmo_chap_disk_1_rep}, the equations of motion simplify to
\begin{equation} \label{eqmo_chap_disk_2}
- \ddot \phi  \mathbf{e}_2 = \left( -r F_{\mathrm{e},1}+ \sum_{i=0}^n m_i K_i \right) \left[\sum_{i=0}^n m_i \widehat{\mathbf{s}_i}^2  -\inertia \right]^{-1} \mathbf{e}_2 = \frac{1}{a_{22}} \left( -r F_{\mathrm{e},1}+ \sum_{i=0}^n m_i K_i \right) \mathbf{e}_2,
\end{equation}
which gives the scalar equation of motion
\begin{equation} \label{eqmo_chap_disk_3}
\ddot \phi = \frac{-1}{a_{22}} \left( -r F_{\mathrm{e},1}+ \sum_{i=0}^n m_i K_i \right).
\end{equation}
From \eqref{eq_s_i_squared},
\begin{equation} \label{eq_a_22}
a_{22} = \sum_{i=0}^n \left\{ m_i \left[-\left( r \sin \phi + \zeta_{i,1} \right)^2-\left( r \cos \phi+ \zeta_{i,3} \right)^2 \right] \right\}  -d_2.
\end{equation}
Plugging \eqref{eq_a_22} into \eqref{eqmo_chap_disk_3} gives the equation of motion for the rolling disk \eqref{eqmo_chap_disk_4}.}
Note that $\left[\sum_{i=0}^n m_i \widehat{\mathbf{s}_i}^2  -\inertia \right]^{-1} \mathbf{e}_2$ is just the middle column of the matrix inverse of $A = \sum_{i=0}^n m_i \widehat{\mathbf{s}_i}^2  -\inertia$, where $\widehat{\mathbf{s}_i}^2$ is given in \eqref{eq_s_i_squared}. Denote the entries of $A$ by
\begin{equation}
A = \sum_{i=0}^n m_i \widehat{\mathbf{s}_i}^2  -\inertia = \begin{bmatrix} a_{11} & a_{12} & a_{13} \\ a_{21} & a_{22} & a_{23} \\ a_{31} & a_{32} & a_{33} \end{bmatrix}.
\end{equation}
Since $\inertia$ is diagonal and from \eqref{eq_s_i_squared}, $a_{12} = a_{21} = a_{23} = a_{32} = 0$, so that
\begin{equation}
A = \sum_{i=0}^n m_i \widehat{\mathbf{s}_i}^2  -\inertia = \begin{bmatrix} a_{11} & 0 & a_{13} \\ 0 & a_{22} & 0 \\ a_{31} & 0 & a_{33} \end{bmatrix}
\end{equation}
and the determinant of $A$ simplifies to 
\begin{equation}
\begin{split}
\det A &= a_{11} a_{22} a_{33}+a_{21} a_{32} a_{13}+a_{31} a_{12} a_{23}-a_{11} a_{32} a_{23}-a_{31} a_{22} a_{13}-a_{21} a_{12} a_{33} \\
&= a_{11} a_{22} a_{33} - a_{31} a_{22} a_{13}
= a_{22} \left( a_{11} a_{33} - a_{31} a_{13} \right).
\end{split}
\end{equation}
From \eqref{eq_s_i_squared},
\begin{equation} \label{eq_a_22}
a_{22} = \sum_{i=0}^n \left\{ m_i \left[-\left( r \sin \phi + \zeta_{i,1} \right)^2-\left( r \cos \phi+ \zeta_{i,3} \right)^2 \right] \right\}  -d_2.
\end{equation}
Using the formula for the inverse of a $3 \times 3$ matrix and \eqref{eq_a_22}, the middle column of the matrix inverse of $ A = \sum_{i=0}^n m_i \widehat{\mathbf{s}_i}^2  -\inertia$ is
\begin{equation} \label{eq_mid_col_inverse}
\begin{split}
\left[\sum_{i=0}^n m_i \widehat{\mathbf{s}_i}^2  -\inertia \right]^{-1} \mathbf{e}_2 &= A^{-1} \mathbf{e}_2
= \begin{bmatrix} a_{11} & 0 & a_{13} \\ 0 & a_{22} & 0 \\ a_{31} & 0 & a_{33} \end{bmatrix}^{-1} \mathbf{e}_2 
= \frac{1}{\det A} \begin{bmatrix} a_{13} a_{32} - a_{12} a_{33} \\ a_{11} a_{33} - a_{13} a_{31} \\ a_{12} a_{31} - a_{11} a_{32}  \end{bmatrix} \\
&= \frac{1}{a_{22} \left( a_{11} a_{33} - a_{31} a_{13} \right)} \begin{bmatrix} 0 \\ a_{11} a_{33} - a_{13} a_{31} \\ 0 \end{bmatrix} 
=  \frac{1}{a_{22}} \begin{bmatrix} 0 \\ 1 \\ 0 \end{bmatrix} 
= \frac{1}{a_{22}} \mathbf{e}_2 \\
&= \frac{1}{\sum_{i=0}^n \left\{ m_i \left[-\left( r \sin \phi + \zeta_{i,1} \right)^2-\left( r \cos \phi+ \zeta_{i,3} \right)^2 \right] \right\}  -d_2} \mathbf{e}_2.
\end{split} 
\end{equation}
Plugging \eqref{eq_mid_col_inverse} into \eqref{eqmo_chap_disk_1} gives the scalar equation of motion \eqref{eqmo_chap_disk_4} for the rolling disk.

\section{Quaternions} \label{app_quaternions}

Quaternions were invented by William Rowan Hamilton in 1843. Good references on quaternions and how they are used to model rigid body dynamics are \cite{Ho2011_pII,graf2008quaternions,stevens2015aircraft,baraff2001physically}. The set of quaternions, which is isomorphic to $\mathbb{R}^4$, is denoted by $\mathbb{H}$. A quaternion $\mathfrak{p} \in \mathbb{H}$ can be expressed as the column vector 
\begin{equation}
\mathfrak{p} =
\begin{bmatrix} p_0 & p_1 & p_2 & p_3 \end{bmatrix}^\mathsf{T}.
\end{equation}
Given a column vector $\bv \in \mathbb{R}^3$, $\bv^\sharp$ is the quaternion 
\begin{equation}
\bv^\sharp = \begin{bmatrix}0 \\ \bv \end{bmatrix}.
\end{equation}
Given a quaternion $\mathfrak{p} \in \mathbb{H}$, $\mathfrak{p}^\flat \in \mathbb{R}^3$ is the column vector such that 
\begin{equation}
\mathfrak{p} = \begin{bmatrix} p_0 \\ \mathfrak{p}^\flat \end{bmatrix}.
\end{equation}
Given a column vector $\bv \in \mathbb{R}^3$, note that 
\begin{equation}
\left(\bv^\sharp\right)^\flat=\bv.
\end{equation}
However, given a quaternion $\mathfrak{p} \in \mathbb{H}$, 
\begin{equation}
\left( \mathfrak{p}^\flat \right)^\sharp = \mathfrak{p} \quad \mathrm{iff} \quad \mathfrak{p} = \begin{bmatrix} 0 \\ \mathfrak{p}^\flat \end{bmatrix}. 
\end{equation}
Given quaternions $\mathfrak{p}=\begin{bmatrix}p_0 \\ \mathfrak{p}^\flat \end{bmatrix},\mathfrak{q}=\begin{bmatrix}q_0 \\ \mathfrak{q}^\flat \end{bmatrix} \in \mathbb{H}$, 
their sum is
\begin{equation}
\mathfrak{p} + \mathfrak{q}=\begin{bmatrix}p_0 \\ \mathfrak{p}^\flat \end{bmatrix}+ \begin{bmatrix}q_0 \\ \mathfrak{q}^\flat \end{bmatrix}=\begin{bmatrix}p_0+q_0 \\ \mathfrak{p}^\flat + \mathfrak{q}^\flat \end{bmatrix},
\end{equation}
their product is
\begin{equation}
\mathfrak{p} \mathfrak{q}=\begin{bmatrix} p_0 \\ \mathfrak{p}^\flat \end{bmatrix} \begin{bmatrix}q_0 \\ \mathfrak{q}^\flat \end{bmatrix}=\begin{bmatrix}p_0q_0-\mathfrak{p}^\flat \cdot \mathfrak{q}^\flat \\ p_0\mathfrak{q}^\flat + q_0\mathfrak{p}^\flat+\mathfrak{p}^\flat \times \mathfrak{q}^\flat \end{bmatrix},
\end{equation}
and their dot product is 
\begin{equation}
\mathfrak{p} \cdot \mathfrak{q}=\begin{bmatrix} p_0 \\ \mathfrak{p}^\flat \end{bmatrix} \cdot \begin{bmatrix}q_0 \\ \mathfrak{q}^\flat \end{bmatrix}=\begin{bmatrix}p_0 & p_1 & p_2 & p_3\end{bmatrix}^\mathsf{T} \cdot \begin{bmatrix}q_0 & q_1 & q_2 & q_3\end{bmatrix}^\mathsf{T}= p_0 q_0 + \mathfrak{p}^\flat \cdot \mathfrak{q}^\flat=p_0 q_0 +p_1 q_1 + p_2 q_2 + p_3 q_3.
\end{equation}
It may be shown that multiplication in $\mathbb{H}$ is associative (i.e. $\mathfrak{p} \left(\mathfrak{q} \mathfrak{r} \right) = \left( \mathfrak{p} \mathfrak{q} \right) \mathfrak{r} \quad \forall \mathfrak{p},\mathfrak{q},\mathfrak{r} \in \mathbb{H} $) but not commutative (i.e. $\mathfrak{p} \mathfrak{q} \ne \mathfrak{q} \mathfrak{p} $ for general $\mathfrak{p},\mathfrak{q} \in \mathbb{H}$). Given $c \in \mathbb{R}$ and a quaternion $\mathfrak{p} = \begin{bmatrix}p_0 \\ \mathfrak{p}^\flat \end{bmatrix} \in \mathbb{H}$, scalar multiplication of $\mathfrak{p}$ by $c$ is
\begin{equation}
c \mathfrak{p} = c \begin{bmatrix}p_0 \\ \mathfrak{p}^\flat \end{bmatrix} = \begin{bmatrix}c p_0 \\ c\mathfrak{p}^\flat \end{bmatrix}.
\end{equation}
Given a quaternion $\mathfrak{p} = \begin{bmatrix}p_0 \\ \mathfrak{p}^\flat \end{bmatrix} \in \mathbb{H}$, its conjugate is
\begin{equation}
\mathfrak{p}^* = \begin{bmatrix} p_0 \\ \unaryminus\mathfrak{p}^\flat \end{bmatrix},
\end{equation}
its magnitude is
\begin{equation}
\left| \mathfrak{p} \right| = \left( \mathfrak{p} \cdot \mathfrak{p} \right)^\frac{1}{2} = \left( p_0^2 +  \mathfrak{p}^\flat \cdot \mathfrak{p}^\flat \right)^\frac{1}{2} ,
\end{equation}
and its inverse is
\begin{equation}
\mathfrak{p}^{-1} = \frac{\mathfrak{p}^*}{\left| \mathfrak{p} \right|^2}.
\end{equation}
In the language of abstract algebra, $\mathbb{H}$ is a four-dimensional associative normed division algebra over the real numbers.
$\mathscr{S} \subset \mathbb{H}$ denotes the set of unit quaternions, also called versors, which is isomorphic to $\mathbb{S}^3 \subset \mathbb{R}^4$. That is,
\begin{equation}
\mathscr{S} \equiv \left\{ \mathfrak{q} = \begin{bmatrix}q_0 & q_1 & q_2 & q_3 \end{bmatrix}^\mathsf{T} \in \mathbb{R}^4 : \left| \mathfrak{q} \right|^2 = \mathfrak{q} \cdot \mathfrak{q} = q_0^2+q_1^2+q_2^2+q_3^2=1  \right\} \subset \mathbb{H}.
\end{equation}
The set of versors $\mathscr{S}$ is useful because it may be utilized to parameterize the set of rotation matrices $SO(3)$. Given a versor 
\begin{equation}
\mathfrak{q} = \begin{bmatrix}q_0 & q_1 & q_2 & q_3 \end{bmatrix}^\mathsf{T} \in \mathscr{S},
\end{equation}
the corresponding rotation matrix $\Lambda \in SO(3)$ is
\begin{equation} \label{eq_versor_rot}
\Lambda = \begin{bmatrix} 1-2\left(q_2^2+q_3^2\right) & 2\left(q_1q_2-q_0q_3 \right) & 2\left(q_1 q_3+q_0q_2 \right) \\ 2\left(q_1 q_2+q_0q_3 \right) & 1-2\left(q_1^2+q_3^2\right) &  2\left(q_2 q_3-q_0q_1 \right) \\ 2\left(q_1 q_3-q_0q_2 \right) & 2\left(q_2 q_3+q_0q_1 \right) & 1-2\left(q_1^2+q_2^2\right) \end{bmatrix} \in SO(3).
\end{equation}
It is easy to see from \eqref{eq_versor_rot}, that the versors 
\begin{equation}
\mathfrak{q} = \begin{bmatrix}q_0 & q_1 & q_2 & q_3 \end{bmatrix}^\mathsf{T} \in \mathscr{S} \quad \mathrm{and} \quad \unaryminus \mathfrak{q} = \begin{bmatrix}\unaryminus q_0 & \unaryminus q_1 & \unaryminus q_2 & \unaryminus q_3 \end{bmatrix}^\mathsf{T} \in \mathscr{S} 
\end{equation}
correspond to the same rotation matrix $\Lambda \in SO(3)$, so that $\mathscr{S}$ is a double covering of $SO(3)$. Given a vector $\bY \in \mathbb{R}^3$, the rotation of $\bY$ by $\Lambda \in SO(3)$ can be realized using the versor $\mathfrak{q} \in \mathscr{S}$ via the Euler-Rodrigues formula
\begin{equation} \label{eq_euler_rod}
\Lambda \bY = \left[\mathfrak{q} \bY^\sharp \mathfrak{q}^{-1} \right]^\flat.
\end{equation}
Since $\mathfrak{q}^{-1} \in \mathscr{S}$ parameterizes $\Lambda^{-1} \in SO(3)$, \eqref{eq_euler_rod} says that the rotation of $\mathbf{y} \in \mathbb{R}^3$ by $\Lambda^{-1} \in SO(3)$ can be realized using the versor $\mathfrak{q}^{-1} \in \mathscr{S}$ via
\begin{equation}
\Lambda^{-1} \mathbf{y} = \left[\mathfrak{q}^{-1} \mathbf{y}^\sharp \mathfrak{q} \right]^\flat.
\end{equation}
Now consider a rigid body, such as a free rigid body, a heavy top, Suslov's problem, a rolling disk, a rolling ball, etc., with orientation matrix $\Lambda \in SO(3)$ (i.e. $\Lambda$ maps the body frame into the spatial frame) and body angular velocity 
\begin{equation}
\bOm \equiv \left[ \Lambda^{-1} \dot \Lambda \right]^\vee =\left[ \Lambda^\mathsf{T} \dot \Lambda \right]^\vee \in \mathbb{R}^3,
\end{equation}
so that 
\begin{equation}
\dot \Lambda = \Lambda \widehat{\bOm}.
\end{equation}
Let $\mathfrak{q} \in \mathscr{S}$ denote a versor corresponding to $\Lambda$. Then it may be shown that
\begin{equation}
\dot {\mathfrak{q}} = \frac{1}{2} \mathfrak{q} \bOm^\sharp.
\end{equation}	

\end{document}